\documentclass[12pt,a4paper,english]{amsart}
    
\title{Cuspidal character sheaves on graded Lie algebras}
\author{Wille Liu}\address{Institute of Mathematics, Academia Sinica, 7F, Astronomy-Mathematics Building, No. 1, Sec. 4, Roosevelt Road, Taipei, Taiwan}\email{wliu@gate.sinica.edu.tw}

\author{Cheng-Chiang Tsai}\address{Institute of Mathematics, Academia Sinica, 6F, Astronomy-Mathematics Building, No. 1, Sec. 4, Roosevelt Road, Taipei, Taiwan,\vskip.05cm
\noindent also Department of Applied Mathematics, National Sun Yat-Sen University, and Department of Mathematics, National Taiwan University}\email{chchtsai@gate.sinica.edu.tw}
\thanks{CCT is supported by NSTC grants 113-2115-M-001-002 and 113-2628-M-001-012}

\author{Kari Vilonen}\address{School of Mathematics and Statistics, University of Melbourne, VIC 3010, Australia, also Department of Mathematics and Statistics, University of Helsinki, Helsinki, Finland}
\email{kari.vilonen@unimelb.edu.au, kari.vilonen@helsinki.fi}
\thanks{KV was supported in part by the ARC grants  FL200100141, DP250100824 and the Academy of Finland}
\author{Ting Xue}
\address{School of Mathematics and Statistics, University of Melbourne, VIC 3010, Australia, also Department of Mathematics and Statistics, University of Helsinki, Helsinki, Finland} 
\email{ting.xue@unimelb.edu.au}
\thanks{TX was supported in part by the ARC grant   DP250100824}
\date{}
    
\usepackage{etoolbox}
\usepackage{enumitem}
\usepackage{aliascnt}
\usepackage{amsmath}
\usepackage{amscd}
\usepackage{amssymb}
\usepackage{amsfonts}
\usepackage{amsthm}
\usepackage{bbm}
\usepackage{cancel}
\usepackage{color}
\usepackage{eucal}
\usepackage{pdfsync}
\usepackage[all,cmtip]{xy}
\usepackage{graphicx}
\usepackage{graphics}
\usepackage[hidelinks]{hyperref}
\usepackage{latexsym}
\usepackage{mathrsfs}
\usepackage{placeins}
\usepackage{bookmark}
\usepackage{mathtools}
\usepackage{stmaryrd}
\usepackage{url}
\usepackage{mathabx}
\usepackage{float}

\usepackage{xcolor}
\usepackage{tikz-cd}
\usepackage[normalem]{ulem}
\usepackage{manfnt}
\usepackage{leftindex}
\newcommand{\li}{\leftindex}

  \usepackage{mathrsfs}

  \usepackage{accents}

\usepackage[backend=biber,doi=false,url=false,isbn=false,url=false,clearlang,giveninits=true,style=alphabetic,maxbibnames=99,maxalphanames=6]{biblatex} 
\DeclareRedundantLanguages{english}{english}
\DeclareSourcemap{
  \maps{
    \map{
      \step[fieldsource=language, fieldset=langid, origfieldval, final]
      \step[fieldset=language, null]
    }
  }
}

\DeclareFieldFormat[article]{pages}{#1}
\DeclareFieldFormat[article]{page}{#1}
\DeclareFieldFormat[article]{volume}{\textbf{#1}}
\DeclareFieldFormat[article]{number}{no.~#1}

\DeclareFieldFormat[book]{volume}{\textbf{#1}}
\renewbibmacro{in:}{%
  \ifentrytype{article}{}{\printtext{\bibstring{in}\intitlepunct}}}

\renewbibmacro*{volume+number+eid}{%
  \printfield{volume}%
  \setunit*{\addcomma\space}%
  \printfield{number}%
  \setunit{\addcomma\space}%
  \printfield{eid}}
  
\numberwithin{equation}{section}

    \newcommand{\N}{\mathbb{N}}
    \newcommand{\Z}{\mathbb{Z}}
    \newcommand{\Q}{\mathbb{Q}}
    \newcommand{\R}{\mathbb{R}}
    \newcommand{\C}{\mathbb{C}}
    \newcommand{\D}{\operatorname{\mathbf{D}}}
    \newcommand{\Cc}{\mathbb{C}^{\times}}

    \newcommand{\cC}{\mathcal{C}}
    
    \newcommand{\cF}{\mathcal{F}}
    \newcommand{\cG}{\mathcal{G}}
    \newcommand{\cH}{\mathcal{H}}
    \newcommand{\cI}{\mathcal{I}}
    \newcommand{\cJ}{\mathcal{J}}
    \newcommand{\cK}{\mathcal{K}}
    \newcommand{\cL}{\mathcal{L}}
    \newcommand{\cM}{\mathcal{M}}
    \newcommand{\cP}{\mathcal{P}}
    
    \newcommand{\cS}{\mathcal{S}}
    \newcommand{\cT}{\mathcal{T}}
    \newcommand{\cU}{\mathcal{U}}

    \newcommand{\cZ}{\mathcal{Z}}
    \newcommand{\cO}{\mathcal{O}}

    \newcommand{\bV}{\mathbf{V}}
    \newcommand{\bX}{\mathbf{X}}
    
    \newcommand{\coker}{\operatorname{coker}}

    \usepackage{calligra}
    \DeclareMathOperator{\scrHom}{R\mathscr{H}\text{\kern -3pt {\calligra\large om}}\,}
    
    \newcommand{\Db}{\mathrm{D}^{\mathrm{b}}}
    \newcommand{\Dbn}{\mathrm{D}^{\mathrm{b},\nil}}
    \newcommand{\Dbc}{\mathrm{D}^{\mathrm{b}}_{\mathrm{c}}}
    
    \newcommand{\ul}{\underline}
    \newcommand{\rmH}{\mathrm{H}}
    \newcommand{\rmHc}{\operatorname{H}_{\mathrm{c}}}
    
    \newcommand{\Ga}{\mathbb{G}_a}
    \newcommand{\Gm}{\mathbb{G}_m}
    \newcommand{\til}{\tilde}
    
    \newcommand{\op}{\mathrm{op}}
    \newcommand{\nil}{\mathrm{nil}}
    
    \newcommand{\ra}{\rightarrow}

    \newcommand{\rc}{\mathrm{c}}
    
    \newcommand{\id}{\operatorname{id}}
    \newcommand{\Lie}{\operatorname{Lie}}
    
    \newcommand{\supp}{\operatorname{supp}}
    \newcommand{\codim}{\operatorname{codim}}

    \newcommand{\fc}{\mathfrak{c}}
    \newcommand{\fd}{\mathfrak{d}}
    \newcommand{\fg}{\mathfrak{g}}
    \newcommand{\fh}{\mathfrak{h}}
    \newcommand{\fl}{\mathfrak{l}}
    \newcommand{\fm}{\mathfrak{m}}
    
    \newcommand{\fp}{\mathfrak{p}}
    \newcommand{\fq}{\mathfrak{q}}
    \newcommand{\ft}{\mathfrak{t}}
    \newcommand{\fu}{\mathfrak{u}}
    \newcommand{\fv}{\mathfrak{v}}
      \newcommand{\fZ}{\mathfrak{Z}}
    
    \newcommand{\fsl}{\mathfrak{sl}}

    \newcommand{\Lm}{\mathfrak{m}}

    \newcommand{\CF}{\mathcal{F}}
    \newcommand{\CG}{\mathcal{G}}

    \newcommand{\CQ}{\mathcal{Q}}

    \newcommand{\CO}{\mathcal{O}}

    \newcommand{\Ad}{\operatorname{Ad}}
    
    \newcommand{\Aut}{\operatorname{Aut}}
    \newcommand{\Char}{\operatorname{Char}}
    \newcommand{\Cone}{\operatorname{Cone}}
    \newcommand{\End}{\operatorname{End}}
    
    \newcommand{\IC}{\operatorname{IC}}
    \newcommand{\pH}{\operatorname{\li^{\mathrm{p}}{\mathcal{H}}}}
    \newcommand{\Perv}{\operatorname{Perv}}
    \newcommand{\Res}{\operatorname{Res}}
    \newcommand{\res}{\operatorname{res}}
    \newcommand{\Ind}{\operatorname{Ind}}
    \newcommand{\NRes}{{\operatorname{Res}}}
    \newcommand{\Nres}{{\operatorname{res}}}
    \newcommand{\NInd}{{\operatorname{Ind}}}
    \newcommand{\Irr}{\operatorname{Irr}}
    \newcommand{\RGamma}{\operatorname{R\Gamma}}
    \newcommand{\RGammac}{\operatorname{R\Gamma_c}}
    \newcommand{\RHom}{\operatorname{RHom}}
    \newcommand{\dgHom}{\operatorname{dgHom}}
    \newcommand{\Hom}{\operatorname{Hom}}

    \newcommand{\mof}{\operatorname{-mod}}
    \newcommand{\dgmod}{\operatorname{-dgmod}}
    \newcommand{\dgMod}{\operatorname{-dgMod}}
    
    \newcommand{\ord}{\operatorname{ord}}
     \newcommand{\Cunip}{\operatorname{mon}}
 
     \newcommand{\FS}{\operatorname{FS}}

     \newcommand{\on}{\operatorname}

    \newcommand{\beqn}{\begin{equation*}}
\newcommand{\eeqn}{\end{equation*}}

\newcommand{\beq}{\begin{equation}}
\newcommand{\eeq}{\end{equation}}

\newcommand{\bega}{\begin{gathered}}
\newcommand{\eega}{\end{gathered}}

\newcommand{\bern}{\begin{eqnarray*}}
\newcommand{\eern}{\end{eqnarray*}}
\newcommand{\ber}{\begin{eqnarray}}
\newcommand{\eer}{\end{eqnarray}}

\newcommand{\Lz}{\underline{L}_{0}}
\newcommand{\Pz}{\underline{P}_{0}}
\newcommand{\flz}{\underline{\fl}_{0}}
\newcommand{\flo}{\underline{\fl}_{1}}

\newcommand{\fpo}{\underline{\fp}_{1}}

\newcommand{\fuo}{\underline{\fu}_{1}}

\theoremstyle{plain}
\newtheorem{theo}{Theorem}[section]
\newtheorem{theointro}{Theorem}

\newaliascnt{lemm}{theo}
\newaliascnt{prop}{theo}
\newaliascnt{coro}{theo}
\newaliascnt{rema}{theo}
\newaliascnt{defi}{theo}
\newaliascnt{exam}{theo}
\newaliascnt{conj}{theo}
\newtheorem{lemm}[lemm]{Lemma}
\newtheorem{prop}[prop]{Proposition}
\newtheorem{coro}[coro]{Corollary}
\newtheorem*{coro*}{Corollary}
\theoremstyle{definition}
\newtheorem{rema}[rema]{Remark}
\newtheorem*{rema*}{Remark}
\newtheorem{defi}[defi]{Definition}

\makeatletter
\renewcommand{\chapterautorefname}{\S\@gobble}
\renewcommand{\sectionautorefname}{\S\@gobble}
\renewcommand{\subsectionautorefname}{\S\@gobble}
\renewcommand{\subsubsectionautorefname}{\S\@gobble}
\makeatother

\aliascntresetthe{lemm}
\aliascntresetthe{prop}
\aliascntresetthe{coro}
\aliascntresetthe{rema}
\aliascntresetthe{defi}
\aliascntresetthe{exam}
\aliascntresetthe{conj}

  \DeclareFontFamily{U}{mathx}{\hyphenchar\font45}
\DeclareFontShape{U}{mathx}{m}{n}{
      <5> <6> <7> <8> <9> <10>
     <10.95> <12> <14.4> <17.28> <20.74> <24.88>
    mathx10
      }{}
\DeclareSymbolFont{mathx}{U}{mathx}{m}{n}
\DeclareMathAccent{\widecheck}{\mathalpha}{mathx}{"71}

\begin{document}
\date{24 December, 2025}
\begin{abstract}
We show in this paper that in the context of graded Lie algebras, all cuspidal character sheaves arise from a nearby-cycle construction followed by a Fourier--Sato transform in a very specific manner. Combined with results of the last two named authors, this completes the explicit description of cuspidal character sheaves for Vinberg's type I graded classical Lie algebras. 
\end{abstract}
\maketitle
\setcounter{tocdepth}{1}
 \tableofcontents

\section*{Introduction}

Let $G$ be a connected reductive complex algebraic group and $\fg$ its Lie algebra. Character sheaves on $G$ were introduced by Lusztig in \cite{L85a} as a sheaf-theoretic tool to study characters of finite groups of Lie type.  In \cite{MV88} Mirkovi\'c and the third-named author gave  a micro-local characterization of character sheaves as simple $G$-equivariant perverse sheaves on $G$ with nilpotent singular support, i.e., the singular support lies in $G\times (\fg^*)^{\nil}\subseteq T^*G$.  In light of this characterization, one can define a character sheaf on the Lie algebra $\fg$ to be a simple $G$-equivariant perverse sheaf on $\fg$ whose singular support lies in $\fg\times (\fg^*)^{\nil}$. Perverse sheaves with nilpotent singular support arise in many contexts. Invariant systems of differential equations, in particular those satisfied by characters and relative characters of representations, are important examples. Character sheaves are irreducible composition factors of such systems.  

In~\cite{L84,L87}, Lusztig gave a complete description of character sheaves on $\fg$. This can be viewed as a generalization of Springer theory in \cite{S78}, where irreducible representations of the Weyl group of $G$ are realized in the cohomology of Springer fibers. In particular, Lusztig determined the cuspidal character sheaves on $\fg$, i.e., those that cannot be obtained as a direct summand, up to a shift, of parabolic induction from smaller groups. He shows that cuspidal character sheaves on $\fg$ are very rare, and that they have nilpotent support modulo the centre of $\fg$.

This paper focuses on character sheaves in the setting of graded Lie algebras. Let $\theta:G\to G$ be an automorphism of order $m\in \Z_{>0}$. Let $G_0$ be the neutral component of the subgroup $G^{\theta}$ and let $\fg_k$ denote the $\zeta^k_m$-eigenspace of $\theta$ for $k\in \Z/m\Z$, where $\zeta_m = e^{2\pi i / m}$. Analogously to the ungraded $m=1$ case, we define a character sheaf on $\fg_1$ to be a simple $G_0$-equivariant perverse sheaf on $\fg_1$ whose singular support lies in $\fg_1\times (\fg_1^*)^{\nil}$. Here $G_0$ acts on $\fg_1$ via the adjoint action.

We give a uniform construction of cuspidal character sheaves on $\fg_1$, i.e., those that can not be obtained as direct summand by proper parabolic induction. We show that all cuspidal character sheaves arise in a unique way from a nearby-cycle construction followed by a Fourier--Sato transform applied to an explicit set of input data. In the case of involutions, $m=2$, the classification and the explicit determination of  character sheaves on $\fg_1$, in particular the cuspidal ones, have been carried out in \cite{VX22,VX23sl,X21} for Lie algebras of classical types by an exhaustion method.  
In~\cite{VX23st}, the last two named authors began the study of character sheaves on graded Lie algebras for $m>2$ with the goal of explicitly determining the cuspidal ones. It was conjectured that when the $G_0$-action on $\fg_1$ is GIT-stable the cuspidal character sheaves are exactly the character sheaves with support equal to all of $\fg_1$. Recall that GIT-stability amounts to the existence of a semisimple element of $\fg_1$ with finite stabiliser in $G_0$. The main result of this paper confirms this conjecture. It also provides a new proof that the list of cuspidal character sheaves in \cite{VX22,VX23sl,X21} for the involution case is complete, avoiding the exhaustion step.

For type I graded classical  Lie algebras (in the sense of Vinberg~\cite{Vi77}), cuspidal character sheaves are constructed explicitly in~\cite{VX23st,VX21,X24}, making use of the nearby-cycle construction in~\cite{GVX20} in a slightly generalized form. The fact that these character sheaves, as anticipated, are cuspidal and constitute all of them follows from the results of this paper. Unlike in the ungraded setting, these cuspidal character sheaves very rarely have nilpotent support modulo centre. They arise in families parametrized by irreducible representations of Hecke algebras attached to complex reflection groups with parameters at certain roots of unity. The classification of cuspidal character sheaves  can be viewed as a first and most important step towards a generalized Springer theory in the graded setting.

In addition, our results confirm that there are no cuspidal character sheaves in Vinberg's type II graded classical Lie algebras, and that all cuspidal character sheaves have nilpotent supports in Vinberg's type III graded classical Lie algebras, as expected in~\cite{X24}. These consequences are discussed in~\cite{LVX}. 

\subsection*{Main results}
Let us explain the results of this paper in a bit more detail, see~\autoref{sec:main} for the precise statements. Let $G, \fg, \theta$ and $m$ be as above. For simplicity, we assume here that $G$ is \emph{semisimple}.

In the ungraded setting, Lusztig shows that every cuspidal character sheaf on $\fg$ is an intersection cohomology sheaf given by a clean local system on a distinguished nilpotent orbit. In the graded setting,  we define a natural stratification of $\fg_1$ into $G_0$-stable subvarieties such that every character sheaf is constructible with respect to this stratification.  Our first main result is:
\begin{theointro}[\autoref{thm:cuspidal-distinguished}]\label{theo:A}
Let $\cF= \IC(\cF^{\circ})$ be a cuspidal character sheaf on $\fg_1$ associated with a  local system $\cF^{\circ}$ on a stratum $\cS\subseteq\fg_1$. 
\begin{enumerate}
\item
The adjoint quotient $\cS/\!/ G_0\subseteq \fg_1 /\!/ G_0$ is of dimension $\dim\fg_1 - \dim \fg_0$. 
\item
The stratum $\cS$ consists of {$G_0$-}distinguished elements (see~\autoref{distinguished}).
\item For every $G_0$-orbit $\rc\subseteq \cS$, the restriction $\cF^{\circ}|_{\rc}$ is a clean local system  and $\on{IC}(\cF^{\circ}|_{\rc})$ has nilpotent  singular support. 
\end{enumerate}
\end{theointro}

The following statements are immediate consequences of the above theorem.
\begin{enumerate}
\item[(i)] If $\dim\fg_1 < \dim \fg_0$, then there are no cuspidal character sheaves on $\fg_1$.
\item[(ii)]  If $\dim\fg_1 = \dim \fg_0$, then every cuspidal character sheaf (if exists) has nilpotent support. Moreover, every cuspidal character sheaf is of the form $\IC(\cF^{\circ})$ for some clean local system $\cF^{\circ}$ on a $G_0$-distinguished nilpotent orbit in $\fg_1$.
\item[(iii)]  If there exists a cuspidal character sheaf with nilpotent support then $\dim\fg_1 = \dim \fg_0$.
\item[(iv)] If the $G_0$-action on $\fg_1$ is GIT-stable, then all cuspidal character sheaves have full support.
\end{enumerate}

Character sheaves on a $\Z$-graded reductive Lie algebra $\fg=\oplus_{i\in\mathbb{Z}}\fg_i$ have been studied in \cite{L95}.  In this case  Lusztig shows that a cuspidal character sheaf on $\fg_1$ is clean and is supported on a \emph{distinguished} nilpotent $G_0$-orbit, the unique open dense orbit in $\fg_1\subseteq\fg^{\nil}$. Moreover, a cuspidal character sheaf on $\fg_1$ always comes from a cuspidal character sheaf on the ungraded $\fg$ via restriction. The  latter were classified in~\cite{L84}. Using these works of Lusztig as input we give in \cite{LVX} a complete classification of cuspidal character sheaves with nilpotent support in our graded setting.  The cuspidal character sheaves with nilpotent support, which we also refer to as bi-orbital cuspidal character sheaves, on $\theta$-stable Levi subalgebras will play an essential role in our nearby-cycle construction, which we explain next.
\par

Our second main result shows that every cuspidal character sheaf arises as  a simple subsheaf and as a simple quotient sheaf of the Fourier--Sato transform of some nearby-cycle sheaf. The inputs of our nearby cycle construction are pairs $(\rc,\chi)$,  which we call nil-supercuspidal data (see~\autoref{def-sc/nsc}), consisting of a $G_0$-distinguished orbit $\rc$ and a certain clean  local system $\chi$ on $\rc$. Roughly, a nil-supercuspidal datum is given by a $\theta$-stable Levi subgroup $M$ together with a cuspidal character sheaf with nilpotent support on the (graded) derived Lie algebra of $M$. \par

For each nil-supercuspidal datum $(\rc,\chi)$, we define the corresponding nearby-cycle sheaf $\cP_\chi$ as follows (see \eqref{eqn-psi} for its precise definition):
\[
\cP_{\chi} := \lim_{t\to 0} \IC(t\,\rc, \chi)\in \Perv_{G_0}(\fg^{\nil}_1);
\]
here $\Perv_{G_0}(\fg^{\nil}_1)$ stands for $G_0$-equivariant perverse sheaves on the null-cone $\fg^{\nil}_1$. 
Instead of just the Fourier--Sato transform it is convenient for us to use the functor $\dag$ (see~\autoref{ssec:functor dag} for its precise definition), which is the composition of the Fourier--Sato transform, Verdier duality and a suitable isomorphism which identifies $G_0$-equivariant sheaves on $\fg_1$ with those on $\fg_{-1}\cong \fg_{1}^*$. We write $\cP_\chi^\dag\in\Perv_{G_0}(\fg_1)$ for the resulting sheaf.

\begin{theointro}[\autoref{thm:main}]\label{theo:B}
\leavevmode
\begin{enumerate}
    \item Every simple subsheaf, or quotient sheaf, of $\cP_{\chi}^{\dag}$ is a cuspidal character sheaf on $\fg_1$.
    \item  Every cuspidal character sheaf on $\fg_1$ arises as a subsheaf, or a quotient sheaf, of $\cP_{\chi}^{\dag}$ for some nil-supercuspidal datum $(\rc,\chi)$. 
    \item Let $(\rc,\chi)$ be a nil-supercuspidal datum and $\cS$ the stratum of $\fg_1$ containing $\rc$. There exists a $G_0$-equivariant local system $\cM_{\chi}$ on $\cS$ such that $$\cP^{\dag}_{\chi} \cong \IC(\cM_{\chi}).$$
    \item Given nil-supercuspidal data $(\rc, \chi)$ and $(\rc', \chi')$, the corresponding nearby-cycle sheaves $\cP_{\chi}$ and $\cP_{\chi'}$ are isomorphic if $(\rc, \chi)$ and $(\rc', \chi')$ are equivalent (in the sense of~\autoref{def:equiv}); otherwise, we have $\Hom(\cP_{\chi}, \cP_{\chi'}[k]) = 0$ for $k\in \Z$. 
\end{enumerate} 
\end{theointro}
\begin{rema*}
Given a cuspidal character sheaf $\cF=\on{IC}(\cF^\circ)$ where $\cF^\circ$ is a local system on the stratum $\cS$, the nil-supercuspidal datum $(\rc,\chi)$ in (2) can be chosen as, in the notation of~\autoref{theo:A}, a $G_0$-orbit $\rc\subseteq \cS$  and a direct summand  $\chi$ of $\cF^{\circ}|_{\rc}$. 
\end{rema*}
As a consequence of the above theorem, we have the following Springer-type correspondence for cuspidal character sheaves, which is a major step towards their classification: \begin{coro*}[\autoref{coro:springer}] There is a natural one-to-one correspondence between the set of isomorphism classes of cuspidal character sheaves on $\fg_1$ and the set
\[
\bigsqcup_{(\rc, \chi)/\sim} \Irr\left(\End(\cM_{\chi})\mof\right)
\]
where $(\rc, \chi)/\sim$ runs over the equivalence classes of nil-supercuspidal data.
\end{coro*}

In~\cite{LVX} we make use of Theorems~\ref{theo:A} and \ref{theo:B} to classify all the gradings that afford cuspidal character sheaves and determine the supports of the cuspidal character sheaves. We also classify the nil-supercuspical data. This reduces the classification of cuspidal character sheaves to the calculation of the local system $\cM_{\chi}$ for each nil-supercuspidal datum $(\rc, \chi)$. This is done in \cite{VX23st} for GIT-stable gradings and in \cite{VX21,X24} for, not necessarily GIT-stable, type I graded classical Lie algebras. 

We will derive \autoref{theo:B}   making use of the following adjunction formula (\autoref{thm-adjunction}):
\[
  \dim\Hom\left(\cP_{\chi} ,\cG[k]\right)=\dim\Hom\left(\cG^{\dag}, \IC(\chi)[N + k]\right)
\]
for every simple perverse sheaf $\cG \in \Perv_{G_0}(\fg^{\nil}_1)$ and every $k\in \Z$, where $N=\dim\fg_1-\dim\fg_0$. Moreover, we show that both sides vanish for $k\neq 0$. As an interesting by-product, we show that the nearby-cycle sheaf $\cP_{\chi}$ is both projective and injective in the category $\Perv_{G_0}(\fg^{\nil}_1)$. \par

It may be useful to regard the nearby-cycle construction $(\rc, \chi)\mapsto \cP_{\chi}$ as a twisted induction, in analogy with the Deligne--Lusztig induction for representations of finite groups of Lie type. It would be interesting to compare with the twisted induction in \cite{NY25}  where character sheaves on loop algebras are constructed partly adapting the construction of Yu \cite{Yu01} for representations of p-adic reductive groups.

\par

\begin{rema*}
Let us point out that the nearby-cycle construction is trivial and produces nothing new when $\dim\fg_1 = \dim \fg_0$. This situation includes the ungraded setting $m = 1$ and the $\Z$-graded setting where $m \to\infty$ as special cases. 
    
\end{rema*}

\subsection*{Connection with affine Springer fibres and double affine Hecke algebras}

Our study of cuspidal character sheaves is partly motivated by the connection between (homogeneous elliptic) affine Springer fibres and finite-dimensional modules of double affine Hecke algebras (see for example~\cite{VV09,OY16}).

Our setting of graded Lie algebras arises naturally from Moy-Prasad filtrations of loop groups and loop Lie algebras (see for example~\cite{LY18}). In more detail, given a simply connected  almost simple group $G$ and a $\Z/m\Z$-grading $\fg_*$ (assumed to arise from an inner automorphism $\theta\in \Ad(G)$ for simplicity) on its Lie algebra $\fg = \Lie G$,  consider the loop group $LG = G(\C(\!(\varpi)\!))$ and the loop Lie algebra $L\fg = \fg\otimes\C(\!(\varpi)\!)$. There exists a $\Cc$-action on $LG$ by an automorphism called $\lambda$, such that there are natural isomorphisms
\[
(LG)^{\lambda} \cong G_0,\quad \li^{\lambda}_i(L\fg) \cong \fg_{i+m\Z} \text{ for $i\in \Z$}. 
\]
There is an affine analogue $\pi: \widetilde{L\fg}\to L\fg$ of the Grothendieck-Springer resolution  and the $\Cc$-action $\lambda$ lifts naturally to $\widetilde{L\fg}$. 
If we take the $\lambda$-homogeneous component of $\pi$, denoted $\li^{\lambda}_i\pi:\li^{\lambda}_i(\widetilde{L\fg})\to \li^{\lambda}_i(L\fg)$, then $\cI:=(\li^{\lambda}_i\pi)_*\C$ is a semisimple complex of sheaves on $\li^{\lambda}_i(L\fg) = \fg_i$ whose simple perverse constituents are character sheaves. Vasserot~\cite{Va05} constructed an action of a double affine Hecke algebra $\mathcal{H}$ on $\cI$, which induces a one-to-one correspondence between simple perverse constituents of $\cI$ and simple objects of a  category $\on{O}_{\lambda}(\mathcal{H})$ of certain $\mathcal{H}$-modules.

In~\cite{LY18} Lusztig and Yun obtain a block decomposition of the equivariant derived category $\Db_{G_0}(\fg_i^{\nil})=\oplus_{\xi}\Db_{G_0}(\fg_i^{\nil})_\xi$ and conjecture that the simple objects in each block $\Db_{G_0}(\fg_i^{\nil})_\xi$ are parametrized by simple modules of a degenerate double affine Hecke algebra $\cH_\xi$ (associated to the block) that lie the category $\cO_\lambda(\cH_\xi)$ of integrable $\cH_\xi$-modules. This conjecture was proved by the first-named author in~\cite{Liu23}. These results generalise Vasserot's work, which corresponds to the {\em principal} Lusztig--Yun (LY-)block. 

We show that the finite dimensional simple modules that lie in $\cO_\lambda(\cH_\xi)$ correspond exactly bijectively with the cuspidal character sheaves that lie in the block $\Db_{G_0}(\fg_1)_\xi$, that is, we have a bijection
$$\Irr^{\on{fd}}(\cO_\lambda(\cH_\xi))\to\Char_{G_0}^{\on{cusp}}(\fg_1)_\xi\,.$$

On the other hand, the works~\cite{VV09,OY16} provide strong evidence of the connection between finite-dimensional simple modules of double affine Hecke algebras and the geometry of affine Springer fibers. 
Assume that $\gamma\in L\fg$ is regular semisimple, topologically nilpotent and elliptic. Kazhdan--Lusztig~\cite{KL88} showed that the affine Springer fibre $\on{Sp}_{\gamma}:=\pi^{-1}(\gamma)$ is of finite type over $\C$. Assume furthermore that $\gamma\in \li^{\lambda}_i(L\fg)$ for some $i>0$ (in other words, the grading $\fg_*$ is GIT-stable). Then, $\lambda$ induces a $\Cc$-action on the fibre $\pi^{-1}(\gamma)$, still denoted by $\lambda$. Varagnolo--Vasserot~\cite{VV09} constructed a representation of $\cH = \cH_{\xi_0}$ (where $\xi_0$ denotes the principal LY-block) in the cohomology group $\rmH^*(\on{Sp}_{\gamma}^{\lambda})$ of the fixed-point subvariety $\on{Sp}_{\gamma}^{\lambda}$. It gives rise to finite-dimensional simple modules in $\on{O}_{\lambda}(\mathcal{H})$. As an application of our results on cuspidal character sheaves, we prove in~\autoref{theo:DAHA-fd} that, conversely, every simple finite-dimensional module in $\on{O}_{\lambda}(\mathcal{H})$ appears as a constituent in $\rmH^*(\on{Sp}_{\gamma}^{\lambda})$ and we generalise this result to all LY-blocks.  
\subsection*{Structure of this paper}

In~\autoref{sec:graded-Lie}, we recollect basic notions we will use in this paper: orbital, anti-orbital, bi-orbital complexes and character sheaves on $\fg_1$,  parabolic induction, restriction, co-restriction, and hyperbolic restriction. We also introduce  cuspidal and supercuspidal orbital and anti-orbital sheaves. 

 In \autoref{sec:support-character-sheaves}, we define and study a stratification on $\fg_1$ along which every character sheaf is constructible. We describe the support of a character sheaf and analyse the restriction of a character sheaf to a generic slice of its support. 
 
In \autoref{sec:main}, we define the nearby-cycle sheaf $\cP_{\chi}$ attached to a pair $(\rc, \chi)$ consisting of a $G_0$-orbit $\rc$ and a simple $G_0$-equivariant local system $\chi$ on $\rc$.    
We state and prove the main theorem \autoref{thm:main} (\autoref{theo:B}) making use of technical results \autoref{prop:parabolic-vanishing} (annihilation of $\cP_{\chi}$  by parabolic restriction), \autoref{thm:cuspidal-distinguished} (\autoref{theo:A}) (and its analog~\autoref{prop:annihlation-distinguished}),  and \autoref{thm-adjunction} (adjunction formula), which are proved in \autoref{sec:parabolic-vanishing}, \autoref{sec:support}, and \autoref{sec:proof-adjunction} respectively. 

 In \autoref{sec:LY} we recollect some basic constructions and results of Lusztig--Yun's theory from \cite{LY17a}, which we make ample use of in the subsequent sections. We deduce a criterion for a supercuspidal character sheaf to be bi-orbital. 
 
 The most technical parts of this paper are in \autoref{sec:support} and \autoref{sec:proof-adjunction}. 
In \autoref{sec:support} we prove \autoref{thm:cuspidal-distinguished} and \autoref{prop:annihlation-distinguished}. The latter applies to  anti-orbital sheaves annihilated by proper parabolic restrictions, in particular, to $\cP_\chi^\dag$. Its proof is simpler and  illustrates the ideas behind the proof of \autoref{thm:cuspidal-distinguished}. The crucial part of the proof of \autoref{thm:cuspidal-distinguished} is the construction called semi-orthogonal approximation. It transforms a given cuspidal character sheaf into a complex annihilated by parabolic restrictions up to a given perverse degree. Unlike in the ungraded setting, cuspidal character sheaves in our setting are in general no longer annihilated by proper parabolic restrictions.  The semi-orthogonal approximation allows us to view cuspidal character sheaves as if they were annihilated by parabolic restrictions.  
In \autoref{sec:proof-adjunction} we prove \autoref{thm-adjunction}. The proof has three main ingredients: Beilinson's construction of unipotent nearby cycles, the hyperbolic restriction theorem applied to spiral restrictions, and cohomological duality for differential graded algebras. \par

In \S\ref{sec:applications} we discuss applications to classification of finite dimensional simple modules of degenerate double affine hecke algebras and cohomology of regular elliptic homogeneous affine Springer fibers.

\section*{Acknowledgment}
 We would like to thank  M. Grinberg, C. Mautner, G. Schwarz, P. Shan, E. Vasserot  for inspiring and valuable discussions.

\section{Perverse sheaves on graded Lie algebras}\label{sec:graded-Lie}
In this section, we recall some previously known constructions and results related to graded Lie algebras and perverse sheaves on them.

\subsection{Basic setup}

Let $G$ be a connected reductive group over $\C$, $\theta\in \Aut(G)$ an automorphism and $m \in \Z_{\ge 1}$ a positive integer such that $\theta^m = \id$. We write $\zeta_m = e^{2\pi i/m}$. Then $\theta$ induces a $\Z/m\Z$-grading on the Lie algebra $\fg$ of $G$
\[
\fg = \bigoplus_{k\in \Z/m\Z} \fg_k,\quad \fg_k = \{x\in \fg\;|\; \theta x = \zeta_m^k x\}.
\]
For the simplicity of notation, we will write $\fg_k = \fg_{k+m\Z}$ for $k\in \Z$.

Let $G_{0} = (G^{\theta})^\circ$ denote the neutral component of the $\theta$-fixed-point subgroup of $G$. Let $\fg^{\nil}\subseteq \fg$ be the nilpotent cone and write $\fg^{\nil}_k = \fg^{\nil}\cap \fg_k$ for each $k\in \Z/m\Z$. The adjoint action of $G$ on $\fg$ restricts to an action of $G_{0}$ on $\fg_k$, which preserves $\fg^{\nil}_k$. By~\cite{Vi77}, the number of $G_{0}$-orbits in $\fg^{\nil}_k$ is finite; moreover, $\fg^{\nil}_k$ coincides with the pre-image of $0$ under the adjoint quotient map $\fg_k\to \fg_k/\!/G_0$.

Let $\Db_{G_{ 0}}(\fg_{k})$ (resp. $\Db_{G_{ 0}}(\fg^{\nil}_{k})$) be the $G_{ 0}$-equivariant derived category of sheaves of $\C$-vector spaces on $\fg_{k}$ (resp. on $\fg_{k}^\nil$) with bounded and constructible cohomology. The perverse t-structures on $\Db_{G_{ 0}}(\fg_{k})$ etc. are normalised so that the functor $\Db_{G_{ 0}}(\fg_{k})\to\Db(\fg_{k})$ forgetting the equivariance is perverse t-exact. 

Let $\Perv_{G_{ 0}}(\fg_{k})\subseteq \Db_{G_{ 0}}(\fg_{k})$ etc. denote the heart of the perverse t-structure. The Verdier duality functor $$\D:\Db_{G_{ 0}}(\fg_{k})^{\op}\to \Db_{G_{ 0}}(\fg_{k})$$ is normalised so that it preserves the perverse t-structure. Given a $G_0$-equivariant local system $\cL$ of $\C$-vector spaces on a smooth $G_0$-stable subvariety $U\subseteq \fg_k$, we denote by $\IC(\cL)\in \Perv_{G_0}(\fg_k)$
the intersection cohomology complex with coefficients in $\cL$, namely $\IC(\cL) = j_{!*}\cL[\dim U]$, where $j:U\to \fg_k$ is the inclusion and $j_{!*}$ is the intermediate extension along $j$. We use similar notation in other contexts. \par

Given a nilpotent element $x\in \fg^{\nil}_k$, a $G_0$-equivariant local system on the orbit $G_0 x$ is equivalent to a representation of $\pi_0(Z_{G_0}(x))$, the component group of the stabiliser. The assignment $(\cO, \cL)\mapsto \IC(\cL)$ yields a bijection from the pairs $(\cO, \cL)$ of a $G_0$-orbit $\cO\subseteq \fg^{\nil}_k$ and an irreducible $G_0$-equivariant local system $\cL$ on $\cO$, defined up to isomorphism, to the isomorphism classes of simple objects of $\Perv_{G_0}(\fg^{\nil}_k)$.

\begin{rema}
Note that although we consider here $\Z/m\Z$-graded Lie algebras, the discussion applies equally well to the cases where $\fg$ is ungraded  $m = 1$ case and the $\Z$-graded $m =\infty$ case. The latter case, where $\fg_k^{\nil} = \fg_k$ for $k\in \Z\setminus\{0\}$, is studied in~\cite{L95}.
\end{rema}

\subsection{Fourier--Sato transform and monodromic complexes}\label{ssec:FS}

The adjoint $G_0$-action on $\fg_k$ can be enriched into a $(G_0\times \Cc)$-action by 
\[
    (G_0\times \Cc)\times \fg_k\to \fg_k,\quad ((g, r), z)\mapsto r^{-2}\Ad_g z.
\]
The weight $-2$ is inserted for the purposes of \autoref{lem:orbital-monodromic} below. 
We define $\Db_{G_{ 0}}(\fg_{ k})^{\Cunip}$  to be the thick triangulated subcategory of $\Db_{G_{ 0}}(\fg_{ k})$ generated by the image of the natural functor $\Db_{G_{ 0}\times \Cc}(\fg_{ k})\to \Db_{G_{ 0}}(\fg_{ k})$ which forgets the $\Cc$-equivariance. Its objects are called \emph{unipotently $\Cc$-monodromic complexes}.
\begin{lemm}\label{lem:orbital-monodromic}
There is an inclusion $\Db_{G_{ 0}}(\fg_{k}^{\nil})\subseteq \Db_{G_{ 0}}(\fg_{k})^{\Cunip}$. 
\end{lemm}
\begin{proof}
The graded Jacobson--Morosov theorem (see \cite[\S 2.3]{LY17a}) implies that every $G_0$-orbit of $\fg^{\nil}_{k}$ is stable under the extra $\Cc$-action and moreover, for each $x\in \fg^{\nil}_{k}$, the inclusion $Z_{G_{0}}(x) \subseteq Z_{G_{0}\times \Cc}(x)$ induces an isomorphism on the component groups 
\[
    \pi_0(Z_{G_{0}}(x)) \cong \pi_0(Z_{G_{0}\times \Cc}(x)), 
\]
see \cite[\S 2.1.f]{L88}. 
Therefore, every simple object of $\Perv_{G_{ 0}}(\fg_{k}^{\nil})$ admits a unique $(G_{ 0}\times \Cc)$-equivariant enhancement. In particular, simple objects of $\Perv_{G_{ 0}}(\fg_{k}^{\nil})$ lie in $\Db_{G_{ 0}}(\fg_{k})^{\Cunip}$. Since they generate $\Db_{G_{ 0}}(\fg_{k}^{\nil})$ as a triangulated category, the lemma follows.
\end{proof}

We fix once and for all 
a non-degenerate $(G, \theta)$-invariant bilinear form on $\fg$ (for example, the trace pairing of any faithful representation of $G\rtimes \langle \theta\rangle$). It induces a $G_0$-equivariant linear isomorphism $(\fg_k)^* \cong \fg_{-k}$. The Fourier--Sato transform yields a perverse-t-exact triangle equivalence
\[
    \FS: \Db_{G_{ 0}}(\fg_{k})^{\Cunip} \to \Db_{G_{ 0}}(\fg_{-k})^{\Cunip}.
\]
For this fact, see \cite[\S 6]{Bry86} and \cite[\S 3.7]{KS90} for example. 

\subsection{Orbital, anti-orbital and bi-orbital complexes}\label{subsec:orbital}
\begin{defi}
A complex $\cK\in \Db_{G_{0}}(\fg_{k})$ is called \emph{orbital} if it lies in the full subcategory $\Db_{G_{0}}(\fg_{k}^{\nil})$; it is called \emph{anti-orbital} if it is isomorphic to the Fourier--Sato transform of an orbital complex on $\fg_{-k}$; it is called \emph{bi-orbital} if it is both orbital and anti-orbital.
\end{defi}
Let $\Dbn_{G_{0}}(\fg_{k})$ denote the full subcategory of $\Db_{G_{0}}(\fg_{k})$ spanned by anti-orbital complexes and define $\Perv_{G_{0}}^{\nil}(\fg_{k}) = \Dbn_{G_{0}}(\fg_{k})\cap \Perv_{G_{0}}(\fg_{k})$. It follows immediately from the definition that the Fourier--Sato transform restricts to an equivalence between orbital and anti-orbital complexes:
\[
\FS: \Db_{G_0}(\fg^{\nil}_k)\xrightarrow{\sim} \Dbn_{G_0}(\fg_{-k}). 
\]
\begin{defi}
A \emph{character sheaf} on $\fg_ {k}$ is a simple object of $\Perv^{\nil}_{G_{0}}(\fg_ {k})$, i.e., a simple anti-orbital perverse sheaf on $\fg_k$.
\end{defi}
The set of isomorphism classes of character sheaves on $\fg_ {k}$ is denoted by $\Char_{G_{0}}(\fg_ {k})$. \par

\begin{rema}\label{rem:geometric-origin}
Note that simple orbital sheaves are intersection cohomology complexes with coefficients in local systems of finite monodromy; therefore, they are of geometric origin. In addition, $\FS$ admits a formula in terms of vanishing cycles, see \cite[(10.3.31)]{KS90}, for example, which implies that it preserves the subcategories of semisimple complexes of geometric origin. Consequently, every semisimple complex in $\Db_{G_0}(\fg^{\nil}_k)$ or $\Dbn_{G_0}(\fg_{k})$ is of geometric origin --- this can also be easily deduced from Lusztig--Yun's construction of orbital and anti-orbital sheaves, see the reminder in \autoref{sec:LY}.
\end{rema}
\subsection{Parabolic induction, restriction, and co-restriction}
 
Let $P\subseteq G$ be a $\theta$-stable parabolic subgroup and let $U\subseteq P$ denote its unipotent radical.  Let $L= P / U$ denote the Levi quotient and let $U_0 = U^{\theta}$. The set of $\theta$-stable Levi factors of $P$ is a principal homogeneous $U_0$-space under conjugation; every such Levi factor is canonically isomorphic to the Levi quotient $L$. 
 We set 
 \[
 \bega
 P_0 = G_0\cap P,\;L_0 = P_0/ U_0,\\
 \fp=\on{Lie}P,\;\fu = \Lie U,\; \fl=\on{Lie}L,\; \fp_k = \fg_k \cap \fp,\;  \fu_k = \fg_k \cap \fu,\;  \fl_k = \fp_k/\fu_k.
 \eega
 \]
 We have the following diagram of stacks:
\[
[\fl_ {k}/L_0] \xleftarrow{a} [\fp_ {k}/P_0]  \xrightarrow{b} [\fg_ {k} / G_{0}].
\]
 Parabolic induction, restriction and co-restriction are defined to be the following functors:  
\beq\label{parabolic functors}
\begin{gathered}
\Ind^{\fg_ {k}}_{\fp_ {k}}: \Db_{L_{0}}(\fl_ {k})\rightleftarrows \Db_{G_{0}}(\fg_ {k}):\Res^{\fg_ {k}}_{\fp_ {k}}, \res^{\fg_ {k}}_{\fp_ {k}}\\ 
\Ind^{\fg_ {k}}_{\fp_ {k}} = b_!a^*[\dim \fu_k + \dim \fu_0] \cong b_*a^![-\dim \fu_k - \dim \fu_0], \\
\Res^{\fg_ {k}}_{\fp_ {k}} = a_!b^*[\dim \fu_k + \dim \fu_0],\; \res^{\fg_ {k}}_{\fp_ {k}} = a_*b^![-\dim \fu_k - \dim \fu_0].
\end{gathered}
\eeq
We have the following adjunctions: 
\beq\label{eqn-adjunctions-parabolic}
\Res^{\fg_{k}}_{\fp_k}\dashv\Ind^{\fg_{k}}_{\fp_k}\dashv\res^{\fg_{k}}_{\fp_k}.
\eeq
The functor $\Ind^{\fg_1}_{\fp_1}$ is normalised in such a way that it commutes with the Verdier duality. Moreover, arguments similar to \cite[\S 10.2, \S 10.5]{L95} show that these functors commute with the Fourier--Sato transform:
\beq\label{eqn-parabolic-FS}
\FS\circ\Ind^{\fg_{k}}_{\fp_k}\cong \Ind^{\fg_{-k}}_{\fp_{-k}}\circ\FS,\quad
\FS\circ\Res^{\fg_{k}}_{\fp_k}\cong \Res^{\fg_{-k}}_{\fp_{-k}}\circ\FS,\quad
\FS\circ\res^{\fg_{k}}_{\fp_k}\cong \res^{\fg_{-k}}_{\fp_{-k}}\circ\FS.
\eeq

It is easy to see that the pre-image of $\fl_k^{\nil}$ under the projection $\fp_k\to \fl_k$ coincides with the intersection $\fp_k\cap \fg^{\nil}_k$. Therefore, the functors $\Ind^{\fg_{k}}_{\fp_k}, \Res^{\fg_{k}}_{\fp_k}$ and $\res^{\fg_{k}}_{\fp_k}$ preserve the subcategories of orbital and anti-orbital complexes.

\begin{rema}
Note that the functors $a_!, b_!, a^!, b^!$ above are normalised in the stacky way. In the more classical language:
\[
    \fl_ {k} \xleftarrow{\til a} \fp_ {k}  \xrightarrow{\til b} \fg_ {k},
\]
we have $\Res^{\fg_ {k}}_{\fp_ {k}} = \til a_!\til b^*[\dim \fu_k - \dim \fu_0]$ and $\res^{\fg_ {k}}_{\fp_ {k}} = \til a_*\til b^![-\dim \fu_k + \dim \fu_0]$.
\end{rema}
\begin{rema}\label{rem:ss}
The semisimple orbital and anti-orbital complexes are of geometric origin (see \autoref{rem:geometric-origin}).  
Thus, the parabolic induction sends semisimple orbital and anti-orbital complexes to semisimple ones. Making use of the hyperbolic restriction theorem (see~\autoref{prop:hyperbolic-parab}), we see that the parabolic restriction and co-restriction also preserve the semisimplicity of orbital and anti-orbital complexes.
\end{rema}

\subsection{Hyperbolic restriction}\label{ssec:hyperbolic}

Let $(\fp, \fl, \fp')$ be a triplet formed by a pair of opposite $\theta$-stable parabolic subalgebras $(\fp, \fp')$ such that the intersection $\fp\cap \fp' = \fl$ is a common $\theta$-stable Levi factor. Recall the functors defined in~\eqref{parabolic functors}.

\begin{prop}\label{prop:hyperbolic-parab}
There is an isomorphism of functors from $\Db_{G_{0}}(\fg_ {1})$ to $\Db_{L_{0}}(\fl_ {1})$:
\[
\Res^{\fg_ {1}}_{\fp_ {1}} \cong \res^{\fg_ {1}}_{\fp'_ {1}}[D],
\]
where $D = \dim \fg_1 - \dim \fg_0 - \dim \fl_1 + \dim \fl_0$.
\end{prop}
\begin{proof} 
By \cite[\S 1.6]{Liu24}, we can choose a cocharacter $\lambda\in\bX_*(G_{0})$ such that 
\[
\fp = \li^\lambda_{\ge 0}{\fg},\quad\fl = \li^\lambda_{0}{\fg},\quad \fp' = \li^\lambda_{\le 0}{\fg}.
\]
 Consider the diagram:
\[
\begin{tikzcd}
\fg_ {1} & \fp_ {1} \arrow{l}{i}\arrow[bend left=60]{d}{p}\\
\fp'_ {1} \arrow{u}{i'}\arrow[bend right=60,swap]{r}{p'} & \fl_ {1}\arrow[swap]{l}{j'}\arrow{u}{j}
\end{tikzcd}.
\]
It can be viewed as a hyperbolic restriction diagram with respect to the $\Cc$-action on $\fg_ {1}$ given by $\Ad_{\lambda(t)}$. Since $j^!=p_!$, $j^*=p_*$, $(j')^!=p'_!$, $(j')^*=p'_*$, the desired isomorphisms follow from the hyperbolic restriction theorem \cite[Theorem1]{B03}. 
\end{proof}
\begin{rema}\label{rema:non-construcitble}
    Neither the constructibility nor the boundedness is needed in the proof of the hyperbolic restriction theorem. In other words, the isomorphism $\Res^{\fg_ {1}}_{\fp_ {1}} \cong \res^{\fg_ {1}}_{\fp'_ {1}}[D]$ holds on the entire unbounded $G_0$-equivariant derived category $\operatorname{D}_{G_0}(\fg_1)$ of sheaves of $\C$-vector spaces, see~\cite{R19}.
\end{rema}

\subsection{Cuspidal and supercuspidal sheaves}
In this subsection we introduce the main objects of this paper: cuspidal and supercuspidal sheaves.

 \subsubsection{Cuspidal sheaves}\label{ssec:cuspidal}
\begin{defi}
\label{cuspidal:def}

A simple orbital (resp. anti-orbital) perverse sheaf $\cF$ on $\fg_ {k}$ is called \emph{cuspidal} if it does not arise as a direct summand of a parabolically induced orbital (resp. anti-orbital) complex. In other words, $\cF$ is cuspidal if there exists no triplet $(P, L, \cK)$, where $P$ is a $\theta$-stable proper parabolic subgroup of $G$ with $\theta$-stable Levi factor $L$ and $\cK\in \Db_{L_{0}}(\fl^{\nil}_ {k})$ (resp. $\cK\in \Dbn_{L_{0}}(\fl_ {k})$), such that $\cF$ is a direct summand of $\Ind^{\fg_{ k}}_{\fp_{ k}} \cK$. 

Note that in this definition, we can assume that the complex $\cK$ is semisimple. By  \autoref{rem:ss} the induction $\Ind^{\fg_{ k}}_{\fp_{ k}} \cK$ is then also semisimple. 
\end{defi}
Cuspidal simple anti-orbital perverse sheaves are called \emph{cuspidal character sheaves}. \par
\begin{rema}
Note that this notion of cuspidality differs from that in representation theory of $p$-adic groups --- a cuspidal  sheaf $\cF$ need not satisfy $\Res^{\fg_ {k}}_{\fp_ {k}}\cF = 0$ for every $\theta$-stable proper parabolic subalgebra $\fp\subseteq \fg$. 
\end{rema}
\subsubsection{Supercuspidal sheaves}\label{ssec:supercuspidal}
\begin{defi}\label{def:supercuspidal}
A simple orbital or anti-orbital perverse sheaf $\cF$ on $\fg_ {k}$ is called \emph{supercuspidal} if $\Res^{\fg_ {k}}_{\fp_ {k}}\cF = 0$ for every $\theta$-stable proper parabolic subgroup $P\subseteq G$. 
\end{defi}

\begin{defi}
\label{distinguished}
We say that an element $x\in \fg_{k}$ is \emph{$G_0$-distinguished} if it is not contained in a $\theta$-stable Levi factor of a $\theta$-stable proper parabolic subalgebra of $G$, or, equivalently, if the neutral component of $(Z_{G}(x)/Z(G))^{\theta}$ is unipotent. 

We say that a $G_0$-orbit $\cO\subseteq\fg_k$ is $G_0$-distinguished if it consists of $G_0$-distinguished elements.
\end{defi}

\begin{rema}
    Note that in~\cite{X24}, a different notion of distinguishedness was used, where commutator in $\fg_{-k}$ instead of $G_0$ was considered.
\end{rema}

The following proposition characterises supercuspidal orbital sheaves:
\begin{prop}
\label{characterisation of supercuspidal}
Given a pair $(\cO, \cL)$, where $\cO\subseteq\fg^{\nil}_ {k}$ is a nilpotent $G_{0}$-orbit and $\cL$ is an irreducible $G_{0}$-equivariant local system on $\cO$, the following conditions are equivalent:
\begin{enumerate}
\item[\rm{(i)}]
$\IC(\cL)$ is supercuspidal;

\item[\rm{(ii)}]
$\rmHc^*(\cO\cap (z + \fu_ {k}), \cL\mid_{\cO\cap (z + \fu_ {k})}) = 0$ for every $\theta$-stable proper parabolic subgroup $P\subseteq G$ with unipotent radical $U$ and for every $z\in \fp_ {k}$;
\item[\rm{(iii)}]
$\cO$ is $G_0$-distinguished and $\cL$ is clean.
\end{enumerate}
\end{prop}
\begin{rema}A pair $(\cO, \cL)$ satisfying the condition (ii) is called a \emph{supercuspidal pair} in  \cite{Liu24}.
\end{rema}
\begin{proof}
The implications (ii)~$\Leftrightarrow$~(iii) are proved in \cite[Theorem 5.1.2]{Liu24}. Given a pair $(\cO, \cL)$ satisfying (iii), the cleanness implies $j_{\cO!}\cL [\dim \cO] \cong \IC(\cL)$, where $j_{\cO}:\cO\hookrightarrow \fg_k$ is the inclusion of orbit. The proper base change theorem then implies the supercuspidality of $\IC(\cL)$. (i) $\Rightarrow$ (iii) is a straightforward consequence of the non-vanishing of Ext-groups between simple objects of $\Perv_{G_0}(\fg^{\nil}_k)$ \cite[Proposition 3.5.1]{Liu24}.
\end{proof}

\begin{rema}\leavevmode
\begin{enumerate}[label=(\roman*)]
\item
 Supercuspidal sheaves are cuspidal --- this follows immediately from the adjunction \eqref{eqn-adjunctions-parabolic}. The converse is, in stark contrast with the ungraded case ($m = 1$) and the $\Z$-graded case ($m = \infty$), not true in general.  
\item The classification of supercuspidal orbital sheaves can be reduced to that of cuspidal local systems on an (ungraded) Lie algebra; the latter problem has been solved by Lusztig in \cite{L84}. This will be done in a future publication.
\end{enumerate}
\end{rema}

In the bi-orbital case, the following criterion, which we will make use of, is proved in \cite[\S 7.1.1]{Liu24}:
\begin{theo}
 \label{Prop3} A bi-orbital character sheaf is cuspidal if and only if it is supercuspidal.   
\end{theo}

\begin{rema}
    Supercuspidal character sheaves are not necessarily bi-orbital. In \autoref{prop:supercuspidal-dim}, we provide a numerical criterion for a supercuspidal character sheaf to be bi-orbital.
\end{rema}

To classify cuspidal orbital/anti-orbital sheaves on $\fg_ {k}$, we can reduce the problem   to the situation where  $k = 1$. From now on, we will let $k=1$ unless otherwise stated.

\subsection{Some notation}Let $x\in \fg_1$. We will write $x = x_s + x_n$ for its Jordan decomposition with $x_s\in \fg_1$ semisimple, $x_n\in \fg_1$ nilpotent and $[x_s,x_n]=0$. We further write
\ber\label{eqn-hz}
&&Z_x= Z(Z_G(x_s)),\ H_x = [Z_G(x_s), Z_G(x_s)],\ \mathfrak{Z}_x=\on{Lie}Z_x,\  \fh_x=\on{Lie}H_x\\
&&\mathfrak{Z}_{x,1}=\mathfrak{Z}_x\cap\fg_1,\ \fh_{x,1}=\fh_x\cap\fg_1,\text{ and } H_{x,0}=(H_x^\theta)^\circ\nonumber.
\eer
Note that $Z_x=Z_{x_s}$ etc.

Let $\lambda\in \bX_*(G_0)$ and $n\in\Z$. We write
\beq
\li^\lambda_{n}{\fg}=\{z\in\fg\mid\on{Ad}_{\lambda(t)}z=t^nz,\,\forall\,t\in\Cc\}, \ \li^\lambda_{\geq n}{\fg}=\bigoplus_{l\geq n}\li^\lambda_{l}{\fg}, \ \li^\lambda_{\leq n}{\fg}=\bigoplus_{l\leq n}\li^\lambda_{l}{\fg}
\eeq
and for $k\in\Z/m\Z$
\beq
\li^\lambda_{n}{\fg}_k=\li^\lambda_{n}{\fg}\cap \fg_k\text{ etc.}
\eeq

\section{Supports of character sheaves}\label{sec:support-character-sheaves}

In this section, we define a Whitney stratification of $\fg_1$ by $G_0$-stable locally closed subsets, analogous to Lusztig's stratification of reductive groups introduced in~\cite[\S3.1]{L84}. It is well-known to experts that such Whitney stratifications exist in great generality, and in particular, in the setting of graded Lie algebras. As we could not find this precise statement in the literature, we include a discussion below. 

In~\autoref{ssec:support}, we define and discuss the supporting stratum for character sheaves and generic strata of the support of anti-orbital complexes. In~\autoref{ssec:restriction-sheaf}, we show that the restriction of an anti-orbital complex to a generic slice is a bi-orbital complex.

\subsection{A stratification of \texorpdfstring{$\fg_1$}{g1}}
\label{stratification}

Let $$f:\fg_1 \to \fg_1/\!/G_0$$ be the quotient map. Then $\fg_1/\!/G_0$ has a Whitney stratification, the Luna stratification. Furthermore, for any stratum $\cT$ in $\fg_1/\!/G_0$ the map $f^{-1}(\cT) \to \cT$ is a locally trivial fibration. These two 
statements hold for any action of a reductive group on a vector space and are proved using the Luna slice theorem. In our setting we refine the latter statement slightly and specify a Whitney stratification of $\fg_1$ making $f$ a stratified map so that $f^{-1}(\cT) \to \cT$ is a locally trivial stratified fibration. The character sheaves, and more generally the anti-orbital complexes, are constructible with respect to this stratification.

Let $\fc\subseteq\fg_1$ be a Cartan subspace, that is, a maximal commuting subspace consisting of semisimple elements.   We define a stratification on $\fc$ as follows. Let  $\Phi(\fg,\fc)\subseteq\fc^*$ denote the roots and consider the root hyperplanes on $\fc$. Two points in $\fc$ are in the same stratum if and only if they are in the same collection of root hyperplanes. 

Let $t\in\fc$. We write $\fc_{t}\subseteq\fc$ for the stratum that $t$ lies in and call it a $\fc$-stratum. Note that, for $t,t'\in\fc$, $\fc_t=\fc_{t'}$ if and only if $Z_G(t)=Z_G(t')$ if and only if $Z_{G_0}(t)=Z_{G_0}(t')$. In particular, $Z_G(\fc_t)=Z_G(t)$ for $t\in\fc$. It is clear that, for any $g\in N_G(\fc)$, $\Ad_g$ sends a $\fc$-stratum to a $\fc$-stratum.

Let $t\in\fc$. Recall the notation
 $\fZ_t=Z(Z_{\fg}(t))$ and $\fh_t=[Z_{\fg}(t),Z_{\fg}(t)]$. Note that $\fZ_{t,1}=\fZ_t\cap\fg_1$ is the (Zariski) closure of $\fc_t$, the intersection of all root hyperplanes in $\fc$ containing $\fc_t$. In particular $\fZ_t=Z(Z_{\fg}(\fZ_t))=Z(Z_{\fg}(\fZ_{t,1}))$. We have
 \beq\label{c-strata}
 \fc_t=\fZ_{t,1}^\circ:=\{z\in\fZ_{t,1}\mid Z_G(z)=Z_G(t)\}.
 \eeq
 According to~\cite{Vi77}, we have $\fg_1/\!/G_0\cong \fc/W_\fc$, where the Weyl group $W_\fc=N_{G_0}(\fc)/Z_{G_0}(\fc)$ is a complex reflection group. The strata $\{\fc_t\}$ induce the Luna stratification on $\fg_1/\!/G_0$. 

 Next, we construct the stratification of $\fg_1$, which refines the decomposition of $\fg_1$ into the pieces $f^{-1}(f(\fc_t))$. 
 To that end, let $n\in\fh_{t,1}=\fh_t\cap\fg_1$ be nilpotent. We define a \emph{stratum} 
 on $\fg_1$ to be 
 \beq
 \cS(t,n):=G_0({\fc_t}+n).
 \eeq

\begin{lemm}\label{lem:strata-conjugate} We have $\cS(t,n)=\cS(t',n')$ if and only if there exists $g\in G_0$ such that ${\fc_t}=\Ad_g{\fc_{t'}}$ and $n=\Ad_gn'$. Otherwise, they are disjoint.
\end{lemm}

\begin{proof} 

Suppose that $\cS({t},n)$ and $\cS({{t'}},n')$ have a non-empty intersection.  
Then there exist $x_s\in \fc_t,\, x_s'\in {\fc_{t'}}$ and $g\in G_0$ such that $\Ad_g(x_s+n)=x_s'+n'$. Since these are Jordan decompositions, we have $\Ad_gx_s=x_s'$ and $\Ad_gn=n'$.  
 By~\eqref{c-strata}, we conclude that $\Ad_g{\fc_t}=\Ad_g\fZ_{x_s,1}^\circ=\fZ_{x_s',1}^\circ=\fc_{t'}$.
 This completes the proof of the lemma.
\end{proof}

Consequently, the map
 \[
 G_0\times^{N_{G_0}(\fZ_t)\cap Z_{G_0}(n)}(\fc_t + n) \to \fg_1,\quad [g: x]\mapsto \Ad_g x
 \]
 is an immersion. Therefore, $\cS(t, n)$ is a connected locally closed subset of $\fg_1$ and non-singular when endowed with the  structure of reduced subscheme.

Let $x=x_s+x_n\in\fg_1$. Then there exists $g\in G_0$ such that $\Ad_gx_s\in\fc$. It follows that $x\in\cS(\Ad_gx_s,\Ad_{g}x_n)$. Thus, $\fg_1$ is a disjoint union of its strata.

For any $x\in\fg_1$, we will denote by $\cS_x$ the unique stratum containing $x$. It is clear from \autoref{lem:strata-conjugate} that $\cS_x$ is independent of the choice of the Cartan subspace $\fc$.

\begin{prop}\label{prop:whitney}
The family $\{\cS_x\}_{x\in \fg_1}$ of locally closed subsets forms a Whitney stratification of $\fg_1$. 
\end{prop}
\begin{proof}
Let $x = x_s + x_n\in \fg_1$. Given a stratum $\cS$ such that $x\in \overline\cS$, we need to show that $\cS_x\subseteq \overline\cS$ and that the pair $(\cS_x, \cS)$ satisfies Whitney's condition (b) in \cite[\S 19]{Wh65}. \par
Consider first the case where $x$ is nilpotent modulo centre, so that $x_s \in Z(\fg)_1$. Choose an $\fsl_2$-triplet $(e, h, f)$ in $\fg$ such that $e = x_n$, $h\in \fg_0$ and $f\in \fg_{-1}$. Averaging the Slodowy slice $\mathfrak{s} = x_n + Z_{\fg_1}(f)$ yields a map
\[
\kappa: Z(\fg)_1\times G_0\times \mathfrak{s} \to \fg_1,\quad \kappa(z, g, y) = z + \Ad_g y.
\]
This map is a smooth fibre bundle over the stratum $\cS_x = Z(\fg)_1\times G_0x$ with fibre $Z_{G_0}(x)$. It suffices thus to check the required properties for the pair $(\kappa^{-1}(\cS_x), \kappa^{-1}(\cS))$. Since 
\[
(\kappa^{-1}(\cS_x), \kappa^{-1}(\cS)) \cong Z(\fg)_1\times G_0 \times (\cS_x\cap \mathfrak{s}, \cS\cap \mathfrak{s}),
\]
it suffices to check the required properties for $(\cS_x\cap \mathfrak{s}, \cS\cap \mathfrak{s})$. As $\cS_x\cap \mathfrak{s} = \{x_n\}$ is reduced to a point, the Whitney's condition (b) for this pair follows from \cite[Lemma 19.3]{Wh65}. \par
Let us turn to the general case. Set $M = Z_G(x_s)$ and $\fm = \Lie M$. Define an open subset
\[
\fm_1^{\circ} = \{y=y_s+y_n\in \fm_1\;\mid\; Z_G(y_s)\subseteq M\}. 
\]
Every element of $\cS_x$ can be conjugated into $\fm_1^{\circ}$.  
Consider the following morphism:
\[
	G_0\times^{M_0}\fm_1^{\circ}\to \fg_1,\quad [g:y]\mapsto \Ad_g y.
\]
It is \'etale and thus open over $\cS_x$. Therefore, by $G_0$-conjugation, it suffices to show that the pair $(\cS_x\cap \fm_1^{\circ}, \cS\cap \fm_1^{\circ})$ satisfies the required conditions.  \par
We claim that $\cS_x\cap \fm_1^{\circ}$ and $\cS\cap \fm_1^{\circ}$ are disjoint unions of subvarieties of $\fm_1^{\circ}$ of the form $\cT \cap \fm_1^{\circ}$, where $\cT$ is a stratum of $\fm_1$ for the adjoint $M_0$-action. We will only prove this for $\cS\cap \fm_1^{\circ}$ since the proof in the other case is similar. From \autoref{lem:strata-conjugate}, we see that $\cS\cap \fm_1^{\circ}$ is a union of such subvarieties, and it remains to show that the union is disjoint. Let $y\in \cS\cap \fm_1^{\circ}$ and let $\cS^M_y$ denote the stratum of $\fm_1$ containing $y$. Then the inclusion $\cS^M_y\hookrightarrow \cS\cap \fm_1^{\circ}$ induces a linear map on the tangent spaces:
\[
T_y\cS^M_y \to T_y(\cS\cap \fm_1).
\]
It is an isomorphism: indeed, since $y\in \fm_1^{\circ}$, we have 
\[
T_y(\cS\cap \fm_1) = ([\fg_0, y_s]\oplus  \fZ_{y, 1}\oplus [Z_{\fg_0}(y_s), y_n])\cap \fm_1= [\fm_0, y_s]\oplus \fZ_{y, 1} \oplus [Z_{\fm_0}(y_s), y_n]   = T_y\cS^M_y.
\]
As this isomorphism holds for every $y\in \cS\cap \fm_1^{\circ}$, it follows that the union in question is disjoint. \par

 Since the elements of $\cS_x\cap \fm_1^{\circ}$ are nilpotent modulo centre in $\fm$, the situation is reduced to the previous case. It follows that $(\cS_x\cap \fm_1^{\circ}, \cS\cap \fm_1^{\circ})$ satisfies Whitney's condition (b). 
\end{proof}

\begin{prop}\label{prop:fibration} For any stratum $\cT$ of $\fg_1 /\!/ G_0$, the pre-image $f^{-1}(\cT)$ is a union of strata of $\fg_1$ and the restriction $f|_{f^{-1}(\cT)}$ is an \'etale-locally trivial fibration over $\cT$. 
\end{prop}
\begin{proof}
Choose any Cartan subspace $\fc\subseteq \fg_1$ and let $t\in \fc$ be such that $f(t)\in \cT$. Set $M = Z_G(t)$. Then, the restriction of $f:\fg_1\to \fg_1/\!/G_0$ to $\fc_t$ yields a finite \'etale cover of $\cT$. It suffices to check the statements after the \'etale base change along $f|_{\fc_t}:\fc_t\to \cT$. We have
\[
\fg_1 \times_{\fg_1/\!/G_0}\fc_t \cong \fc_t\times G_0\times^{M_0}\fm_1^{\nil}. 
\]
For any $N_{G_0}(\fm)$-orbit $Z$ in $\fm_1^{\nil}$, the subset $\fc_t\times G_0\times^{M_0}Z$ is the pre-image of a stratum of $\fg_1$ contained in $f^{-1}(\cT)$. Therefore, the stratification on $\fg_1\times_{\fg_1/\!/G_0}\fc_t$ induced from $\fg_1$ is trivial over $\fc_t$, which proves the proposition.  
\end{proof}

\subsection{Supporting stratum of a character sheaf}\label{ssec:support}

By definition, the singular support of a character sheaf on $\fg_1$ is contained in the union of the conormal bundles $T^*_\CO\fg_{-1}$ where $\CO$ runs through the nilpotent $G_0$-orbits on $\fg_{-1}^{\nil}$. For such an orbit we have
\beq
\label{conormal}
 T^*_\CO\fg_{-1} =\{(e,x)\in \fg_{-1}\times\fg_1\mid e\in\CO, \ \ [e,x] = 0\} \,.
\eeq
 We write 
 \beqn
\tilde\cO=\{x\in \fg_1\mid\exists\, e\in \cO,\; [e,x] = 0\}\subseteq\fg_1
\eeqn
for its projection in $\fg_1$. We have
\begin{lemm}
The  $\tilde\cO$ is a union of strata.
\end{lemm}
\begin{proof}
Let $x=x_s+x_n\in\tilde\cO$. We have to show that $\cS_x\subseteq\tilde\cO$. We fix a Cartan subspace $\fc\subseteq \fg_1$ containing $x_s$. 

Given any $x' = x'_s + x'_n\in \cS_x$, there exists $g\in G_0$ such that $\Ad_g x'_s\in {\fc_{x_s}}$ and $\Ad_g x'_n = x_n$. Pick $e\in \cO$ such that $[e,x] = 0$. Then $[e,x_s] = [e,x_n] = 0$. Since $e\in Z_\fg(x_s)=Z_\fg(\fc_{x_s})$, it follows that 
\[
[\Ad_{g^{-1}}e,x'] = \Ad_{g}^{-1}([e,\Ad_{g}x'_s] + [e,x_n]) = 0.
\]
Therefore, $x'\in \tilde\cO$.
  \end{proof}   
The closure of the constructible set $\tilde\cO$ is irreducible and hence contains a unique stratum with the same closure, which we denote by $\widecheck\CO$. It follows immediately that $\overline{T^*_{\widecheck\CO}\fg_1} = \overline{T^*_\CO\fg_{-1}}$.  Hence the singular support of any character sheaf lies in the closure of the union of the conormal bundles $T^*_{\widecheck\CO}\fg_{1}$ where $\CO$ runs through the nilpotent $G_0$-orbits on $\fg_{-1}^{\nil}$.

From the discussion above and \autoref{prop:whitney}, we conclude:
\begin{prop}\label{prop:stratum-support} Every anti-orbital complex on $\fg_1$ is constructible with respect to the stratification $\{\cS_x\}_{x\in \fg_1}$. Moreover, if $\cF$ is a character sheaf on $\fg_1$, then there exist a nilpotent $G_0$-orbit $\CO$ in $\fg_{-1}$ and a $G_0$-equivariant local system $\cF^{\circ}$ on $\widecheck\CO$ such that $\cF\cong\IC(\cF^{\circ})$.
\end{prop}

We will call $\widecheck\CO$ the \emph{supporting stratum} of $\cF$.

\subsection{Restriction along a generic slice}\label{ssec:restriction-sheaf}

\begin{defi}\label{def-generic strata} Let $Z\subseteq \fg_1$ be a Zariski-closed subset  which is a union of strata of $\fg_1$.  A stratum $\cS\subseteq Z$ is called \emph{generic} in $Z$ if the image of $\cS$ in $Z /\!/ G_0$ is open and $\cS$ is maximal among the strata contained in $Z$.
\end{defi}\par

Let $\cG\in \Dbn_{G_0}(\fg_1)$ and let $\cS\subseteq \supp(\cG)$ be a generic stratum. Choose $x=x_s+x_n\in\cS$ and  write  $M = Z_G(x_s)$, $\Lm=\on{Lie}M$. Recall the notation $\fh_x = [\fm,\fm]$ (see~\eqref{eqn-hz}).

We consider the restriction $\cG|_{x_s + \fh^{\nil}_{x,1}}$ and regard it as an object of $\Db_{Z_{G_0}(x_s)}(\fh_{x,1})$ via translation.

\begin{prop} \label{prop:localtriv} The locus $x_s+\fh_{x,1}$ is non-characteristic for $\CG$ in an open neighbourhood of $x_s+\fh^{\nil}_{x,1}$.
\end{prop}

\begin{proof} 
We claim that there exists an open subset $X\subseteq x_s + \fh_{x, 1}$ containing $x_s + \fh_{x, 1}^{\nil}$ such that $X\cap \supp(\CG)\subseteq x_s+\fh^{\nil}_{x,1}$.   

Pick a Cartan subspace $\fc\subseteq \fg_1$ containing $x_s$ and  set $U = \{y\in \fc\;|\; Z_G(y)\subseteq M\}$. By definition, $U$ is a union of $\fc$-strata and is an open neighbourhood of $\fc_{x_s}$ in $\fc$. Let $f: \fm_{1}\to \fm_{1} /\!/ M_{0}$ denote the adjoint quotient map. The image $f(U)$ is open in $\fm_1 /\!/ M_0$ because the restriction $f|_{\fc}:\fc \to \fm_1 /\!/ M_0$ is a ramified cover. We define $X = (x_s + \fh_{x, 1})\cap f^{-1}(f(U))$. It is an $M_0$-invariant open subset of $x_s + \fh_{x, 1}$ containing $x_s + \fh_{x, 1}^{\nil}$.  \par
Let $y\in \supp(\cG) \cap X$. We shall prove that $y\in x_s + \fh^{\nil}_{x,1}$. Since $x_s + \fh^{\nil}_{x,1}$ is $M_0$-stable, up to replacing $y$ with an $M_0$-conjugate, we may assume that the semisimple part $y_s$ lies in $(x_s + \fh_{x,1})\cap U$. If $y_s\notin \fc_{x_s}$, then $y_s$ is contained in a $\fc$-stratum strictly bigger than $\fc_{x_s}$, which contradicts the genericity hypothesis of $\cS$. It follows that $y_s \in \fc_{x_s}$, which in turn  implies $y_s\in \fc_{x_s}\cap (x_s + \fh_{x,1}) = \{x_s\}$. Therefore, $X$ satisfies the required condition.

To prove the non-characteristic property, it remains to show that for $(y,\xi)\in T^*_{x_s+\fh_{x,1}}\fg_1\cap \operatorname{SS}(\CG)$ with $y\in x_s+\fh^{\nil}_{x,1}$, we have $\xi=0$. As $(y,\xi)\in T^*_{x_s+\fh_{x,1}}\fg_1$, we conclude that $\xi|_{\fh_{x,1}}=0$. As $\supp(\cG)$ is contained in the union of the conormal bundles $T^*_{\cO}\fg_{-1}$ for nilpotent $G_0$-orbits $\cO\subseteq\fg_{-1}$, we conclude from~\eqref{conormal} that  $[\xi,y]=0$. Because $x_s$ is the semisimple part of $y$, we also have $[\xi,x_s]=0$. As $\xi\in\fg_{-1}^{\nil}$, we must have $\xi\in\fh_{x,-1}$. Recalling that $\xi|_{\fh_{x,1}}=0$, we conclude that $\xi=0$.\end{proof}

    Thus, we obtain by  \cite[Proposition 5.4.13 \& Corollary 10.3.16(ii)]{KS90}, for example:
\begin{coro}  \label{coro:restriction-biorbital}
The singular support $\operatorname{SS}(\cG|_{x_s + \fh^{\nil}_{x,1}})\subseteq (x_s+\fh_{x,1}^{\nil})\times\fh_{x,-1}^{\nil}$, i.e., $\cG|_{x_s + \fh^{\nil}_{x,1}}$ is a $Z_{G_0}(x_s)$-equivariant bi-orbital complex on $\fh_{x,1}$. Furthermore, if $\cG$ is perverse,  then $\cG|_{x_s + \fh^{\nil}_{x,1}}$ is concentrated in perverse cohomological degree $-\codim_{\fg_1}\fh_{x,1}$. 
\end{coro}

\begin{rema}\label{rem:restriction}
Let $\cG\in \Dbn_{G_0}(\fg_1)$ and let $\cS\subseteq \supp \cG$ be a generic stratum and set $\tilde\cS = f^{-1}(f(\cS))$. Then, we have for $x\in \cS$:
\[
    F_x:= f^{-1}(f(x)) \cong G_0\times^{Z_{G_0}(x_s)}(x_s + \fh^{\nil}_{x,1}).
\]
Since $\tilde \cS\to f(\cS)$ is a stratified submersion by \autoref{prop:fibration} and $\cG$ is constructible along the stratification on $\fg_1$, the restriction $\cG|_{\tilde \cS}$ is locally acyclic over $f(\cS)$. Therefore, \autoref{coro:restriction-biorbital} can be viewed as the restriction of $\cG|_{\tilde\cS}$ to the fibre $F_x$, which is transversal to every stratum of $\tilde\cS$. 
\end{rema}

\section{Main theorems and key propositions}\label{sec:main}
We state our main theorems and some key propositions  in this section. In particular, in our main theorem (\autoref{thm:main}) we show  that all cuspidal character sheaves arise from Fourier--Sato transforms of nearby-cycle sheaves associated with  nil-supercuspidal data (see~\autoref{def-sc/nsc}) and that the latter sheaves are IC-extensions of local systems. Furthermore, we show that non-equivalent nil-supercuspidal data (see~\autoref{def:equiv}) give rise to mutually orthogonal nearby-cycle sheaves. The proof of \autoref{thm:main} is given in \autoref{ssec:proof of main thm} using results proved in later sections.

We show that the nearby-cycle sheaves associated with  nil-supercuspidal data are both projective and injective in $\Perv_{G_0}(\fg^{\nil}_1)$, see~\autoref{thm:vanishing}. We make use of a key technical result: the adjunction formula~(\autoref{thm-adjunction}). Moreover, we show that the supporting stratum of a cuspidal character sheaf consists of $G_0$-distinguished elements and  determine the dimension of its image under the adjoint quotient map in~\autoref{thm:cuspidal-distinguished}.  In \autoref{ssec:stable} we deduce some consequences for (GIT)-stable gradings.

It will be convenient for us to pass between sheaves on $\fg_1$ and $\fg_{-1}$. Therefore, we fix a Chevalley involution $\tau_0:G_0\to G_0$ which restricts to $t\mapsto t^{-1}$ on a fixed maximal torus $T_0\subseteq G_0$, and a linear isomorphism $\tau:\fg_1 \xrightarrow{\sim} \fg_{-1}$, such that $(\tau_0, \tau)$ intertwines the adjoint actions of $G_0$ on $\fg_1$ and $\fg_{-1}$.

\subsection{The functor \texorpdfstring{$\dag$}{Lg}} \label{ssec:functor dag}
 
For technical reasons, we introduce the following contravariant functor:
\beq\label{eqn-dag}
 \dag:\Db_{G_0}(\fg^{\nil}_1)^{\op}\xrightarrow{\sim} \Dbn_{G_0}(\fg_1),\quad \cF \mapsto \cF^{\dag} := \operatorname{FS}(\tau_* \D \cF),
\eeq
 where $\operatorname{FS}:\Db_{G_0}(\fg^{\nil}_{-1})\to \Dbn_{G_0}(\fg_1)$ is the Fourier--Sato transform (see \autoref{subsec:orbital}) and $\D$ is the Verdier duality on $\Db_{G_0}(\fg^{\nil}_1)$. It is a perverse-t-exact triangle equivalence. 

 \begin{lemm}\label{lem:parab-dag}
Let $P\subseteq G$ be a $\theta$-stable parabolic subgroup and $L\subseteq P$ a $\theta$-stable Levi factor containing $T_0$ (fixed above). Let $P'\subseteq G$ be the parabolic subgroup opposite to $P$ with common Levi factor $L$. Then there are isomorphisms of functors
\[
\dag\circ \Ind^{\fg_1}_{\fp_1} \cong  \Ind^{\fg_1}_{\fp'_1}\circ\dag: \Db_{L_0}(\fl^{\nil}_1)^{\op}\xrightarrow{\sim} \Dbn_{G_0}(\fg_1)
\]
\[
\dag\circ \Res^{\fg_1}_{\fp_1} \cong  \Res^{\fg_1}_{\fp_1}\circ\dag[D]: \Db_{G_0}(\fg_1^{\nil})^{\op}\xrightarrow{\sim} \Dbn_{L_0}(\fl_1)\,
\]
where $D = \dim \fg_1 - \dim \fg_0 - \dim \fl_1 + \dim \fl_0$.
\end{lemm}
\begin{proof}
By \cite[\S 1.6]{Liu24}, we can choose a cocharacter $\lambda: \Cc\to T_0$  such that 
    \[
        \fp = \li^\lambda_{\ge 0}{\fg},\quad \fl = \li^\lambda_{0}{\fg},\quad \fp' = \li^\lambda_{\le 0}{\fg}.
    \]
    Since $\tau_0$ acts by inversion on $T_0$, we have $\tau_0\circ\lambda = -\lambda$. Therefore, $\tau(\fp_1) = \li^{-\lambda}_{\ge 0}{\fg}_{-1} = \li^{\lambda}_{\le 0}{\fg}_{-1} = \fp'_{-1}$ and similarly $\tau_0(P_0) = P'_0$, $\tau(\fl_1) = \fl_{-1}$ and $\tau_0(L_0) = L_0$. It follows that $\tau_* \circ \Ind^{\fg_1}_{\fp_1} \cong  \Ind^{\fg_{-1}}_{\fp'_{-1}} \circ\tau_*$. The first isomorphism follows because the Fourier--Sato transform commutes with parabolic induction by \eqref{eqn-parabolic-FS} and the parabolic induction is normalised so that it commutes with Verdier duality. The second isomorphism is proved in an analogous manner making use of~\autoref{prop:hyperbolic-parab}.
\end{proof}

\par

\subsection{The nearby cycles}Let $\rc\subseteq\fg_1$ be a $G_0$-orbit and $\chi$ a $G_0$-equivariant local system on $\rc$.  The nearby-cycle construction, which we introduce now, deforms the perverse sheaf $\IC(\chi)$ to the nilpotent cone $\fg_1^\nil$. It is a straightforward generalisation of the construction in~\cite{Grin98,GVX23,GVX20}. \par
 
Recall the quotient map $$f:\fg_{ 1}\to \fg_{ 1} /\!/G_{0}.$$ We let $x\in\rc$ and write
\beq\label{eqn-cz}
\cZ = \{(z, r)\in \fg_{1}\times \C\;|\; f(z) = rf(x)\}
\eeq
with $(G_{0}\times \Cc)$-action given by $$(g, t)(z, r) = (t^{-1}\Ad_g z, t^{-1}r).$$ 
Let   
$$\pi:\cZ\to \C$$ be the projection onto the second component. Then, we have $\pi^{-1}(0) \cong \fg_{1}^{\nil}$.
Set $\cU = \pi^{-1}(\Cc)$. The $\Cc$-action on $\cU$ is free and induces an isomorphism
\[
     \Cc\times \pi^{-1}(1)\xrightarrow{\sim} \cU,\quad (t, z)\mapsto (tz, t).
\]
Consider the $G_{0}$-equivariant perverse sheaf $\IC(\chi)\in \Perv_{G_0}(\pi^{-1}(1))$. It extends to a $(G_{0}\times \Cc)$-equivariant perverse sheaf on $\cU$, denoted by $\cL_{\chi}$. 
We form the nearby-cycle sheaf
\beq\label{eqn-psi}
\cP_{\chi}:=\psi_{\pi}\cL_\chi.
\eeq
 By the t-exactness of $\psi_\pi$ (see~\cite[\S 10.3.13]{KS90} for example), $\cP_{\chi}$ is a $G_0$-equivariant perverse sheaf on $\fg^{\nil}_1 = \pi^{-1}(0)$. 

 \subsection{Dependence on input data \texorpdfstring{$(\rc,\chi)$}{<PDF>}}

In this subsection we discuss how the nearby cycle construction depends on the input data $(\rc,\chi)$. We may and will assume that $\chi$ is irreducible. Here we consider general pairs. In the next subsection we specialize to the pairs used in the rest of the paper. 

Let $\rc$ and $\rc'$ be two $G_0$ orbits on $\fg_1$ and let $\chi$ and $\chi'$ be irreduicble $G_0$-equivariant local systems on $\rc$ and $\rc'$, respectively. We assume that $\rc$ and $\rc'$ are contained in the same stratum $\cS=G_0\cdot(\fc_s+n)$; here $s\in \fg_1$ is semisimple and $n\in Z_{\fg_1}(s)^{\nil}$. We want to understand the relationship between  $\cP_{\chi}$ and $ \cP_{\chi'}$. 

We set 
$$
\fd=\overline{\fc_s}\text{ and } W_\fd=N_{G_0}(\fd)/Z_{G_0}(\fd).
$$  
Choose $x = x_s + n\in \rc$ and $x' = x'_s + n\in \rc'$ such that $x_s,x'_s\in\fc_s$. Then, we have $Z_G(x_s) = Z_G(x'_s)$ and $Z_{G_0}(x) = Z_{G_0}(x')$. Therefore, we can identify $\pi_0(Z_{G_0}(x))^{\wedge}$, the set of irreducible representations of $\pi_0(Z_{G_0}(x))$, with the analogous $\pi_0(Z_{G_0}(x'))^{\wedge}$.

Consider the fibration $\cS \to \fd^{\circ}/W_{\fd}$, $\fd^\circ=\fc_s$, 
and its fibre $G_0 x$ through $x$. Let us write $\tilde\cS \to \fd^{\circ}$ for the fibration under base change. Note that this fibration can be trivialised by using the section $\fd^{\circ} \to \fd^{\circ} + n$. 
We have the following diagram
\[
\begin{tikzcd}
    1 \ar{r}& \pi_0(Z_{G_0}(x)) \ar{r}\ar{d}{=}& \pi_1^{G_0}(\tilde \cS, x) \ar{r}\ar{d}& \pi_1(\fd^{\circ}, s) \ar{r}\ar{d}&  1 \\
    1 \ar{r}& \pi_0(Z_{G_0}(x)) \ar{r}& \pi_1^{G_0}(\cS, x) \ar{r}& B_{\fd} \ar{r} &  1,
\end{tikzcd}
\]
where $B_{\fd} = \pi_1(\fd^{\circ}/ W_{\fd}, s)$, which is by definition the braid group of $W_{\fd}$.
The top row is exact because $\tilde\cS \to \fd^{\circ}$ is split. The second and third vertical morphisms are inclusions of normal subgroups because they come from Galois covers with Galois group $W_{\fd}$. It follows that the bottom row is also exact.  \par

The conjugation action of $\pi_1^{G_0}(\cS, x)$ on $\pi_0(Z_{G_0}(x))$ induces an action of $B_{\fd}$ on $\pi_0(Z_{G_0}(x))^{\wedge}$. 
Because $\tilde\cS \to \fd^{\circ}$ is split, this action factors through $W_{\fd}$. \par

\begin{defi}\label{def:equiv}
We say that two pairs $(\rc, \chi)$ and $(\rc', \chi')$ are \emph{equivalent} and write $(\rc, \chi)\sim (\rc', \chi')$ if 
\begin{enumerate}
    \item $\rc$ and $\rc'$ are contained in the same stratum $\cS$
    \item $\chi$ and $\chi'$, viewed as elements of $\pi_0(Z_{G_0}(x))^{\wedge}$ ($x\in\rc$),  are $W_{\fd}$-conjugate.
\end{enumerate} 
\end{defi}
We have
\begin{theo}
    If the pairs $(\rc, \chi)$ and $(\rc', \chi')$ are equivalent then $\cP_{\chi}$ is isomorphic to $ \cP_{\chi'}$.
\end{theo}
    
    This theorem can be proved using arguments similar to~\cite[Section 3.2]{GVX23}. We will not do so here as the arguments are very different from those used in this paper.

    We will prove a special case we will be using below as well as an orthogonality statement for the nil-supercuspidal data we introduce below.

\subsection{The supercuspidal and nil-supercuspidal data} To state our main results, we introduce the following definition:

 \begin{defi}\label{def-sc/nsc}
     Let $\rc\subseteq\fg_1$ be a $G_0$-orbit  and $\chi$ an irreducible $G_0$-equivariant local system on  $\rc$.  
\begin{enumerate}
\item[(i)]     We say that $(\rc,\chi)$ is a supercuspidal datum if the orbit $\rc$ is $G_0$-distinguished, and $\chi$ is clean, that is, $(j_\rc)_!\chi\cong (j_{\rc})_*\chi$, where $j_\rc:\rc\to\fg_1$.

\item[(ii)] We say that $(\rc,\chi)$ is a nil-supercuspidal datum if it is supercuspidal and \linebreak$\IC(\chi|_{\rc\cap(x_s+\fh_{x,1}^{\nil})})\in \on{Perv}_{H_{x,0}}(\fh_{x,1}^{\nil})$ is a sum of bi-orbital cuspidal sheaves for any $x\in \rc$.
\end{enumerate}
     \end{defi}

     \begin{rema}
     Note that in the definition (ii) of nil-supercuspidal datum the cleaness assumption is unnecessary as it follows from \autoref{characterisation of supercuspidal}. 
 \end{rema}

 We have the following vanishing property:
 \begin{prop}\label{prop:supercuspidal}
     Let $(\rc,\chi)$ be a supercuspidal datum. Then $\on{IC}(\chi)$ is annihilated by parabolic restriction along every $\theta$-stable proper parabolic subgroup of $G$.
 \end{prop}
The bi-orbital condition for nil-supercuspidal data can be replaced with the following simple numerical condition:
 \begin{prop}\label{cor:dim est-biorbitalcuspidal}
        Let $(\rc,\chi)$ be a supercuspidal datum. Let $\cS$ be the stratum containing $\rc$. Then, we have $$\dim\cS-\dim\rc \le \dim\fg_1-\dim\fg_0+\dim Z(\fg)_0.$$ Moreover, the equality holds if and only if $(\rc, \chi)$ is a nil-supercuspidal datum.
    \end{prop}
 
The proofs of~\autoref{prop:supercuspidal} and~\autoref{cor:dim est-biorbitalcuspidal} will be given in~\autoref{ssec:proof-supercuspidal} and~\autoref{ssec:bi-orbital}, respectively.

\subsection{The main theorem}We can now state our main result:

\begin{theo}\label{thm:main}\leavevmode 
\begin{enumerate}
    \item[{\rm(i)}] Let $(\rc,\chi)$ be a supercuspidal datum. Then the irreducible subsheaves of  $\cP_{\chi}^\dag$ are cuspidal character sheaves. 
    \item[{\rm(ii)}]  Let $\cF\in\on{Char}_{G_0}(\fg_1)$ be a cuspidal character sheaf. Then $\cF$ arises as a subsheaf of $\cP_{\chi}^\dag$ for some nil-supercuspidal datum $(\rc, \chi)$.
 More precisely, let $\widecheck{\cO}$ denote the supporting stratum of $\cF$ (see \autoref{prop:stratum-support}) and write $\cF =\IC(\cF^{\circ})$ for some $G_0$-equivariant local system on $\widecheck\cO$. Then, $(\rc, \chi)$ can be chosen so that $\rc\subseteq\widecheck\cO$ and $\chi$ is any simple summand of $\cF^{\circ}|_{\rc}$. 

 \item[{\rm(iii)}] Let $(\rc, \chi)$ be a nil-supercuspidal datum and let $\cS$ denote the stratum containing $\rc$. Then, there exists a $G_0$-equivariant local system $\cM_{\chi}$ on $\cS$ such that $\cP^{\dag}_{\chi} \cong \IC(\cM_{\chi})$. In particular, the irreducible subsheaves of $\cP^{\dag}_{\chi}$ are supported on $\cS$.

 \item[{\rm(iv)}] Let $(\rc, \chi)$ and $(\rc', \chi')$ be nil-supercuspidal data on $\fg_1$. If $(\rc, \chi)$ and $(\rc', \chi')$ are equivalent, then $\cP_{\chi}$ and $\cP_{\chi'}$ are isomorphic. Otherwise, we have $\Hom(\cP_{\chi}, \cP_{\chi'}[k]) = 0$ for every $k\in \Z$. 

 \end{enumerate}

\end{theo}
\begin{rema}
    The Verdier duality yields $\D\cP_{\chi} \cong \cP_{\chi^*}$. Therefore, the statements of \autoref{thm:main} remain true with ``quotient sheaf" replacing ``subsheaf".   
\end{rema}
An immediate consequence of our theorem is the following corrollary.
\begin{coro}[Springer-type correspondence]\label{coro:springer} There is a natural one-to-one correspondence 
\[
\{\cF\in \Char_{G_0}(\fg_1)\mid \text{$\cF$ cuspidal}\,\}\leftrightarrow\bigsqcup_{(\rc, \chi)/\sim} \Irr\left(\End(\cM_{\chi})\mof\right),
\]
given by $\cF\mapsto \Hom(\cP_{\chi}^{\dag}, \cF)$, for some (any) nil-supercuspidal datum $(\rc, \chi)$ such that\linebreak $\Hom(\cP_{\chi}^{\dag}, \cF)\neq 0$.
\end{coro}
 \begin{rema} 

 \autoref{coro:springer} 
 reduces the determination of cuspidal character sheaves on a given graded Lie algebra to two subproblems: 
\begin{enumerate}
    \item classification of equivalent classes of nil-supercuspidal data
    \item calculation of $\End(\cP_{\chi})$ and its irreducible modules, or equivalently, the calculation of irreducible direct summands of the cosocle of $\on{IC}(\cM_\chi)$.
\end{enumerate} 

The first problem is solved in \cite{LVX}. 
The second problem is solved in a large number of cases in~\cite{VX23st,VX21} with only the spin groups and some exceptional groups remaining.
\end{rema}

\subsection{Structure of the proof of \autoref{thm:main}}\label{ssec:proof of main thm}
In this subsection, we reduce the proof of \autoref{thm:main} to results proved in the later sections. 
We begin with the following vanishing property of the nearby-cycle sheaf, whose proof we defer to~\autoref{sec:parabolic-vanishing}:
 \begin{prop}\label{prop:parabolic-vanishing}
  
 Let $(\rc,\chi)$ be a supercuspidal datum. 
 Then, the nearby-cycle sheaf $\cP_{\chi}$ is annihilated by parabolic restriction along every $\theta$-stable proper parabolic subgroup.
 \end{prop}
By the first adjunction in~\eqref{eqn-adjunctions-parabolic}, this proposition implies that every simple quotient of $\cP_{\chi}$ is cuspidal. Since~\autoref{lem:parab-dag} implies that the functor $\dag$ respects cuspidality,  part (i) of \autoref{thm:main} follows, i.e., that the simple subsheaves of $\cP_{\chi}^{\dag}$ are cuspidal.

\subsubsection{Part (ii) of~\autoref{thm:main}}To address part (ii) of \autoref{thm:main}, we make use of the result below, whose proof we defer to~\autoref{sec:support}. Recall  that $x_s$ and $x_n$ denote the semisimple part and nilpotent part of $x$, respectively, and $\fh_{x,1}=[Z_\fg(x_s),Z_\fg(x_s)]\cap\fg_1$.
\begin{theo}\label{thm:cuspidal-distinguished}
 Let $\cF\in\on{Char}_{G_0}(\fg_1)$ be a cuspidal character sheaf with supporting stratum $\widecheck{\cO}$ (see \autoref{prop:stratum-support}). 
 \begin{enumerate}
     \item[{\rm(i)}] The stratum $\widecheck{\cO}$ 
 consists of $G_0$-distinguished elements and 
 \beqn
 \dim f(\widecheck\cO) = \dim \fg_1 - \dim \fg_0 + \dim Z(\fg)_0.
 \eeqn
 \item[{\rm(ii)}] The restriction $\CF|_{x_s + \fh^{\nil}_{x,1}}[-d]$ is a direct sum of bi-orbital supercuspidal sheaves on~$\fh_{x,1}$, where $x=x_s+x_n\in\widecheck{\cO}$ and $d= \on{codim}_{\fg_1}\fh_{x,1}$.
 
 \end{enumerate}
 \end{theo}
 We record the following immediate but useful consequence of~\autoref{thm:cuspidal-distinguished} for future use:
 \begin{coro}
 If $\dim \fg_1 - \dim \fg_0 + \dim Z(\fg)_0 < 0$, then there is no cuspidal character sheaf on $\fg_1$. If $\dim \fg_1 - \dim \fg_0 + \dim Z(\fg)_0 = 0$, then every cuspidal character sheaf on $\fg_1$ is bi-orbital and supercuspidal. \qed
 \end{coro}

 As a precursor of~\autoref{thm:cuspidal-distinguished}, we will prove in~\autoref{ssec:supercuspidal dis} the following proposition which applies to  anti-orbital perverse sheaves annihilated by all proper parabolic restrictions, in particular, to $\cP_\chi^\dag$. 
 
\begin{prop}\label{prop:annihlation-distinguished}
Let $\cF\in \Perv_{G_0}^{\nil}(\fg_1)$. Suppose that $\cF$ is annihilated by parabolic restriction along every $\theta$-stable proper parabolic subgroup of $G$. Then, given a generic stratum (see \autoref{def-generic strata}) $\cS\subseteq \supp \cF$ the following statements hold:
\begin{enumerate}
    \item[{\rm(i)}] The stratum $\cS$ consists of $G_0$-distinguished elements and
    $$\dim f(\cS) = \dim \fg_1 - \dim \fg_0 + \dim Z(\fg)_0$$ 
    
    \item[{\rm(ii)}]
    The restriction $\cF|_{\rc}[-\dim \cS]$ is a clean local system for any $G_0$-orbit $\rc\subseteq \cS$.
\end{enumerate}
\end{prop}

The second ingredient for the proof of part (ii) of~\autoref{thm:main} is the following important formula, whose proof we defer to~\autoref{sec:proof-adjunction}:

\begin{theo}[Adjunction formula]\label{thm-adjunction}
 Suppose that $Z(G)_0 = 1$. Let $(\rc,\chi)$ be a nil-supercuspidal datum. Then
\[
  \dim\Hom\left(\cP_{\chi} ,\cG[k]\right)=\dim\Hom\left(\cG^{\dag}, \IC(\chi)[N + k]\right)
\]
for every $\cG\in \Irr\Perv_{G_0}(\fg^{\nil}_1)$ and every $k\in \Z$, where $N=\dim\fg_1-\dim\fg_0$. Moreover, both sides vanish for $k\neq 0$.
 
\end{theo}

We will now deduce part (ii) of \autoref{thm:main} from these ingredients. 

Observe first that we may replace $G$ with $G / Z(G)_0$ in the statement without loss of generality. Indeed, the change of equivariance along the quotient map $G\to G/Z(G)_0$ yields an equivalence $\Perv_{G_0 / Z(G)_0}(\fg_1)\cong\Perv_{G_0}(\fg_1)$, which commutes with the nearby-cycle construction as well as the parabolic induction/restriction/co-restriction and the functor $\dag$. Consequently, we may assume that $Z(G)_0 = 1$. 

Let $\cF$ be a cuspidal character sheaf and $\widecheck\CO$ its supporting stratum. It follows from part (i) of \autoref{thm:cuspidal-distinguished} that $\widecheck\CO$ consists of $G_0$-distinguished elements.  
We choose $x=x_s+x_n\in\widecheck\CO$. We consider the fibre $F_x=f^{-1}(f(x))$ and the orbit $\rc:=G_0x\subseteq F_x$. We first observe that the isomorphism in \autoref{rem:restriction} yields a t-exact triangle equivalence
\beq
\label{restriction}
\Db_{G_0}(F_x) \cong \Db_{Z_{G_0}(x_s)}(x_s + \fh^{\nil}_{x,1}),\quad \cK\mapsto \cK|_{x_s + \fh^{\nil}_{x,1}}[N-d]
\eeq
 (since $\dim G_0/Z_{G_0}(x_s) = \dim \fg_0 - \dim \fh_0 = d - N$). By part (ii) of \autoref{thm:cuspidal-distinguished},  $\cF|_{x_s + \fh^{\nil}_{x,1}}[-d]$ is a direct sum of bi-orbital supercuspidal sheaves. In particular, these summands are clean. The t-exactness of~\eqref{restriction} implies that $\cF|_{F_x}[-N]$ is a direct sum of simple clean sheaves and therefore
 $\cF|_{\rc}[-N-\dim\rc]$ 
 is a direct sum of clean irreducible $G_0$-equivariant local systems concentrated in cohomological degree $0$. Let us choose an irreducible summand $\chi$ of $\cF|_{\rc}[-N-\dim\rc]$. By construction, $(\rc,\chi)$ is a nil-supercuspidal datum. We will now form $\cP_{\chi}$. Let $\cG\in\Irr\Perv_{G_0}(\fg^{\nil}_1)$ be such that $\cG^{\dag}=\cF$. Applying \autoref{thm-adjunction} we obtain that
\[\dim\Hom\left(\cF, \cP_{\chi}^\dag\right)=\dim\Hom(\cP_{\chi},\cG)\]\[=\dim\Hom(\cF,\IC(\chi)[N])=\dim\Hom(\cF|_\rc,\chi[N+\dim\rc])\neq 0\,.\]
We conclude that there is a non-trivial morphism $\cF\to \cP_{\chi}^\dag$. \qed

We record the following corollary of \autoref{thm-adjunction} as an interesting by-product:
\begin{coro}\label{thm:vanishing}
Suppose that $Z(G)_0 = 1$. Let $(\rc,\chi)$ be a nil-supercuspidal datum.  
We have
\[
\Hom(\cP_{\chi}, \cG[k]) = 0\text{ and } \Hom(\cG, \cP_{\chi}[k]) = 0, \text{ for $k> 0$ and for  every $\cG\in \Perv_{G_0}(\fg^{\nil}_1)$}.
\]
 In particular, $\cP_{\chi}$ is projective and injective in $\Perv_{G_0}(\fg^{\nil}_1)$. 
\end{coro}
\begin{proof}
Let $\cG\in \Perv_{G_0}(\fg^{\nil}_1)$. Since Verdier duality commutes with nearby cycles, we have
\[
    \Hom(\cP_{\chi}, \cG[k]) \cong \Hom(\D\cG,  \D\cP_{\chi}[k]) \cong \Hom(\D\cG, \cP_{\chi^*}[k]).
\]
Thus, it suffices to prove the first vanishing property. When $\cG$ is simple, this follows from the second assertion of~\autoref{thm-adjunction}. The general case can be deduced by d\'evissage. 
\end{proof}

\subsubsection{Part (iii) of~\autoref{thm:main}}{}To address part (iii) of \autoref{thm:main}, we make use of the following:
\begin{coro}\label{lemm:supp-stratum-cusp}
Let $\cF$ be a cuspidal character sheaf with supporting stratum $\cS$. Then
$\cF|_{f^{-1}(f(\cS))}$ is a clean extension of $\cF|_\cS$.
\end{coro}

\begin{proof}
In the proof of part (ii) of \autoref{thm:main} we show that the restriction of $\cF$ to the fibre $F_x$ of $f$ at $f(x)\in f(\cS)$ is clean for every $x\in \cS$. Taking into account \autoref{prop:fibration}, we conclude that the corollary holds. 
\end{proof}

We will now prove part (iii) of \autoref{thm:main} making use of~\autoref{lemm:supp-stratum-cusp},~\autoref{prop:supercuspidal},  \autoref{cor:dim est-biorbitalcuspidal} and~\autoref{prop:annihlation-distinguished}. The proofs of the propositions are deferred to~\autoref{sec:parabolic-vanishing},~\autoref{sec:LY} and~\autoref{sec:support} respectively.

We may assume without loss of generality that $Z(G)_0 = 1$. Recall that $\cP_{\chi}^{\dag}$ is anti-orbital, and thus \autoref{prop:stratum-support} implies that it is constructible along the stratification defined in \autoref{stratification}. Let $\cF\subseteq \cP_{\chi}^{\dag}$ be a simple subsheaf and let $\cS'$ denote its supporting stratum.  Since $\cF$ is cuspidal by \autoref{thm:main}, it follows from \autoref{thm:cuspidal-distinguished} that $\dim f(\cS')=N$, where $N = \dim \fg_1 - \dim \fg_0$. Since $\cS$ contains $\rc$, \autoref{cor:dim est-biorbitalcuspidal} implies that  $\dim f(\cS)=N$, and hence $\dim f(\cS) =  \dim f(\cS').$ 

From the cleanness of $\chi$ and \autoref{thm-adjunction}, we see that 
\[
\dim\Hom(\cF|_{\rc},\chi[N+\dim \rc]) = \dim\Hom(\cF,\IC(\chi)[N]) = \dim\Hom(\cF,\cP_{\chi}^{\dag})\neq 0,
\]
so that $\cF|_{\rc}\neq 0$ and hence $\cS\subseteq \bar \cS'$. Thus $f(\cS) = f(\cS')$ and we further conclude from \autoref{lemm:supp-stratum-cusp} that $\cS = \cS'$. We have shown that $\cS$ is the supporting stratum of every simple subsheaf of $\cP_{\chi}^{\dag}$. By Verdier duality, $\cS$ is also the supporting stratum of every simple quotient sheaf of $\cP_{\chi}^{\dag}$. \par

We will now show that $\cS$ is open in the support of $\cP_{\chi}^{\dag}$. Let $\cS'$ be any generic stratum of $\supp \cP_{\chi}^\dag$ such that  $\cS \subseteq \overline{\cS'}$. We will show that $\cS = \cS'$.   \autoref{lem:parab-dag}, \autoref{prop:parabolic-vanishing} and \autoref{prop:annihlation-distinguished} imply that $\dim f(\cS) =N=  \dim f(\cS')$ so that $f(\cS) = f(\cS')$. Let $\cL\subseteq \cP_{\chi}^{\dag}|_{\cS'}[-\dim \cS']$ be a simple sub local system. Then, $\IC(\cL)$ is a composition factor of $\cP_{\chi}^{\dag}$, and in particular anti-orbital. So $\IC(\cL)$ is a character sheaf.   Let $\rc'\subseteq\cS'$ be a $G_0$-orbit and $\chi'\subseteq \cL|_{\rc'}$ a simple summand. Then $(\rc',\chi')$ is a supercuspidal datum by \autoref{prop:annihlation-distinguished} applied to $\cP_{\chi}^{\dag}$. Thus \autoref{prop:supercuspidal} implies that $\Hom(\IC(\chi'), \cF[k]) = 0$ for every non-cuspidal character sheaf $\cF$ and every $k\in \Z$, whereas the cleanness of $\cL|_{\rc'}$ implies
\[
\Hom(\IC(\chi'), \IC(\cL)[-\codim_{\cS'}\rc']) \cong \Hom(\chi', \cL|_{\rc'})\neq 0.
\]
 Therefore, $\IC(\cL)$ is a cuspidal character sheaf. We claim  that $\IC(\cL)$  is a subsheaf of $\cP_{\chi}^{\dag}$ which will imply that  $\cS' = \cS$. Let $j':\cS'\to \fg_1$ denote the inclusion. We obtain by adjunction a morphism $\alpha:\li^{\on{p}}j'_!(\cL[\dim \cS'])\to \cP_{\chi}^{\dag}$. By \autoref{lemm:supp-stratum-cusp}, the natural map $\beta:\li^{\on{p}}j'_!(\cL[\dim \cS'])\to \IC(\cL)$ satisfies $f(\supp(\ker(\beta)))\subseteq \partial(f(\cS))$. As all irreducible subsheaves of $\cP_{\chi}^{\dag}$ are supported on $\cS$ we have $\alpha(\ker(\beta))=0$ and hence $\IC(\cL)$  is a subsheaf of $\cP_{\chi}^{\dag}$.\par

We have shown that $\cP_{\chi}^{\dag}$ is a perverse sheaf supported on the closure of $\cS$, extending the local system $\cM_{\chi}:= \cP_{\chi}^{\dag}|_{\cS}[-\dim \cS]$ and without non-zero subsheaf nor quotient sheaf concentrated in $\partial \cS$. It follows that $\cP_{\chi}^{\dag} \cong \IC(\cM_{\chi})$.\qed

\subsubsection{Part (iv) of~\autoref{thm:main}}\label{ssec:orthogonality}

 In the proof part (iv) of~\autoref{thm:main} we make use of the following lemma whose proof is straightforward. 
 \begin{lemm}\label{lem:rep-normal}
Let $\Gamma$ be a group and $N\subseteq \Gamma$ a finite normal subgroup. Let $\operatorname{Rep} \Gamma$ denote the category of $\C\Gamma$-modules. Then, we have:
\begin{enumerate}
\item
The quotient $H = \Gamma/N$ acts via conjugacy on the set $N^{\wedge} \subset \on{Rep} N$ of isomorphism classes of irreducible representations of $N$. 
\item
There is a decomposition into abelian subcategories
\[
\operatorname{Rep} \Gamma = \bigoplus_{[\chi]\in N^{\wedge}/H }\operatorname{Rep}_{[\chi]} \Gamma,
\]
where $\operatorname{Rep}_{[\chi]} \Gamma$ is the subcategory of $\C\Gamma$-modules whose restrictions to $N$ have all simple constituents in the $H$-orbit $[\chi]\subseteq N^{\wedge}$.
\item
If $V\in \operatorname{Rep}_{[\chi]} \Gamma$ is finite dimensional, then we have
\[
\Res^{\Gamma}_{N} V \cong \bigoplus_{\tau\in [\chi]} \tau^{\oplus d},\quad \text{where $d = \frac{\dim V}{\# [\chi] \dim \chi}$.}
\]
\end{enumerate}
\end{lemm}

By \autoref{thm:main}~(iii), we write $\cP_{\chi}^{\dag} = \IC(\cM_{\chi})$ and $\cP_{\chi'}^{\dag} = \IC(\cM_{\chi'})$ for some $G_0$-equivariant local systems $\cM_{\chi}$ on $\cS$ and $\cM_{\chi'}$ on $\cS'$. Assume first that $(\rc, \chi)$ and $(\rc', \chi')$ are equivalent. Choose $x= s + n\in \rc$ and $x' = s' + n\in \rc'$ such that $Z_G(s) = Z_G(s')$. The inclusions $\rc,\rc'\subseteq  \cS$ induce injective homomorphisms:
\[
 \pi_0(Z_{G_0}(x))\cong \pi^{G_0}_1(\rc, x)\to \pi^{G_0}_1(\cS, x),\quad \pi_0(Z_{G_0}(x'))\cong \pi^{G_0}_1(\rc, x')\to \pi^{G_0}_1(\cS, x').
\]
Choosing a path in $\fd^{\circ} = \{y\in \fd\mid Z_G(y) = Z_G(s)\}$ 
connecting $s$ and $s'$, we get an isomorphism
\[
 \pi^{G_0}_1(\cS, x)\cong \pi^{G_0}_1(\cS, x'),
\]
which identifies the image of $\pi_0(Z_{G_0}(x))$ with that of $\pi_0(Z_{G_0}(x'))$. It follows that,  for any $G_0$-equivariant local system $\cL$ on $\cS$, $\cL|_{\rc}$ can be identified with $\cL|_{\rc'}$, and the multiplicity of $\chi$ in $\cL|_{\rc}$ coincides with that of $\chi'$ in $\cL|_{\rc'}$. Combining with the adjunction formula~\autoref{thm-adjunction}, we get
\[
\dim\Hom(\cL,\cM_{\chi}) = \dim\Hom(\cL|_{\rc},\chi) = \dim\Hom(\cL|_{\rc'},\chi') = \dim\Hom(\cL,\cM_{\chi'})
\]
for every $G_0$-equivariant local system $\cL$ on $\cS$ such that $\IC(\cL)$ is a character sheaf on $\fg_1$. In particular, the socles of $\cP_{\chi}^\dagger$ and  $\cP_{\chi'}^\dagger$ are isomorphic. Since $\cP_{\chi}$ and $\cP_{\chi'}$ are both projective in $\Perv_{G_0}(\fg_1^{\nil})$, we conclude that $\cP_{\chi} \cong \cP_{\chi'}$. \par
Assume now that $(\rc, \chi)$ and $(\rc', \chi')$ are not equivalent. \autoref{thm:vanishing} implies that $\Hom(\cP_{\chi}, \cF[k]) = 0$ for every $\cF\in \Perv_{G_0}(\fg_1^{\nil})$ and every $k\in \Z-\{0\}$. In particular, we have $\Hom(\cP_{\chi}, \cP_{\chi'}[k]) = 0$ for $k\neq 0$. We need to show that
\[
\Hom(\cP_{\chi}, \cP_{\chi'}) = \Hom(\cP_{\chi'}^{\dagger}, \cP_{\chi}^{\dagger}) = 0
\]
Let $\cS$ (resp. $\cS'$) denote the stratum containing $\rc$ (resp. $\rc'$). Since $\cP_{\chi}^{\dagger} = \IC(\cM_{\chi})$ and $\cP_{\chi'}^{\dagger} = \IC(\cM_{\chi'})$, we have $\Hom(\cP_{\chi'}^{\dagger}, \cP_{\chi}^{\dagger}) = 0$ if $\cS \neq \cS'$. We may thus assume $\cS = \cS'$ and it suffices to show that $\Hom(\cM_{\chi'}, \cM_{\chi}) = 0$. We may further assume that $\rc = \rc'$. Let $\cL$ be a simple quotient of $\cM_{\chi}$. Then, the adjunction formula~\autoref{thm-adjunction} shows that $\cL|_{\rc}$ contains $\chi$. Let $x\in\rc$. Applying~\autoref{lem:rep-normal} to $\Gamma = \pi_1^{G_0}(\cS, x)$, $N = \pi_0(Z_{G_0}(x))$ and $V = \cL_x$ we see that every simple factor of $\cL|_{\rc}$ is $W_{\fd}$-conjugate to $\chi$. Since this holds for every simple quotient of $\cM_{\chi}$, the same lemma also implies that every simple factor of $\cM_{\chi}|_{\rc}$ is $W_{\fd}$-conjugate to $\chi$. Similarly, every simple factor of $\cM_{\chi'}|_{\rc}$ is $W_{\fd}$-conjugate to $\chi'$. As $\chi$ and $\chi'$ are not $W_{\fd}$-conjugate by our assumption, we have 
\[
\Hom(\cM_{\chi'}, \cM_{\chi})\subseteq \Hom(\cM_{\chi'}|_{\rc}, \cM_{\chi}|_{\rc}) = 0.
\]
 \qed

\subsection{The case of stable gradings}\label{ssec:stable}
We call a vector $z\in \fg_1$ \emph{stable} (in the sense of geometric invariant theory) if it is semisimple and the stabiliser $Z_{G_0}(z)$ is finite. We call the grading $\fg_*$ \emph{stable} if $\fg_1$ admits a stable vector. The cuspidal character sheaves on a stably graded Lie algebra are studied in \cite{VX23st} making use of the nearby-cycle construction in \cite{GVX20}. The corollary below combined with \autoref{thm:main} shows that the construction therein produces indeed all cuspidal character sheaves in the stable case. In particular, the character sheaves given in~\cite[\S7.2]{VX23st} are precisely the cuspidal character sheaves for stably graded Lie algebras of classical types. This confirms Conjecture 7.8 of~\cite{VX23st}. 
\begin{coro}\label{coro:GIT-stable}
Suppose that the grading $\fg_*$ is stable. Then, a character sheaf $\cF\in \Char_{G_0}(\fg_1)$ is cuspidal if and only if $\cF$ has full support, i.e., $\supp\cF = \fg_1$.
\end{coro}
\begin{proof}
Assume that $\cF$ is cuspidal and let $\widecheck\cO$ denote the supporting stratum of $\cF$.  Then, \autoref{thm:cuspidal-distinguished}~(i) implies that 
\[
\dim f(\widecheck{\cO}) = \dim \fg_1 - \dim \fg_0 + \dim Z(\fg)_0. 
\]
The stability of $\fg_*$ implies that $Z(\fg)_0 = 0$ and that
\[
\dim \fg_1 /\!/ G_0 = \dim \fg_1 - \dim \fg_0.
\]
It follows that $f(\widecheck{\cO})$ is open and dense in $\fg_1/\!/ G_0$. Since the generic fibre of $f$ is a closed $G_0$-orbit in $\fg_1$ by stability, we deduce that $\supp \cF = \fg_1$. \par
Assume conversely that $\cF$ is not cuspidal. Then, there exists a $\theta$-stable proper parabolic subgroup $P\subseteq G$ such that $\supp \cF \subseteq G_0\fp_1$.  Note that the latter is a closed subset of $\fg_1$ because the map $G_0\times^{P_0}\fp_1\to \fg_1$ induced by the adjoint action is proper and $G_0\fp_1$ is its image.  Let $\lambda:\Cc\to G_0$ be a cocharacter such that $\fp = \li^\lambda_{\ge 0}{\fg}$. If $z\in G_0\fp_1$, then $z\in \Ad_g\fp_1$ for some $g\in G_0$ and hence $z\in \li^{\lambda'}_{\ge 0}{\fg}_1$, where $\lambda' = \Ad_g^{-1}\circ \lambda$. By the Hilbert--Mumford criterion, $z$ is not stable. Thus, $\supp \cF\subseteq G_0\fp_1$ consists of non-stable vectors and hence, $\supp\cF\neq \fg_1$. 
\end{proof}

\section{Annihilation under parabolic restriction (Proof of \autoref{prop:parabolic-vanishing})}\label{sec:parabolic-vanishing}

 In this section, we prove \autoref{prop:parabolic-vanishing}.  The proof has two ingredients: \autoref{prop:hyperbolic-parab} (the hyperbolic restriction theorem \cite[Theorem1]{B03}) and a vanishing property of supercuspidal data (\autoref{prop:supercuspidal}).

For a $\theta$-stable parabolic subgroup $P$ with unipotent radical $U$ and $\theta$-stable Levi factor $L$, we write $U_0=U^\theta$, $\fp=\on{Lie}P$, $\fu = \Lie U$, $\fl=\on{Lie}L$, and $\fp_1=\fp\cap\fg_1$ etc.

\subsection{Proof of~\autoref{prop:supercuspidal}} \label{ssec:proof-supercuspidal}
 
Let $P\subseteq G$ be a $\theta$-stable parabolic subgroup with unipotent radical $U$.  The cleanness of $\chi$ implies that 
\[
(\Res^{\fg_1}_{\fp_1}\IC(\chi))_{y} = \mathrm{R}\Gamma_c(\rc\cap (y + \fu_{  1}), \chi|_{\rc\cap (y + \fu_{  1})})[\dim \rc + \dim \fu_1 - \dim \fu_0]
\]
for $y\in \fl_1$. Therefore, the vanishing of $\Res^{\fg_1}_{\fp_1}\IC(\chi)$ is equivalent to that
\beq\label{eqn-vanishing1}
\mathrm{R}\Gamma_c(\rc\cap (x + \fu_{1}), \chi|_{\rc\cap (x + \fu_{1})})=0\,
\eeq
for every $x\in \rc\cap \fp_1$.

To prove~\autoref{prop:supercuspidal}, we begin with the following lemmas:

    \begin{lemm}\label{lem:hp}
  Let $s\in\fp_1$ be a semisimple element. Then, there exists a $\theta$-stable Levi factor $\fl\subseteq \fp$ such that $s\in \fl_1$. 
    \end{lemm}
    \begin{proof}
   
    Let $S$ be the minimal algebraic subgroup of $G$ whose Lie algebra contains $s$. Then, $S$ is a $\theta$-stable subtorus contained in $P$. Let $T\subseteq P$ be a maximal $\theta$-stable subtorus containing $S$. We claim that $T$ is a maximal torus of $P$ --- indeed, $Z_P(T) / T$ must be a unipotent group since otherwise the pre-image of any $\theta$-stable maximal torus of $Z_P(T) / T$ (which exists by Steinberg's theorem) under the quotient map $Z_P(T)\to Z_P(T) / T$ would be a bigger $\theta$-stable subtorus of $P$. It follows that there is a unique Levi factor $L$ of $P$ containing $T$ and its Lie algebra $\fl = \Lie L$ contains $s$. The uniqueness of $L$ implies that it is $\theta$-stable. 
   \end{proof} 
\begin{lemm}\label{lem:Ux}
 
Let $x=x_s+x_n\in \fp_{1}$.  There is an isomorphism
\begin{equation}\label{eq:iso-x+u}
 U_{0}\times^{Z_{U_{0}}(x_s)}(x_n + Z_{\fu_{ 1}}(x_s))\xrightarrow{\sim} x + \fu_{ 1}, \quad [u:x_n + z] \mapsto \Ad_u(x+z).
\end{equation}
\end{lemm}
\begin{proof}
Choose a $\theta$-stable Levi factor $\fl\subseteq \fp$ containing $x_s$, as in  \autoref{lem:hp}. Let $\lambda: \Cc\to G_{0}$ be a cocharacter such that $\fp = \li^\lambda_{\ge 0}{\fg}$, $\fl = \li^\lambda_0{\fg}$ and $\fu = \li^{\lambda}_{>0}{\fg}$. \par
We first show that for each $n \in \Z$:
\begin{equation}\label{eq:decomp-g-1}
\li^\lambda_{\ge n}{\fg}_{ 1} = [x, \li^\lambda_{\ge n}{\fg}_{0}] + \li^\lambda_{\ge n}{Z}_{\fg_{ 1}}(x_s).
\end{equation}
Since $x_n\in \fp_{  1}$, we may write $x_n = \li^{\lambda}_0{x}_n + \li^{\lambda}_{>0}{x}_n$, where $\li^{\lambda}_0{x}_n\in Z_{\fl_{  1}}(x_s)$ and $\li^{\lambda}_{>0}{x}_n\in Z_{\fu_{  1}}(x_s)$. Since $\li^{\lambda}_0{x}_n$ is nilpotent (it lies in the closure of $G_0x_n$ because $\lim_{t\to 0}\Ad_{\lambda(t)}x_n = \li^\lambda_{0}{x}_n$), we may find an $\fsl_2$-triple $(\li^{\lambda}_0{x}_n, h, f)$ with $h\in Z_{\fl_0}(x_s)$ and $f\in Z_{\fl_{-1}}(x_s)$ and a cocharacter $\mu:\Cc\to Z_{L_{0}}(x_s)$ such that $d\mu(1) = h$. In particular, $\li^\lambda_{0}{x}_n\in \li^{\mu}_2{\fl}_{ 1}$ and $x_s\in \li^{\mu}_0{\fl}_{1}$ hold. 
It follows that for each $j,k\in \Z$, we have
\begin{equation}\label{eq:decomp-g-2}
    [x_s, \li^\mu_{k}{}\li^\lambda_{j}{\fg}_{0}]\subseteq \li^\mu_{k}{}\li^\lambda_{j}{\fg}_{1},\quad  [\li^{\lambda}_0{x}_n, \li^\mu_{k}{}\li^\lambda_j{\fg}_{0}]\subseteq \li^\mu_{(k+2)}{}\li^\lambda_{j}{\fg}_{1}.
\end{equation}
Since $x_s$ is semisimple, there is a decomposition for each $j,k\in \Z$:
\begin{equation}\label{eq:decomp-g-3}
\li^\mu_{k}{}\li^\lambda_{j}{\fg}_{1} = [x_s, \li^\mu_{k}{}\li^\lambda_{j}{\fg}_{0}] \oplus \li^\mu_{k}{}\li^\lambda_{j}{Z_{\fg_{  1}}(x_s)}.
\end{equation}
For any fixed $j\in \Z$, we can easily deduce from \eqref{eq:decomp-g-2} and \eqref{eq:decomp-g-3} by descending induction on $k\in \Z$ that:
\[
\li^\mu_{\ge k}{}\li^\lambda_{j}{\fg}_{  1} = [x_s + \li^\lambda_0{x}_n, \li^\mu_{\ge k}{}\li^\lambda_{j}{\fg}_{0}] + \li^\mu_{\ge k}{}\li^\lambda_{j}{Z_{\fg_{  1}}(x_s)}.
\]
When $k \ll 0$, the above identity becomes  
\begin{equation}\label{eq:decomp-g-4}
\li^\lambda_{j}{\fg}_{  1} = [x_s + \li^\lambda_0{x}_n, \li^\lambda_{j}{\fg}_{0}] + \li^\lambda_{j}{Z_{\fg_{  1}}(x_s)}.
\end{equation}
On the other hand, there is an inclusion
\[
[\li^{\lambda}_{>0}{x}_n, \li^\lambda_{\ge j}{\fg}_{0}]\subseteq \li^\lambda_{\ge j + 1}{\fg}_1.
\]
for each $j\in \Z$. Combining \eqref{eq:decomp-g-4} and the above inclusion, we deduce \eqref{eq:decomp-g-1} by descending induction on $j\in \Z$. \par

For $j > 0$ and $z\in \li^\lambda_{\ge j}{\fg}_{1} + Z_{\fu_{  1}}(x_s)$, we can find $g\in U_{0}$ such that $$\Ad_g (x+z)\in x + \li^\lambda_{\ge j + 1}{\fg}_{1}+  Z_{\fu_{  1}}(x_s).$$ Indeed, \eqref{eq:decomp-g-1} implies that there exists $y\in \li^\lambda_{\ge j}{\fg}_{0}$ and $z'\in Z_{\fu_{  1}}(x_s)$ such that $z = [x,y] + z'$. Then $g = \exp(y)$ is a desired element. 
By induction, it follows that given $z\in \fu_{  1}$, we can find $g_i\in U_{0}$ for every $i\ge 1$ such that $\Ad_{g_j}\cdots\Ad_{g_1}(x + z) \in x + \li^\lambda_{\ge j + 1}{\fg}_{  1} + Z_{\fu_{  1}}(x_s)$ for every $j \ge 1$. When $j$ is large, the element $g = g_j\cdots g_1\in U_{0}$ satisfies $\Ad_{g}(x + z) \in x + Z_{\fu_{  1}}(x_s)$; therefore, the map \eqref{eq:iso-x+u} is surjective. \par
Let $z\in Z_{\fu_{  1}}(x_s)$. The element $\li^\lambda_{0}{x}_n$ is nilpotent and lies in the closure of $G_0(x_n + z)$ because $\lim_{t\to 0}\Ad_{\lambda(t)}(x_n + z) = \li^\lambda_{0}{x}_n$; hence, $x_n + z$ is nilpotent and commute with $x_s$; therefore, $x_s$ is the semisimple part of $x + z$. It follows that the pre-image of $x + z$ under the map \eqref{eq:iso-x+u} is the singleton $\{[1:x_n + z]\}$. Similar arguments show that the map is \'etale at $[1:x_n + z]$. This proves that \eqref{eq:iso-x+u} is an isomorphism. 
\end{proof}

We now proceed to prove~\autoref{prop:supercuspidal}. Let $P\subseteq G$ be a $\theta$-stable proper parabolic subgroup of $G$ with unipotent radical $U$ and let $x\in\rc\cap\fp_1$. We shall prove~\eqref{eqn-vanishing1}. Recall $Z_x = Z(Z_{G}(x_s))$ and $H_x  = [Z_{G}(x_s),Z_{G}(x_s)]$ defined in \eqref{eqn-hz}. By \autoref{lem:hp}, we can choose a $\theta$-stable Levi factor $L\subseteq P$ such that $x_s\in \fl_1 = (\Lie L)_1$. We choose a cocharacter $\lambda:\Cc\to G_0$ such that $\fp = \li^{\lambda}_{\ge 0}{\fg}$ and $\fl =\li^{\lambda}_{0}{\fg}$. In particular, we have $\lambda(\Cc)\subseteq Z_{G_0}(x_s)$. 
 Note that $\fZ_{x,0} = Z(\fg)_{0}$ since $x$ is $G_0$-distinguished by assumption. The hypothesis $\fp \neq \fg$ then implies that $\lambda(\Cc)\not\subseteq Z_{x,0}$. This in turn implies that $\fh_x\cap \fl = Z_{\fh_x}(\lambda)\neq \fh_x$. Hence, $\fh_x\cap \fp\neq \fh_x$ is a proper parabolic subalgebra of $\fh_x$. 
 
 Consider the restriction $\chi' = \chi|_{x_s + H_{x,0} x_n}$. When translated by $-x_s$, we obtain a clean local system on the orbit $\rc' := H_{x, 0}x_n$. Since $\rc'$ is a $H_{x,0}$-distinguished nilpotent orbit in $\fh_{x, 1}$, \autoref{characterisation of supercuspidal} implies that 
     \[
         \mathrm{R}\Gamma_c(\rc'\cap (\fh_{x,1}\cap \fu_1), \chi'|_{\rc'\cap (\fh_{x,1}\cap \fu_1)}) = 0.
     \]
     Intersecting both sides of \eqref{eq:iso-x+u} with $\rc=G_{0}x$,  we obtain:
     \[
     U_{0}\times^{H_{x, 0}\cap U_0}(\rc'\cap (\fh_{x,1}\cap \fu_1))\xrightarrow{\sim} \rc\cap (x + \fu_{  1}).
     \]
     In particular, $\rc\cap (x + \fu_{  1})$ is an affine bundle over $\rc'\cap (\fh_{x,1}\cap \fu_1)$ with fibre $U_0 / (H_{x, 0}\cap U_0)$. Hence, we have
     \[
         \bega
             \mathrm{R}\Gamma_c(\rc\cap (x + \fu_{  1}), \chi|_{\rc\cap (x + \fu_{  1})}) \\
             \cong \mathrm{R}\Gamma_c(\rc'\cap (\fh_{x,1}\cap \fu_1), \chi'|_{\rc'\cap (\fh_{x,1}\cap \fu_1)})[2\dim U_0 / (H_{x, 0}\cap U_0)] = 0.
         \eega
     \]
This completes the proof of \autoref{prop:supercuspidal}. \qed

\subsection{Proof of \autoref{prop:parabolic-vanishing}}

Let $P\subseteq G$ be a $\theta$-stable proper parabolic group with unipotent radical $U$ and $L\subseteq P$ a $\theta$-stable Levi factor. Let $P'\subseteq G$ be the parabolic subgroup opposite to $P$ containing $L$.

Let $x\in\rc$. Consider the following commutative diagrams:
    \[
     \begin{tikzcd}
     \fl_{1} & \fl_1^{\nil}\arrow{l}[swap]{\iota_{L}}\\
       \fp_{1} \arrow{d}{i}\arrow{u}[swap]{p} & \fp^{\nil}_{1}\arrow{d}{i}\arrow{u}[swap]{p} \arrow{l}[swap]{\iota_{P}}\\
       \fg_{1} & \fg^{\nil}_{1} \arrow{l}[swap]{\iota}\\
     \end{tikzcd}
     \quad 
     \begin{tikzcd}
     \cU_{L} \arrow{r}{j_{L}} & \fl_{ 1} & \fl^{\nil}_{ 1} \arrow{l}[swap]{\iota_{L}}\\
      \cU_{P'} \arrow{u}[swap]{p'_{\cU}}\arrow{d}{i'_{\cU}}\arrow{r}{j_{P'}} & \fp'_{ 1}\arrow{d}{i'}\arrow{u}[swap]{p'} & \fp'^{\nil}_{ 1} \arrow{l}[swap]{\iota_{P'}}\arrow{d}{i'}\arrow{u}[swap]{p'} \\
      \cU \arrow{r}{j} \arrow{d}{\pi_{\cU}} &  \fg_{ 1}\arrow{d}{f} &  \fg^{\nil}_{ 1}\arrow{l}[swap]{\iota}\arrow{d} \\
      \Cc  \arrow{r}{r\mapsto rf(x)} &  \fg_{ 1}/\!/G_{ 0} &  \{0\}\arrow{l}
     \end{tikzcd}
     \]
     where $\cU = \{(z, r)\in \fg_{ 1}\times \Cc\;|\; f(z) = rf(x)\}$, $\cU_{P'} = \cU\times_{\fg_{ 1}} \fp'_{ 1}$, $\cU_L = \cU\times_{\fg_{ 1}} \fl_{ 1}$ and $j$, $\pi_{\cU}$ are the canonical projections. 
     
     Let $\cI$ be the (infinite-rank) local system on $\Cc$ associated with the regular representation of the fundamental group $\pi_1(\Cc, 1)$.  
     By~\cite[(8.6.2)]{KS90} for example, the nearby cycles can be written as $\psi_\pi \cL_\chi= \iota^*j_*\tilde \cL_\chi[-1]$, where $\tilde \cL_\chi = \scrHom(\pi_{\cU}^* \cI, \cL_\chi)$. 
     We have 
     \begin{equation}\label{eq:res-psi}
     \begin{gathered}
     \Res^{\fg_{ 1}}_{\fp_{ 1}} \cP_{\chi}[-\dim \fu_1 + \dim \fu_0 + 1] = p_!i^* \iota_!\iota^*j_*\tilde \cL_\chi \cong \iota_{L !}\iota_L^*p_!i^* j_*\tilde \cL_\chi\\ 
     \cong \iota_{L !}\iota_L^*p'_*i'^! j_*\tilde \cL_\chi \cong \iota_{L !}\iota_L^*j_{L*}p'_{\cU*}i_{\cU}'^! \tilde \cL_\chi,
     \end{gathered}
     \end{equation}
     where the isomorphism $p_!i^*\cong p'_*i'^!$ comes from  \autoref{prop:hyperbolic-parab}.  Given $\zeta = (z, r)\in \cU_{L}$, the $!$-fibre of $p'_{\cU*}i_{\cU}'^! \tilde\cL_\chi$ at $\zeta$ is given by
     \begin{equation} \label{eq:fibre-res-psi}
     \RHom(\cI_{r}, \mathrm{R}\Gamma_c((z + \fu_{ 1})\cap \rc, \chi\mid_{(z + \fu_{ 1})\cap \rc}))),
     \end{equation}
     where $\cI_r$ denotes the $*$-fibre of $\cI$ at $r$. By \eqref{eqn-vanishing1}, the term in \eqref{eq:fibre-res-psi} vanishes.  
 It follows that $p'_{\cU*}i_{\cU}'^! \tilde \cL_\chi = 0$, which concludes the proof. \qed
 \begin{rema}
 The isomorphisms in \eqref{eq:res-psi} show that the nearby-cycle functor commutes with hyperbolic restriction. This is a formal consequence of the hyperbolic restriction theorem and has been shown by Nakajima in \cite[5.4.1(2)]{N17}. Note that, thanks to \autoref{rema:non-construcitble}, the hyperbolic restriction theorem is applicable to the complex $j_* \tilde\cL_\chi$ even though the latter is not constructible. 
 \end{rema}

\section{Spiral induction and bi-orbital cuspidal character sheaves}\label{sec:LY}
In this section, we recall briefly the notion of spiral induction and restriction introduced by Lusztig--Yun in \cite{LY17a} as well as some of their main results. We explain its relationship with hyperbolic restriction and how the spiral induction interacts with the functor $\dag$ from \autoref{ssec:functor dag}. We deduce a necessary condition for the existence of bi-orbital cuspidal character sheaves.

\subsection{Spirals and splittings}
For $\epsilon\in \{\pm 1\}$, $\lambda\in \bX_*(G_{0})$, $r\in\Z_{>0}$ and $n\in\Z$, define
\[
   \li^\epsilon{\ul\fp}_n = \li^{\lambda}_{\ge \epsilon nr}{\fg}_{ n},\quad\li^\epsilon{\ul\fu}_n = \li^{\lambda}_{> \epsilon nr}{\fg}_{ n},\quad \li^\epsilon{\ul\fl}_n = \li^{\lambda}_{\epsilon nr}{\fg}_{ n},
\]
and let
$$\li^\epsilon{\ul\fl}_* = \bigoplus_{n\in \Z} \li^\epsilon{\ul\fl}_n.$$
Let $^\epsilon\ul\fl$ denote the underlying (ungraded) Lie algebra of $^\epsilon\ul\fl_*$; it is a Lie subalgebra of $\fg$.
Let $\ul P_0 = \exp(^\epsilon\ul\fp_0)$ and $\ul L_0 = \exp(^\epsilon\ul\fl_0)$;  a parabolic subgroup and a Levi subgroup of $G_{0}$, respectively. The collection of subspaces $\li^\epsilon{\ul\fp}_*=\{\li^\epsilon{\ul\fp}_n\}_{n\in\Z}$ arising in this way is called an \emph{($\epsilon$-)spiral} of $\fg_*$. We call $\li^\epsilon{\ul\fu}_*$ the \emph{radical} of $\li^\epsilon{\ul\fp}_*$ and $\li^\epsilon{\ul\fl}_*$ a \emph{splitting} of $\li^\epsilon{\ul\fp}_*$. \par

\subsection{Spiral induction, restriction and co-restriction}\label{ssec:spiral-ind}
Let $\epsilon\in \{\pm 1\}$. Given an $\epsilon$-spiral $(\ul\fp_*, \Pz)$ on $\fg_*$ and a splitting $(\ul\fl_*, \Lz)$ of $(\ul\fp_*, \Pz)$, consider the following diagram of quotient stacks
\[
[\flo / \Lz]\xleftarrow{a} [\fpo / \Pz]\xrightarrow{b} [\fg_{1} / G_{ 0}],
\]
where $a$ is induced by the projections $\fpo\to \flo$ and $\Pz\to \Lz$, whereas $b$ is induced by the natural inclusions. Let $\ul\fu_*$ be the radical of $\ul\fp_*$.
Spiral induction, spiral restriction and spiral co-restriction are defined to be the following functors:
 \beqn
    \begin{gathered}
    \Ind^{\fg_ {1}}_{\ul\fp_ {1}}: \Db_{\ul L_{0}}(\ul\fl_ {1})\rightleftarrows \Db_{G_{0}}(\fg_ {1}):\Res^{\fg_ {1}}_{\ul\fp_ {1}}, \res^{\fg_ {1}}_{\ul\fp_ {1}}\\ 
    \Ind^{\fg_ {1}}_{\ul\fp_ {1}} = b_!a^*[\dim \ul\fu_1 + \dim \ul\fu_0] \cong b_*a^![-\dim \ul\fu_1 - \dim \ul\fu_0], \\
    \Res^{\fg_ {1}}_{\ul\fp_ {1}} = a_!b^*[\dim \ul\fu_1 + \dim \ul\fu_0],\; \res^{\fg_ {1}}_{\ul\fp_{1}} = a_*b^![-\dim \ul\fu_1 - \dim \ul\fu_0].
    \end{gathered}
 \eeqn
They satisfy adjunctions: 
\beq\label{eqn-adjunctions}
\NRes^{\fg_{1}}_{\ul\fp_1}\dashv\NInd^{\fg_{1}}_{\ul\fp_1}\dashv\Nres^{\fg_{1}}_{\ul\fp_1}.
\eeq
The cohomological shifts are chosen so that the spiral induction commutes with the Verdier duality, see \cite[4.1(d)]{LY17a}. \par
When $\epsilon=1$, we have $\ul\fp_1\subseteq \fg^{\nil}_1$, so that the image of $\Ind^{\fg_{1}}_{\fpo}$ lies in the full subcategory of orbital complexes $\Db_{G_{ 0}}(\fg^{\nil}_{1})\subseteq \Db_{G_{ 0}}(\fg_{1})$; when $\epsilon=-1$, the image of $\Ind^{\fg_{1}}_{\fpo}$ lies in the full subcategory of anti-orbital complexes $\Dbn_{G_{ 0}}(\fg_{1})\subseteq \Db_{G_{ 0}}(\fg_{1})$, see \cite[\S 7.2, \S 7.7]{LY17a}.

\subsection{LY cuspidal supports and block decomposition}\label{ssec:LYSupp}
 A \emph{Lusztig--Yun (LY) cuspidal support} on $\fg_{1}$ is the datum $\xi = (\underline{M}, \ul\fm_*, \cC)$ consisting of a connected reductive subgroup $M\subseteq G$, a $\Z$-grading $\ul\fm_*$ on its Lie algebra $\ul\fm = \Lie \ul M$ arising as a splitting of some spiral of $\fg_*$ and a cuspidal local system $\cC$ on the dense open $\ul M_0$-orbit of $\ul\fm_1$ in the sense of \cite{L95}.

 We have the following theorem due to Lusztig--Yun:
\begin{theo}[{\cite[Theorem 0.6]{LY17a}}]\label{thm:orthogonal-decomposition}
There are orthogonal decompositions
\[
\Db_{G_{ 0}}(\fg^{\nil}_{1}) = \bigoplus_{\xi}\Db_{G_{ 0}}(\fg^{\nil}_{1})_{\xi}\quad\text{and} \quad \Dbn_{G_{ 0}}(\fg_{1}) = \bigoplus_{\xi}\Dbn_{G_{ 0}}(\fg_{1})_{\xi},
\]
where $\xi$ runs over all $G_{ 0}$-conjugacy classes of LY cuspidal supports on $\fg_{1}$. 
\end{theo}
The block $\Db_{G_{ 0}}(\fg^{\nil}_{1})_{\xi}\subseteq \Db_{G_{ 0}}(\fg^{\nil}_{1})$ (resp. $\Dbn_{G_{ 0}}(\fg_{1})_{\xi}\subseteq \Dbn_{G_{ 0}}(\fg_{1})$) is defined to be the thick triangulated subcategory generated by the family $\{\Ind^{\fg_{1}}_{\fpo}\IC(\cC)\}_{\ul\fp_*}$, where $\ul\fp_*$ runs over all $1$-spirals (resp. $(-1)$-spirals) of $\fg_*$, such that $\ul\fm_*$ is a splitting of $\ul\fp_*$.

\subsection{Hyperbolic restriction and spiral restriction}\label{ssec:hyperbolic-spiral}

Let $(\ul\fp_*, \ul\fl_*, \ul\fp'_*)$ be a triplet formed by a $1$-spiral $\ul\fp_*$, a $(-1)$-spiral $\ul\fp'_*$ and a common splitting $\ul\fl_*$ of $\ul\fp_*$ and $\ul\fp'_*$ such that $\ul\fp_n\cap \ul\fp'_n = \ul\fl_n$ for all $n\in\Z$. We call $(\ul\fp_*, \ul\fl_*, \ul\fp'_*)$ an \emph{opposite spiral triplet}.  \par

Recall the subcategory $\Db_{G_{ 0}}(\fg_{ 1})^{\Cunip}\subseteq \Db_{G_{ 0}}(\fg_{ 1})$ of unipotently $\Cc$-monodromic complexes introduced in \autoref{ssec:FS}. 
\begin{prop}\label{prop:hyperbolic-spiral}
There are isomorphisms of functors from $\Db_{G_{ 0}}(\fg_{ 1})^{\Cunip}$ to $\Db_{\Lz}(\flo)$:
\[
\Res^{\fg_{1}}_{\fpo} \cong \res^{\fg_{1}}_{\fpo'}[D],\quad \Res^{\fg_{1}}_{\fpo'}\cong \res^{\fg_{1}}_{\fpo}[ D]; \quad D = \dim \fg_1 - \dim \fg_0 - \dim \ul\fl_1 + \dim \ul\fl_0.
\]
\end{prop}
\begin{proof}
We choose a cocharacter $\lambda\in\bX_*(G_{ 0})$ and $r\in\Z_{>0}$ such that 
\[
\ul\fp_n = \li^\lambda_{\ge rn}{\fg}_n,\quad \ul\fl_n = \li^\lambda_{rn}{\fg}_n,\quad \ul\fp'_n= \li^\lambda_{\le rn}{\fg}_n
\] 
for every $n\in \Z$. Consider the diagram:
\[
\begin{tikzcd}
\fg_{1} & \fpo \arrow{l}{i}\arrow{d}{p}\\
\fpo' \arrow{u}{i'}\arrow{r}{p'} & \flo.
\end{tikzcd}
\]
It can be viewed as a hyperbolic restriction diagram with respect to the $\Cc$-action on $\fg_1$ defined by the following cocharacter of $G_0\times \Cc$:
\[
(2\lambda, r)\in \bX_*(G_0\times \Cc),\quad t\mapsto (\lambda(t^2), t^r).
\]
The hyperbolic restriction theorem (\cite[Theorem 1]{B03}) then yields the desired isomorphisms, as in the proof of \autoref{prop:hyperbolic-parab}. 
\end{proof}

We note that  \autoref{prop:hyperbolic-spiral} is applicable to orbital complexes, by \autoref{lem:orbital-monodromic}.  
 \par

\subsection{The functor \texorpdfstring{$\dag$}{\dag} and spiral induction} 
Recall the functor $\dag = \FS\circ\tau_*\circ \D$ from \autoref{ssec:functor dag}, where the isomorphism $\fg_1\cong \fg_{-1}$ is chosen with respect to the Chevalley involution $\tau_0:G_0\to G_0$ which inverts a fixed maximal torus $T_0\subseteq G_0$. Let $(\ul L, \ul\fl_*, \cC)$ be a LY cuspidal support. Let $\ul\fp_*$ be a $1$-spiral and $\ul\fp'_*$ a $(-1)$-spiral such that $(\ul\fp_*, \ul\fl_*, \ul\fp'_*)$ forms an opposite spiral triplet in the sense of \autoref{ssec:hyperbolic-spiral}.

\begin{prop}
\label{prop:involution-series}   There exists an isomorphism
		 \[
			 \left(\Ind^{\fg_1}_{\fpo} \IC(\cC)\right)^{\dag} \cong \Ind^{\fg_1}_{\fpo'} \IC(\cC).
		 \]
\end{prop}

  Conjugating $(\ul L, \ul\fl_*, \cC)$ by a suitable element of $G_0$, we may assume that $\Lz$ contains the maximal torus $T_0$.  
  Let $\lambda\in \bX_*(T_0)$ and $r > 0$ be such that $\ul\fp_n = \li^{\lambda}_{\ge rn}{ \fg}_n$,  $\ul\fp'_n = \li^{\lambda}_{\le rn} {\fg}_n$ and $\ul\fl_n = \li^{\lambda}_{rn} {\fg}_n$ hold for $n\in \Z$.
  The involution $\tau$ (resp. $\tau_0$) restricts to $\tau:\ul\fl_1\xrightarrow{\sim}\ul\fl_{-1}$  (resp. $\tau_0:\Lz\xrightarrow{\sim} \Lz$). By abuse of notation, we will denote 
   \[
       \dag = \FS\circ\tau_*\circ\D: \Db_{\ul L_0}(\ul\fl_1)^{\op}\to \Db_{\ul L_0}(\ul\fl_1).
   \]
	 \begin{lemm}
	 \label{lem:involution-induction} 
     There is an isomorphism of functors
     \[
     \Ind^{\fg_1}_{\fpo'}\circ\dag\cong \dag\circ\Ind^{\fg_1}_{\fpo}:\Db_{\Lz}(\flo)^{\op}\to \Db_{G_0}(\fg_1).
     \]
     \end{lemm}
	 \begin{proof} 
   
  Recall first that spiral induction commutes with Fourier--Sato transform by \cite[7.6(a)]{LY17a}.  The proof now proceeds completely analogously to that of~\autoref{lem:parab-dag}. 
\end{proof}
	\begin{lemm}\label{lem:involution-cuspidal}
	 
 There exists an isomorphism $\IC(\cC)^{\dag}\cong \IC(\cC)$. 
	\end{lemm}

	\begin{proof}Since $\fl_*$ is $\Z$-graded (graded by $\lambda/r$), we have $Z(\ul L) \subseteq \ul L_0$; thus $Z(\ul L)$ is contained in the fixed maximal torus $T_0\subseteq \Lz$. In particular, $\tau_0$ induces the inverse map on $Z(\ul L)$, so the central characters of $\IC(\cC)$ and $\tau_*\IC(\cC)$ are inverse to each other. On the other hand, Verdier duality inverts and Fourier--Sato transform preserves the central character. Therefore, $\IC(\cC)$ and $\IC(\cC)^{\dag} = \on{FS}\tau_*\D\IC(\cC)$ are cuspidal perverse sheaves on $\flo$ with the same central character. By a result of Lusztig in \cite{L84}, the isomorphism classes of cuspidal local systems on a given Lie algebra are determined by their central characters. Therefore, $\IC(\cC)$ and $\IC(\cC)^{\dag}$ are isomorphic. 
	\end{proof}
	
 \autoref{prop:involution-series} now follows immediately from Lemmas~\ref{lem:involution-induction} and~\ref{lem:involution-cuspidal}.

\subsection{Bi-orbital cuspidal sheaves and proof of ~\autoref{cor:dim est-biorbitalcuspidal}}\label{ssec:bi-orbital}

Let $\cF$ be a supercuspidal orbital sheaf on $\fg_1$. By \cite[Theorem 5.1.2]{Liu24}, there exists a LY cuspidal support  $\xi = (\ul L, \ul\fl_*, \cC)$ on $\fg_*$ and a $1$-spiral $\ul\fp_*$ of $\fg_*$ which contains $\ul\fl_*$ as a splitting, such that $Z(\ul L)_0 = Z(G)_0$ and $\cF \cong \Ind^{\fg_1}_{\ul\fp_1}\IC(\cC)$. Let $\ul\fp'_*$ be the $(-1)$-spiral such that $(\ul\fp_*, \ul\fl_*, \ul\fp'_*)$ forms an opposite spiral triplet and put $\cF' = \Ind^{\fg_1}_{\ul\fp'_1} \IC(\cC)$.  \autoref{prop:involution-series} gives us an isomorphism $\cF'\cong \cF^{\dag}\in \Perv^{\nil}_{G_0}(\fg_1)$.

    \begin{lemm}\label{lem:gdim} The graded dimension of $\RHom(\cF, \cF')$ is given by
    \[
    \sum_{k\in \Z} v^k\dim\Hom(\cF, \cF'[k]) = v^{\dim \fg_1 - \dim \fg_0 + \dim Z(\fg)_0}(1 - v^2)^{-\dim Z(\fg)_0}.
    \]
 
    \end{lemm}
    \begin{proof}
    Let $\ul\fu_*$ and $\ul\fu'_*$ denote the radicals of $\ul\fp_*$ and $\ul\fp'_*$. Applying \cite[Proposition 6.4]{LY17a}, we see that the left-hand side equals 
    \begin{equation}\label{eq:dim}
        (1 - v^2)^{-\dim Z(\fl)}\sum_{g}v^{\tau(g)},\quad\tau(g) =\dim\frac{\ul\fu'_1 + \Ad_g\ul\fu_1}{\ul\fu'_1\cap \Ad_g\ul\fu_1} -\dim\frac{\ul\fu'_0 + \Ad_g\ul\fu_0}{\ul\fu'_0\cap \Ad_g\ul\fu_0},
    \end{equation}
    where $g$ runs over a certain subset of representatives for $(\ul P'_0,\ul P_0)$-double cosets in $G_0$ such that $\Ad_g\ul\fp_*$ and $\ul\fp'_*$ admit a common splitting. We claim that any such $g$ must lie in $\ul P'_0\ul P_0$. Let $S\subseteq \ul L_0$ be a maximal torus. We may find $\lambda: \Cc\to S$ and $r\in \Z_{>0}$ such that 
    \[
    \ul\fp_n = \li^{\lambda}_{\ge rn}{\fg}_n,\quad \ul\fp'_n = \li^{\lambda}_{\le rn}{\fg}_n,\quad \ul\fl_n = \li^{\lambda}_{\ge rn}{\fl}_n.
    \]
    Since the set of splittings of $\ul \fp'_*$ forms a principal homogeneous $\ul U'_0$-space by \cite[\S2.7(a)]{LY17a}, there exists $u'\in\ul U'_0$ such that $\ul\fl_* \subseteq \Ad_{u'g} \ul\fp_*$. Similarly, there exists $u\in\ul U_0$ such that $\Ad_{u'gu} \ul\fl_* = \ul\fl_*$. Since the maximal tori of $\ul L_0$ are $\ul L_0$-conjugate to each other, there exists $h\in \ul L_0$ such that $\Ad_{g'} S = S$ with $g' = hu'gu$. It follows that the cocharacter $\mu:\Cc\to S$ defined by $\mu(t) = \lambda(t)^{-1}\Ad_{g'}\lambda(t)$ centralises $\ul\fl_*$. The condition $Z(\ul L)_0 = Z(G)_0$ implies that the image of $\mu$ lies in $Z(G)_0$ and thus $\Ad_{g'} \ul\fp_* = \ul\fp_*$, which implies $g'\in \ul P_0$ (since $\ul P_0$ is a parabolic subgroup of $G_0$) and proves the claim. \par
    
    By the discussion above, the summation in \eqref{eq:dim} is taken over the singleton $\{1\}$ and therefore equals $v^d(1 - v^2)^{-\dim Z(\fg)_0}$, where 
    \[
    d = \dim (\fuo + \underline{\fu}'_1) - \dim (\underline{\fu}_0 + \underline{\fu}'_0) = \dim \fg_1 - \dim \fg_0 - \dim \underline{\fl}_1 + \dim \underline{\fl}_0. 
    \]
    On the other hand, since $\ul\fl_1$ admits a cuspidal local system, the open $\ul L_0$-orbit in $\ul\fl_1$ is distinguished (see~\cite[\S 4.4]{L95}) and therefore $\dim\underline{\fl}_0 = \dim \underline{\fl}_1 + \dim Z(\ul\fl)_0$ see \cite[5.7.5]{C93}).  Hence, we have $d = \dim \fg_1 - \dim \fg_0 + \dim Z(\fg)_0$. 
    \end{proof} 
    \begin{prop}\label{prop:supercuspidal-dim}
        Suppose that there exists a supercuspidal (orbital or anti-orbital) sheaf on $\fg_1$. Then we have $\dim \fg_1 \ge \dim\fg_0- \dim Z(\fg)_0$. Moreover, the supercuspidal sheaf is bi-orbital if and only if the equality holds.
    \end{prop}
    \begin{proof}
        As was explained above, any supercuspidal orbital (resp. anti-orbital) sheaf is of the form $\cF = \Ind^{\fg_1}_{\ul\fp_1}\IC(\cC)$ (resp. $\cF' =  \Ind^{\fg_1}_{\ul\fp'_1}\IC(\cC)$) for some LY cuspidal support $\xi = (\ul L, \ul\fl_*, \cC)$. 
            
        By \autoref{lem:gdim}, we have 
        \[
        \Hom(\cF, \cF'[N]) \neq 0\quad \text{and}\ \  \Hom(\cF, \cF'[k]) = 0 \quad \text{for $k < N$},
        \] 
        where $N = \dim \fg_1 - \dim \fg_0 + \dim Z(\fg)_0$.
        Since both $\cF$ and $\cF'$ are perverse sheaves, we have $\Hom(\cF, \cF'[k]) = 0$ for $k < 0$ which forces  $N\ge 0$. If $\cF$ is not bi-orbital, then $\cF\not\cong\cF'$, so $\Hom(\cF, \cF') = 0$ and thus the strict inequality $\dim \fg_1 \ge \dim\fg_0- \dim Z(\fg)_0$ holds. If $\cF$ is bi-orbital, then $\Hom(\cF, \cF'[N]) \neq 0$ implies that $\cF$ and $\cF'$ lie in the same LY block, namely $\Db_{G_0}(\fg^{\nil}_1)_{\xi}$. Since $\cF$ is the only simple object of $\Perv_{G_0}(\fg^{\nil}_1)_{\xi}$ up to isomorphism (see \cite[Proposition 4.3.3]{Liu24}, for example) it follows that $\cF'\cong \cF$ and in particular $\dim\Hom(\cF, \cF') = 1$, which is equivalent to the condition $N = 0$.
    \end{proof}
 
    The following corollary shows that the existence of bi-orbital cuspidal character sheaves imposes a strong condition on the grading on $\fg$:
    \begin{coro}\label{coro:dim-biorbital}
    If there exists a bi-orbital cuspidal character sheaf on $\fg_1$, then $\dim\fg_0 - \dim Z(\fg)_0=\dim\fg_1$. Moreover, $\cF\cong \cF^{\dag}$ for every bi-orbital cuspidal character sheaf $\cF$. 
    \end{coro}
    \begin{proof}
     By \autoref{Prop3}, every bi-orbital cuspidal character sheaf $\cF\in \Char_{G_0}(\fg_1)$ is supercuspidal. The first statement now follows from \autoref{prop:supercuspidal-dim}. To prove the second statement, recall the isomorphism $\cF'\cong \cF^{\dag}$. Using  \autoref{lem:gdim} and the first part of the corollary, one concludes that $\dim \Hom(\cF, \cF')\neq 0$, proving the second statement.
    \end{proof} 
 
We deduce~\autoref{cor:dim est-biorbitalcuspidal} as a corollary.
    \begin{proof}[Proof of~\autoref{cor:dim est-biorbitalcuspidal}]
     Let $x\in\rc$ and recall the notation in~\eqref{eqn-hz}.   The equality
\beq\label{dimension est}
\dim \fg_{1} - \dim \fg_{0} = \dim \fh_{x,1} - \dim \fh_{x,0} + \dim \fZ_{x,1}-\dim Z(\fg)_0.
\eeq
follows from the decomposition $\fg_{1} = [\fg_{0}, x_s] \oplus \fZ_{x,1} \oplus \fh_{x,1}$, the isomorphism $[\fg_{0}, x_s]\cong \fg_{0} / Z_{\fg_0}(x_s)$, and noting that  $x$ being $G_0$-distinguished implies that $Z_{\fg_0}(x_s)=Z(\fg)_0\oplus \fh_{x,0}$.
 By  \autoref{prop:supercuspidal-dim} we have $\dim\fh_{x,0} \le \dim\fh_{x,1}$ and therefore
   \beqn
  \dim\fZ_{x,1}\le \dim\fg_1-\dim\fg_0+\dim Z(\fg)_0;
   \eeqn
   the equality holds if and only if $(\rc, \chi)$ is nil-supercuspidal. Since $\dim\cS-\dim\rc=\dim f(\cS)=\dim \fZ_{x,1}$, the proposition follows.
    \end{proof}

 \section{Distinguishedness of support (Proof of \autoref{thm:cuspidal-distinguished})}\label{sec:support}

In this section we prove \autoref{thm:cuspidal-distinguished}. We will begin with proving \autoref{prop:annihlation-distinguished} for anti-orbital sheaves annihilated by all proper parabolic restrictions in~\autoref{ssec:supercuspidal dis}, from which \autoref{thm:cuspidal-distinguished} for supercuspidal character sheaves follows immediately. In this case the arguments are much simpler and better illustrate the ideas behind the proof for the general case.    Since a cuspidal character sheaf $\cF$ may not be annihilated by parabolic restrictions in general, we construct in~\autoref{sec:approx} an approximation of $\cF$ by a complex which vanishes under parabolic restrictions up to a given perverse cohomological degree, which allows us to run arguments similar to the supercuspidal case.

\subsection{Supercuspidal case (proof of \autoref{prop:annihlation-distinguished})}\label{ssec:supercuspidal dis}

Recall from \autoref{stratification} that $\cS_x\subseteq \fg_1$  denotes the stratum containing $x$. 
We will make use of the following two lemmas. The first one is an analogue of ~\cite[Lemma 2.9]{L84}:
\begin{lemm}\label{Lem:Lus84-2.9}   
 Let $P\subseteq G$ be a $\theta$-stable parabolic subgroup with $\theta$-stable Levi factor $L$ and unipotent radical $U$. Then, for any $x\in \fl_1$, the following holds:
 \begin{enumerate}
    \item 
        For every $x'\in x + \fu_1$, we have $x\in \overline{\cS_{x'}}$.
     \item
        $\cS_x\cap (x+\fu_1)=U_0x$.
 \end{enumerate}
\end{lemm}  	
\begin{proof}   
Let $\lambda:\Cc\ra G_0$ be a cocharacter such that $\fp=\li^{\lambda}_{\ge 0}{\fg}$ with Levi $\fl=\li^{\lambda}_0{\fg}$ and $\fu=\li^{\lambda}_{>0}{\fg}$. 
For any $x'\in x+\fu_1$,  we have $\lim_{t\to 0} \Ad_{\lambda(t)} x' = x$, and hence $x$ is contained in the closure of the $L_0$-orbit of $x'$. As $x'$ and $x$ have the same image in the adjoint quotient $\cS /\!/G_0$, either they are conjugate under $G_0$ or we have $\dim Z_{G_0}(x')<\dim Z_{G_0}(x)$. In the latter case, $x'$ lies in a bigger stratum than $\cS_x$. It follows that $\cS_x\cap (x + \fu_1)\subseteq G_0x$. \par
We claim that for any integer $d>0$ and any $x'\in G_0x\cap (x+\li^{\lambda}_{\ge d}{\fg}_1)\subseteq G_0x\cap (x+\fu_1)$, there exists $u\in U_0$ such that $\Ad_ux'\in x+\li^{\lambda}_{>d}{\fg}_1$. Since $\fu=\li^{\lambda}_{\ge 1}{\fg}$, this claim with induction on $d$ proves the lemma. Recall that $G_0x\subseteq\fg_1$ is smooth and its tangent space at $x$ is $[\fg_0,x]$. The intersection of the orbit $G_0x$ with $x+\fu_1$  has tangent space $[\fg_0,x]\cap\fu_1=[\fu_0,x]=[\li^{\lambda}_{\ge 1}{\fg}_0,x]$ at $x$.
    	
Let us pick $x'\in G_0x \cap (x +\li^{\lambda}_{\ge d}{\fg}_1)$ which does not lie in $ x+\li^{\lambda}_{>d}{\fg}_1$. Note that $\Ad_{\lambda(t)}x'$ is a curve in $G_0x\cap(x+\fu_1)$ which approaches $x$ as $t\ra 0$. The tangent direction of the curve at $t=0$, when appropriately scaled,  is $\pi_d(x')$, the projection of $x'$ to $\li^{\lambda}_d{\fg}$. From the previous paragraph we conclude that $\pi_d(x')\in [\fu_0,x]$. Suppose $\pi_d(x')=[v,x]$ for some $v\in\fu_0$. Let $u:=\exp(v)\in U_0$. Then $\Ad_ux\in x+\li^{\lambda}_{>d}{\fg}_1$ as desired.
 \end{proof}

\begin{lemm}\label{lem:semisimple}
Let $\Db_{G_0}(\fg^{\nil}_1)_{\mathrm{sc}}$ denote the sum of the blocks of $\Db_{G_0}(\fg^{\nil}_1)$ containing a supercuspidal orbital sheaf. Then, $\Perv_{G_0}(\fg^{\nil}_1)_{\mathrm{sc}} = \Db_{G_0}(\fg^{\nil}_1)_{\mathrm{sc}}\cap \Perv_{G_0}(\fg^{\nil}_1)$ is a semisimple abelian category.
\end{lemm}
\begin{proof}
By the orthogonal decomposition \autoref{thm:orthogonal-decomposition}, it suffices to show that $\Perv_{G_0}(\fg^{\nil}_1)_{\xi} = \Db_{G_0}(\fg^{\nil}_1)_{\xi}\cap\Perv_{G_0}(\fg^{\nil}_1)$ is a semisimple category for every LY cuspidal support $\xi$ such that $\Perv_{G_0}(\fg^{\nil}_1)_{\xi}$ contains a supercuspidal orbital sheaf. Let $\xi$ be such an LY cuspidal support and $\cG\in \Perv_{G_0}(\fg_1^{\nil})_\xi$ a supercuspidal orbital sheaf. By \autoref{characterisation of supercuspidal}, there is a clean $G_0$-equivariant local system $\cL$ on a distinguished $G_0$-orbit $j:\cO\hookrightarrow \fg_1^{\nil}$ such that $\cG\cong j_! \cL[\dim \cO]$. Calculating in the equivariant derived category, we have:
\[
    \Hom^{\bullet}(\cG, \cG) \cong \Hom^{\bullet}(\cL, \cL) \cong (\rmH_{Z_{G_0}(z)^{\circ}}^{\bullet}\otimes \cL_z^*\otimes \cL_z)^{Z_{G_0}(z)/Z_{G_0}(z)^{\circ}} \quad \text{for $z\in \cO$}.
\]
 Since the last term vanishes in odd degrees, we have $\Hom(\cG, \cG[1]) = 0$. This implies that $\Perv_{G_0}(\fg^{\nil}_1)_{\xi}$ is semisimple because $\cG$ is the only simple object of $\Perv_{G_0}(\fg^{\nil}_1)_{\xi}$ up to isomorphism (see \cite[Proposition 4.3.3]{Liu24}). 
\end{proof}

\begin{proof}[Proof of \autoref{prop:annihlation-distinguished}]
Let $x\in \cS$ and let $P\subseteq G$ be any $\theta$-stable parabolic subgroup with $\theta$-stable Levi factor $L$ such that $x\in \fl_1$. We show that $P = G$. Let $U$ be the unipotent radical of $P$. Then, \autoref{Lem:Lus84-2.9} implies that $\supp\cF\cap (x + \fu_1) = U_0x$.  Since $U_0$ is unipotent, this implies that $\supp\cF\cap (x + \fu_1)$ is an affine space. It follows that 
\[
\RGamma_c(x + \fu_1, \cF|_{x + \fu_1}) \cong \cF_x[-2\dim U_0 x] \neq 0.
\]
On the other hand, 
\[
\cF_x\cong (\Res^{\fg_1}_{\fp_1}\cF)_{\bar x}[\dim \fu_0 - \dim \fu_1];
\]
the latter term is zero unless $P = G$ by the assumption on $\cF$. Therefore, $P = G$ as claimed. We conclude that $\cS$ consists of $G_0$-distinguished elements. \par
Let $\cF' = \cF|_{x_s + \fh_{x,1}^{\nil}}[\dim\fh_{x,1} - \dim\fg_1]$, regarded as a sheaf on $\fh_{x,1}^{\nil}$ via translation by $x_s$. By \autoref{coro:restriction-biorbital}, $\cF'$ is bi-orbital and perverse. We claim that $\cF'$ is annihilated by parabolic restriction along every $\theta$-stable proper parabolic subgroup of $H_x$. Let $Q\subseteq H_x$ be a $\theta$-stable proper parabolic subgroup. Choose a cocharacter $\lambda: \Cc \to H_{x,0}$ such that $Q = \li^{\lambda}_{\ge 0}{H}_x$ and let $V = \li^\lambda_{>0}{H}_x$ be the unipotent radical and $M=  (H_x)^{\lambda}$ the Levi factor. Put $P = \li^{\lambda}_{\ge 0}{G}$ and $U=  \li^{\lambda}_{> 0}{G}$. For each  $y\in \fq_1\cap \fh^{\nil}_{x,1}$, consider the following map :
    \[
        j_y:U_0\times^{V_0}(y + \fv_1)\to x_s + y + \fu_1,\quad  [u : y + z]\mapsto \Ad_u(x_s + y+z).
    \]
    It is a $U_0$-equivariant isomorphism by \autoref{lem:Ux} (applied to $x_s + y$). Since $U_0/ V_0$ is an affine space, the map $j_y$ induces an isomorphism
    \[
        \RGammac(y + \fv_1, \cF) \cong 
        \RGammac(x_s + y + \fu_1, \cF|_{x_s + y + \fu_1})[2\dim U_0 / V_0].
    \]
    The assumption on $\cF$ implies $\Res^{\fg_1}_{\fp_1} \cF = 0$. Therefore, $\Res^{\fh_{x,1}}_{\fq_1}\cF' = 0$. By the adjunction \eqref{eqn-adjunctions-parabolic}, this implies that every simple quotient of $\cF'$ is bi-orbital and cuspidal, hence supercuspidal by \autoref{Prop3}. Then, \autoref{lem:semisimple} implies that $\cF'$ is a direct sum of bi-orbital supercuspidal sheaves. In particular, \autoref{characterisation of supercuspidal} implies that $\cF'|_{H_0 x_n}[-\dim H_0 x_n]$ is a clean local system. Therefore, $\cF|_{\rc}[-\dim \cS]$ is clean, where $\rc = G_0x$. This proves part (ii). \par

    Finally, given any simple constituent $\chi\subseteq \cF|_{\rc}[-\dim \cS]$, the pair $(\rc, \chi)$ is a nil-supercuspidal datum, so \autoref{cor:dim est-biorbitalcuspidal} implies $\dim \cS/\!/G_0 = \dim \fg_1 - \dim \fg_0 + \dim Z(\fg)_0$, completing the proof of part (i).   
\end{proof}

\subsection{Semi-orthogonal approximation}\label{sec:approx} 
In the rest of this section, $\cF$ will be a cuspidal character sheaf on $\fg_{  1}$. We write $\Perv^{\nil}_{G_{0}}(\fg_{  1})^{\mathrm{old}}$ for the full subcategory of $\Perv^{\nil}_{G_{0}}(\fg_{  1})$ spanned by objects whose composition factors are non-cuspidal character sheaves (see~\autoref{ssec:cuspidal}).\par

We begin by constructing a semi-orthogonal approximation $\{\cF^{(n)}\}_{n\in\N}$ of $\cF$ as follows.

\begin{prop}[semi-orthogonal approximation]\label{lem:new-approx} There exist a complex $\cF^{(n)}\in \mathrm{D}^{\mathrm{b}}_{G_{0}}(\fg_{  1})$ and a morphism $\alpha_{n+1}:\cF^{(n+1)}\to \cF^{(n)}$,  for each $n\in \N$, satisfying the following conditions:

\begin{enumerate}
\item
$\cF^{(0)} = \cF$
\item
$\Hom(\cF^{(n)}, \cG[k]) = 0$ for $\cG\in \Perv^{\nil}_{G_{0}}(\fg_{  1})^{\mathrm{old}}$ and for $k\le n$ 
\item
$\mathrm{Cone}(\alpha_{n+1})[-n-1] \in \Perv^{\nil}_{G_{0}}(\fg_{  1})^{\mathrm{old}}$.
\end{enumerate}
\end{prop}

The proof relies on the following lemmas:

\begin{lemm} \label{prop:PervG} The category $\Perv^{\nil}_{G_{0}}(\fg_{  1})^{\mathrm{old}}$ is equivalent to $A\mof$ for some finite-dimensional algebra $A$.
\end{lemm}
\begin{proof} 
Since $G_{0}$ acts on $\fg^{\nil}_{-1}$ with finitely many orbits, \cite[\S 4.3]{V94} shows that there exists a projective generator $\cP\in \Perv_{G_{0}}(\fg^{\nil}_{-  1})$. Via the Fourier--Sato transform $\FS:\Perv_{G_{0}}(\fg^{\nil}_{-  1})\xrightarrow{\sim} \Perv^{\nil}_{G_{0}}(\fg_{  1})$, we have an equivalence $\Perv^{\nil}_{G_{0}}(\fg_{  1})\cong B\mof$ with $B =  \End(\cP)^{\mathrm{op}}$. Since $\Perv^{\nil}_{G_{0}}(\fg_{  1})^{\mathrm{old}}$ is a Serre subcategory of $\Perv^{\nil}_{G_{0}}(\fg_{  1})$, there exists an idempotent $e\in B$ and an equivalence $\Perv^{\nil}_{G_{0}}(\fg_{  1})^{\mathrm{old}} \cong (B/BeB)\mof$.
\end{proof}

\begin{lemm}\label{lem:modification} Let $\cK\in \mathrm{D}^{\mathrm{b}}_{G_{0}}(\fg_{  1})$ be a complex whose perverse cohomology sheaves lie in $\Perv^{\nil}_{G_{0}}(\fg_{  1})$ in all degrees. Suppose that $\Hom(\cK, \cG[k]) = 0$ for $\cG\in \Perv^{\nil}_{G_{0}}(\fg_{  1})^{\mathrm{old}}$ and for $k< 0$. Then there exists $\cK'\in \Perv^{\nil}_{G_{0}}(\fg_{  1})^{\mathrm{old}}$  
and $\varphi\in \Hom(\cK, \cK')$ satisfying the following conditions: for every $\cG\in \Perv_{G_{0}}^{\nil}(\fg_{  1})^{\mathrm{old}}$
\begin{enumerate}
\item
$\Hom(\varphi, \cG):\Hom(\cK', \cG)\to \Hom(\cK, \cG)$ is bijective;
\item
$\Hom(\varphi, \cG[1]):\Hom(\cK', \cG[1])\to \Hom(\cK, \cG[1])$ is injective.
\end{enumerate}
\end{lemm}
\begin{proof} 
Define the following functor:
\[
T: \Perv^{\nil}_{G_{0}}(\fg_{  1})^{\mathrm{old}} \to \mathrm{Vect},\quad T(\cG) = \Hom(\cK, \cG).
\]
By our assumption on $\cK$, the functor $T$ is left exact. \autoref{prop:PervG} shows that $\Perv^{\nil}_{G_{0}}(\fg_{  1})^{\mathrm{old}}$ has enough injective objects, so that \cite[\S 2.4]{MV87} is applicable to $T$; it follows that $T$ is representable by some object $\cK'\in \Perv^{\nil}_{G_{0}}(\fg_{  1})^{\mathrm{old}}$ via certain $\varphi'\in T(\cK') = \Hom(\cK, \cK')$.
Then $(\cK', \varphi)$ satisfies the requirements: the first condition is due to the fact that $(\cK', \varphi)$ represents $T$; the second condition follows from the first one by d\'evissage. 
\end{proof}

\begin{proof}[Proof of \autoref{lem:new-approx}] For convenience, we set $\alpha_0:\cF\to 0$. We construct $\cF^{(n)}$ and $\alpha_n$ by induction on $n\in \N$. The case $n = 0$ is trivial. Let $n \geq 0$ and suppose that $\cF^{(n)}$ and $\alpha_{n}$ have been constructed. Applying \autoref{lem:modification} to $\cK_n := \cF^{(n)}[-n-1]$, we obtain $(\cK_n',\varphi_n)$. We define $\cF^{(n+1)} = \Cone(\varphi_n)[n]$
so that it fits into an exact triangle
\beq\label{eqn-triangle}
\cF^{(n+1)}\xrightarrow{\alpha_{n+1}} \cF^{(n)}\xrightarrow{\varphi_n[n+1]} \cK_n'[n+1]\xrightarrow{[1]}.
\eeq
The conditions on $(\cK', \varphi)$ imply that $\Hom(\cF^{(n+1)}, \cG[k]) = 0$ for $\cG\in \Perv^{\nil}_{G_{0}}(\fg_{  1})^{\mathrm{old}}$ and $k\le n+1$. 
\end{proof}
It is easy to see from the construction in~\eqref{eqn-triangle} that 
\begin{subequations}
\label{semiorthogonal property}
\begin{equation} \label{semiorthogonal property-a}
\text{$\pH^k\cF^{(n)} = 0$ unless $k\in [-n+1,0]$ or $k = n = 0$;}
\end{equation}
\begin{equation} \label{semiorthogonal property-b}\text{
$
\pH^k(\alpha_{n+1}):\pH^k\cF^{(n+1)}\to \pH^k\cF^{(n)}
$}
\text{ is an isomorphism unless $k = -n$;}
\end{equation}
\begin{equation} \label{semiorthogonal property-c}\text{
$
\pH^{-n}(\alpha_{n+1}):\pH^{-n}\cF^{(n+1)}\to \pH^{-n}\cF^{(n)}
$}
\text{ is surjective with kernel in $\Perv^{\nil}_{G_{0}}(\fg_{  1})^{\mathrm{old}}$.}
\end{equation}
\end{subequations}

Thus, \autoref{prop:stratum-support} implies that $\{\supp \cF^{(n)}\}_{n\in \N}$ is an increasing sequence of Zariski-closed subsets which are unions of strata of $\fg_1$. It follows that the sequence stabilises for $n \ge n_0$ for some $n_0\in \N$. From now on we fix such an $n_0\in \N$ and let $\cS$ be a generic stratum of $\supp \CF^{(n)}$, $n\geq n_0$, see~\autoref{def-generic strata}.  \par

We will now state two propositions and deduce \autoref{thm:cuspidal-distinguished} from them in the next subsection. The propositions will be proved  in subsequent subsections.

\begin{prop}\label{lem:supp-distinguished-new} 
Every element of $\cS$ is $G_0$-distinguished.
    
\end{prop} 

To state the next proposition let $x=x_s+x_n\in\cS$. We consider the restriction to a generic slice $\cF^{(n)}|_{x_s + \fh_{x, 1}^{\nil}}$ studied in \autoref{ssec:restriction-sheaf}. \autoref{prop:localtriv} implies that for every $n\ge n_0$:
\[
    \pH^{\bullet}(\cF^{(n)}|_{x_s + \fh_{x, 1}^{\nil}}) \cong \pH^{\bullet + \delta}(\cF^{(n)})|_{x_s + \fh_{x, 1}^{\nil}}[-\delta],\quad \delta = \codim_{\fg_1}(\fh_1).
\]
Therefore, recalling~\eqref{semiorthogonal property}, there exists an integer $a$ such that for $n\ge n_0$:
\[
    \pH^a(\CF^{(n)}|_{x_s + \fh^{\nil}_{x,1}})\neq 0\quad \text{and} \quad
    \pH^{>a}(\CF^{(n)}|_{x_s + \fh^{\nil}_{x,1}})= 0.
\]
    
\begin{prop}\label{lem:support-reduction-cuspidal} There exists $n_1 \ge n_0$ such that $\pH^{a}(\CF^{(n)}|_{x_s + \fh^{\nil}_{x,1}})$ is a non-zero direct sum of bi-orbital supercuspidal character sheaves on $\fh_{x,1}$, for all $n\geq n_1$. 
  \end{prop}

 \subsection{Proof of \autoref{thm:cuspidal-distinguished}} The proof can be easily reduced to the case $Z(G)_0 = 1$ and we will make that assumption. Let $\cS$ be a generic stratum of $\supp\cF^{(n)}$, $n\geq n_0$. Let $x=x_s+x_n\in\cS$ and set $\rc=G_0x$. \autoref{lem:supp-distinguished-new} implies that $\rc$ is $G_0$-distinguished. 
 Let $n_1\in \N$ and $a\in \Z$ be as in \autoref{lem:support-reduction-cuspidal}  so that $\pH^{>a}(\CF^{(n)}|_{x_s + \fh^{\nil}_{x,1}}) = 0$ and $\pH^a(\CF^{(n)}|_{x_s + \fh^{\nil}_{x,1}})$ is a non-zero sum of bi-orbital supercuspidal sheaves on $\fh_{x,1}$ for every $n\ge n_1$. Applying \autoref{characterisation of supercuspidal}, we conclude that the restriction $\rmH^{a'}(\cF^{(n)}|_{\rc})\neq 0$ is a clean local system for $a' = a - \dim {H_0}x_n$. 
 
\par
The next step is to pass from $\cF^{(n)}$ for $n\geq n_1$ to $\cF^{(0)}=\cF$. To that end, consider the morphism 
\[
    \beta: \pH^{a}(\cF^{(n)}|_{x_s + \fh_{x,1}^{\nil}}) \to \pH^{a}(\cF^{(n-1)}|_{x_s + \fh_{x,1}^{\nil}})
\]
induced by $\alpha_{n}:\cF^{(n)}\to \cF^{(n-1)}$ for all $n\geq 1$. By \eqref{semiorthogonal property}, it is an isomorphism when $a \neq 1-n -\delta$ and surjective otherwise. As $\pH^{a}(\cF^{(n)}|_{x_s + \fh_{x,1}^{\nil}})$ is a sum of  bi-orbital supercuspidal sheaves on $\fh_{x,1}$ for every $n\ge n_1$ and $\beta$ is always a surjection, we conclude that this also holds for any $n\geq 1$. \par

We will argue next that $\beta$ is an isomorphism for $n\geq 1$. We only need to consider the case $a = 1 - n - \delta$, and let us now suppose that $\beta$ is not injective in this case. Set $\cK = \ker(\pH^{1-n}(\alpha_{n}))\in \Perv^{\nil}_{G_0}(\fg_1)^{\mathrm{old}}$. As $\pH^{a}(\cF^{(n)}|_{x_s + \fh_{x,1}^{\nil}})$ is a direct sum of bi-orbital supercuspidal sheaves on $\fh_{x,1}$ for every $n\ge 1$, we conclude that $\ker(\beta) = \cK|_{x_s + \fh_{x, 1}^{\nil}}[-\delta]$ is a non-zero sum of bi-orbital supercuspidal sheaves on $\fh_{x,1}$. This implies that $\cS$ is a generic stratum of $\supp \cK$. Let $\chi$ be an irreducible direct factor of $\cK|_{\rc}[-\dim \cS]$. We conclude that $(\rc, \chi)$ is a nil-supercuspidal datum. Then, the cleanness of $\chi$,~\autoref{prop:supercuspidal} and the adjunction~\eqref{eqn-adjunctions-parabolic} imply that 
\[
\Hom(\cK|_{\rc}, \chi[\dim \cS]) = \Hom(\cK, \IC(\chi)[\dim \cS-\dim \rc]) = 0, 
\]
which contradicts the assumption that $\chi$ is a factor of $\cK|_{\rc}[-\dim \cS]$. Thus, we have argued that  $\beta$ is an isomorphism for $n\geq 1$.  
        
From the discussion above, we conclude that $\alpha_{1}\cdots \alpha_{n}$ induces an isomorphism \linebreak $\pH^{a}(\cF^{(n)}|_{x_s + \fh_{x,1}^{\nil}}) \cong \pH^{a}(\cF|_{x_s + \fh_{x,1}^{\nil}})$ for any $n\geq 1$  so that $\CF_{x} \neq 0$ and thus $\cS\subseteq \supp \CF$. Since this holds for every generic stratum $\cS\subseteq \supp \CF^{(n)}$, we see that $\supp \CF^{(n)} = \supp \CF$. In particular, \autoref{lem:supp-distinguished-new} implies that every element of the supporting stratum of $\CF$ is $G_0$-distinguished. Moreover, \autoref{coro:restriction-biorbital} implies that the restriction $\CF\mid_{x_s + \fh^{\nil}_{x,1}}[-\delta]$ is a sum of bi-orbital supercuspidal sheaves.

 Let $\cS\subseteq \supp \cF$ denote the supporting stratum of $\cF$. Finally, applying~\autoref{cor:dim est-biorbitalcuspidal}, we conclude that $\dim f(\cS)=\dim\fg_1-\dim\fg_0$. This completes the proof of~\autoref{thm:cuspidal-distinguished}. \qed

\begin{rema}
    The above proof actually shows that $f(\supp\Cone(\alpha_{n}))\subsetneq f(\supp\cF)$ for every $n\ge 0$. Namely, the semi-orthogonal approximation $\{\cF^{(n)}\}_{n\in \N}$ proceeds by successively attaching non-cuspidal perverse sheaves whose support has \emph{strictly smaller} semisimple part than $\supp\cF$. 
\end{rema}

\subsection{Proof of \autoref{lem:supp-distinguished-new}} 
 We have the following vanishing property:
\begin{lemm}\label{lem:parab-res-new} There exists $d\ge 0$ such that 
        \[
            \pH^{\ge d-n}\Res^{\fg_1}_{\fp_1} \CF^{(n)} = 0
        \]
        holds for every $n\in \N$ and every $\theta$-stable proper parabolic subgroup $P$. 
        \end{lemm}
\begin{proof} The adjunction~\eqref{eqn-adjunctions-parabolic} implies that the number $d$ can be chosen to be the maximum of perverse cohomological amplitude of parabolic inductions. 
\end{proof}

 We now proceed to prove   \autoref{lem:supp-distinguished-new}.
     Suppose on the contrary that an element $x\in\cS$ is not $G_0$-distinguished. Then there exists a proper $\theta$-stable parabolic subgroup with $\theta$-stable Levi decomposition $P=LU$ such that $x\in\fl_1$. By \autoref{Lem:Lus84-2.9}, we have 
     \[
      \supp \cF^{(n_0)}\cap (x + \fu_1) = G_0x\cap (x + \fu_1) = U_0x,
      \]
     and the last term is isomorphic to an affine space by the unipotency of $U_0$; let $r$ denote its dimension. We have for $n\ge n_0$:
    \[
        (\NRes^{\fg_1}_{\fp_1}\CF^{(n)})_x = \RGammac(x+\fu_1, \cF^{(n)}|_{x+\fu_1})[\dim\fu_1 - \dim\fu_0] \cong \cF^{(n)}_x[\dim\fu_1 - \dim\fu_0-2r].
    \]
     We choose an $e\in \Z$  such that $\rmH^e(\CF^{(n_0)}_x) \neq 0$ and $\rmH^{>e}(\CF^{(n_0)}_x) = 0$. The properties in \eqref{semiorthogonal property} and the genericity of $\cS$ imply that $\rmH^e(\CF^{(n)}_x) \neq 0$ and $\rmH^{>e}(\CF^{(n)}_x) = 0$, for $n\ge n_0$. On the other hand, \autoref{lem:parab-res-new} implies that there exists $d\ge 0$ such that $\pH^{\ge d - n}\NRes^{\fg_1}_{\fp_1} \CF^{(n)} = 0$ for $n \ge n_0$. Let $n$ be a large integer (larger than $n_0$ and $d - e + \dim\fu_1 - \dim\fu_0-2r$). It follows that we have
     \[
     \pH^{\ge e - \dim\fu_1 + \dim\fu_0+2r}\NRes^{\fg_1}_{\fp_1} \CF^{(n)} = 0
     \]
     and thus
     \[
     \rmH^e(\cF^{(n)}_x) \cong \rmH^{e - \dim\fu_1 + \dim\fu_0+2r}(\NRes^{\fg_1}_{\fp_1} \CF^{(n)})_x = 0,
     \]
a     contradiction. Therefore, the element $x\in \cS$ must be $G_0$-distinguished.  \qed

\subsection{Proof of \autoref{lem:support-reduction-cuspidal}}We begin with some lemmas.

\begin{lemm}\label{lem:criterion-cuspidal}
Let $a\in \Z$ and let $\cK\in \Dbn_{G_0}(\fg_1)$. Suppose that $\pH^{n} \cK = 0$ for $n > a$ and that $\pH^n\Res^{\fg_1}_{\fp_1}\cK = 0$ for every proper $\theta$-stable parabolic subgroup $P\subseteq G$ and for $n \ge a - d$, where $d$ is the maximal perverse cohomological amplitude of parabolic inductions. Then, the cosocle of $\pH^a \cK\in \Perv^{\nil}_{G_0}(\fg_1)$ is a direct sum of cuspidal character sheaves.
\end{lemm}
\begin{proof}
Let $\cJ$ denote the cosocle of $\pH^a\cK$, and let $\cL\in \Char_{G_0}(\fg_1)^{\mathrm{old}}$ be a non-cuspidal character sheaf. We need to show that $\Hom(\cJ, \cL)= 0$. We have  natural morphisms 
\[
\cK[a]\to \pH^a\cK\to \cJ. 
\]
Set $\cJ' = \Cone(\cK[a]\to \cJ)$. Then, by construction, $\pH^{n}\cJ' = 0$ for $n\ge 0$ and we have an exact sequence
\[
0= \Hom(\cJ', \cL)\to \Hom(\cJ, \cL)\to\Hom(\cK[a], \cL)\,.
\]
On the other hand, $\cL$ is a direct summand of $\pH^{n'}\Ind^{\fg_1}_{\fp_1}\cL'$ for a $\theta$-stable proper parabolic subgroup $P\subseteq G$, a $\theta$-stable Levi factor $L$, a character sheaf $\cL'\in \Char_{L_0}(\fl_1)$ and an integer $n'\in [-d, d]$; hence, we get
\[
\Hom(\cK[a], \cL)\subseteq \Hom(\cK[a], \Ind^{\fg_1}_{\fp_1}\cL'[n'])\cong \Hom(\Res^{\fg_1}_{\fp_1}\cK[a-n'], \cL') = 0
\]
by our assumption. It follows that $\Hom(\cJ, \cL)\subseteq \Hom(\cK[a], \cL) = 0$ as asserted. Therefore, $\cJ$ is a direct sum of cuspidal character sheaves.
\end{proof}

 \begin{lemm}\label{lem:support-reduction-vanishing} 
     Let $b\in \Z$. Suppose that $\cK\in \Dbn_{G_0}(\fg_1)$ satisfies that $\pH^{n} \NRes^{\fg_1}_{\fp_1}\cK = 0$ for every $\theta$-stable proper parabolic subgroup $P\subseteq G$ and for every $n > b$. Let $x\in \supp \cK$ be such that the stratum $\cS_x$ is generic in $\supp \cK$. Then,  the restriction $\cK' = \cK|_{x_s + \fh^{\nil}_{x,1}}$, viewed as an object of $\Db_{H_{x,0}}(\fh^{\nil}_{x,1})$ via translation, satisfies $\pH^{n} \NRes^{\fh_{x,1}}_{\fq_1}\cK' = 0$ for every $\theta$-stable proper parabolic subgroup $Q\subseteq H_x$ and for every $n > b + c$, where $c\in \N$ is a constant depending only on $\fg$. 
\end{lemm}
\begin{proof}[Proof of \autoref{lem:support-reduction-cuspidal}]
    Let $Q\subseteq H_x$ be a $\theta$-stable proper parabolic subgroup. Choose a cocharacter $\lambda: \Cc \to H_{x,0}$ such that $Q = \li^{\lambda}_{\ge 0}{H}_x$ and let $V = \li^\lambda_{>0}{H}_x$ be the unipotent radical and $M=  (H_x)^{\lambda}$ the Levi factor. Put $P = \li^{\lambda}_{\ge 0}{G}$ and $U=  \li^{\lambda}_{> 0}{G}$. By assumption, we have $\pH^{>b}\NRes^{\fg_1}_{\fp_1}\cK = 0$ and thus 
    \[
        \rmH^{>b + c}_{\mathrm{c}}(y + \fu_1, \cK|_{y + \fu_1}) = 0
    \]
    for every $y\in \fp_1$, where $c$ is a constant depending only on $\fg$.    As in the proof of \autoref{prop:annihlation-distinguished}, the following map :
    \[
        j_y:U_0\times^{V_0}(y + \fv_1)\to x_s + y + \fu_1,\quad  [u : y + z]\mapsto \Ad_u(x_s + y+z).
    \]
    for each $y\in \fq_1\cap \fh^{\nil}_{x,1}$ induces an isomorphism
    \[
        \rmHc^{*}(y + \fv_1, j_y^*\cK') \cong 
        \rmHc^{* + 2\dim U_0 / V_0}(x_s + y + \fu_1, \cK|_{x_s + y + \fu_1}).
    \]
    Therefore, $\rmHc^{n}(y + \fv_1, j_y^*\cK') = 0$ for $n > b+ c + 2\dim U_0 / V_0$ and for every $y\in \fq_1\cap \fh^{\nil}_{x,1}$. Thus, $\pH^{n} \NRes^{\fh_{x,1}}_{\fq_1}\cK' = 0$ for $n > b + c'$ for some constant $c'$ depending only on $\fg$. 
\end{proof}

We now proceed to prove \autoref{lem:support-reduction-cuspidal}. 
It follows from \autoref{lem:parab-res-new} that for some $d_1\in \N$, we have $\pH^{\ge d_1-n} (\NRes^{\fg_1}_{\fp_1} \CF^{(n)}) = 0$ for all $n\in \N$ and all $\theta$-stable proper parabolic $P\subseteq G$.  
Recall that $a\in \Z$ is the unique integer such that the following holds for $n \ge n_0$:
\[
\pH^{> a}(\cF^{(n)}|_{x_s + \fh_{x, 1}^{\nil}}) = 0 \quad\text{and}\quad \pH^{a}(\cF^{(n)}|_{x_s + \fh_{x, 1}^{\nil}}) \neq 0.
\]
Then, \autoref{lem:support-reduction-vanishing} implies that for some $d_2\in \N$, we have $\pH^{\ge d_2-n} (\NRes^{\fh_{x,1}}_{\fq_1}(\cF^{(n)}|_{x_s + \fh_{x, 1}^{\nil}})) = 0$ for all $n\ge n_0$ and all $\theta$-stable proper parabolic subgroups $Q\subseteq H_{x}$. We may choose $n_1 \ge n_0$ to be large enough such that $\cF^{(n)}|_{x_s + \fh_{x, 1}^{\nil}}$ satisfies the conditions of \autoref{lem:criterion-cuspidal} for $n\ge n_1$. In particular, the cosocle of $\pH^a(\cF^{(n)}|_{x_s + \fh_{x, 1}^{\nil}})$ is a sum of bi-orbital cuspidal character sheaves. By \autoref{Prop3}, they are supercuspidal.  
Applying  \autoref{lem:semisimple} below to the category $\Db_{H_{x,0}}(\fh^{\nil}_{x,1})$ shows that $\pH^a(\CF^{(n)}|_{x_s + \fh_{x, 1}^{\nil}})$ is semisimple. This completes the proof.\qed

\section{Proof of the Adjunction formula (\autoref{thm-adjunction})}
\label{sec:proof-adjunction}

In this section, we prove \autoref{thm-adjunction}. We make use of the following auxiliary propositions:
 \begin{prop}\label{prop:dim-symmetry}
 
Assume that $Z(G)_0 = 1$. Let $(\rc, \chi)$ be a nil-supercuspidal datum. Then, $\Hom(\cF, \IC(\chi)[N + k]) = 0$ for every $k\in \Z \setminus\{0\}$ and every character sheaf $\cF\in \Char_{G_0}(\fg_1)$.
 \end{prop}
 
\begin{prop}
 \label{prop:hom-psi-2}
    Assume that $Z(G)_0 = 1$. Let $(\rc,\chi)$ be a supercuspidal datum, $(\ul L, \ul\fl_*, \cC)$ a LY cuspidal support on $\fg_{  1}$ (see~\S\ref{ssec:LYSupp}), $\ul\fp_*$ a $1$-spiral and $\ul\fp'_*$ a $(-1)$-spiral such that $(\ul\fp_*, \ul\fl_*, \ul\fp'_*)$ is an opposite spiral triplet (\autoref{ssec:hyperbolic-spiral}). Then, there is a canonical isomorphism
     \[
         \Hom\left(\cP_{\chi} ,\Ind_{\ul\fp_1}^{\fg_{  1}}\IC(\cC)[k]\right)\cong\Hom\left(\Ind_{\fpo'}^{\fg_{  1}}\IC(\cC), \IC(\chi)[N - k]\right)^*
     \]
    for each $k\in \Z$, where $N = \dim \fg_{  1} - \dim \fg_{  0}$.
 \end{prop}

 We will apply the above propositions to deduce \autoref{thm-adjunction}  in the next subsection and defer the proof of the propositions to later subsections.

\subsection{Proof of \autoref{thm-adjunction}}
        Let us consider the Grothendieck group of semisimple orbital complexes defined in~\cite{LY17a}. Let $\CQ$ be the full subcategory of $\Db_{G_0}(\fg^{\nil}_1)$ spanned by semisimple complexes and let $\operatorname{K}_0(\CQ)$ be the split Grothendieck group of $\CQ$. Then, $\operatorname{K}_0(\CQ)$ acquires a $\Z[v^{\pm 1}]$-module structure on which $v$ acts by cohomological shift $[-1]$. Set 
		\[
			\bV = \Q(v)\otimes_{\Z[v^{\pm 1}]} \operatorname{K}_0(\CQ).
		\] \par
  
        Consider the $\Q(v)$-linear forms:
       \[
       h,g:\bV\to \Q(v)
       \]
       defined by the following formulas on $\cG\in \CQ$ and extended $\Q(v)$-linearly to $\bV$:
       \[
       h([\cG]) = \sum_{k\in \Z}v^k \dim\Hom(\cP_\chi ,\cG[k]),\quad g([\cG]) = \sum_{k\in \Z}v^k \dim\Hom(\cG^{\dag}, \IC(\chi)[N + k]).
       \]
       
Using the notation above, the adjunction formula is equivalent to the following identity: $g([\cG]) = h([\cG])$ for $\cG\in \Irr \Perv_{G_0}(\fg_1^{\nil})$.  \par

        Let $\xi = (\ul L, \ul\fl_*, \cC)$ be an LY cuspidal support on $\fg_1$ and $\ul\fp_*$ a $1$-spiral of $\fg_*$ with splitting $\ul\fl_*$.
       Put $\cK = \Ind^{\fg_1}_{\ul\fp_{1}}\IC(\cC)$.  \autoref{prop:involution-series} implies that $\cK^{\dag} \cong \Ind^{\fg_1}_{\ul\fp'_{1}}\IC(\cC)$.  By the relative hard Lefschetz theorem~\cite[\S 6.2.10]{BBDG18}, there exists an isomorphism of semisimple anti-orbital sheaves, for each $k\ge 0$:
        \[
        \pH^{-k} \Ind^{\fg_1}_{\ul\fp_1'}\IC(\cC) \cong \pH^{k} \Ind^{\fg_1}_{\ul\fp_1'}\IC(\cC)\,.
        \]
         Making use of~\autoref{prop:dim-symmetry}, we get
        \beq\label{eqn:dim-symmetry}
        \bega
        \dim \Hom(\cK^{\dag}, \IC(\chi)[N - k] ) = \dim \Hom(\pH^{k}\cK^{\dag}, \IC(\chi)[N] ) \\
         = \dim \Hom(\pH^{-k}\cK^{\dag}, \IC(\chi)[N] )
         = \dim \Hom(\cK^{\dag}, \IC(\chi)[N + k] ) 
        \eega
        \eeq
        for every $k\in \Z$. \par
     
       From \eqref{eqn:dim-symmetry} and~\autoref{prop:hom-psi-2}, we conclude that $h([\cK]) = g([\cK])$.
       By~\cite[8.4(b)]{LY17a}, the space $\bV$ is spanned over $\Q(v)$ by the classes $\{[\Ind^{\fg_1}_{\ul\fp_1}\IC(\cC)]\}_{\xi, \ul\fp_*}$, where $\xi = (\ul L, \ul\fl_*, \cC)$ runs over all LY cuspidal supports on $\fg_1$ and $\ul\fp_*$ runs over all $1$-spirals having $\ul\fl_*$ as splitting factor, up to $G_0$-conjugation. It follows that $h = g$ on $\bV$ as claimed, which proves the first assertion. \par
       Given $\cG\in \Irr \Perv_{G_0}(\fg_1^{\nil})$, the vanishing $\Hom(\cG^{\dag}, \IC(\chi)[N + k]) = 0$ for $k\notin\Z$ follows from~\autoref{prop:dim-symmetry}. This proves the second assertion. \qed

 \subsection{Proof of \autoref{prop:dim-symmetry}}
 Let $\cF\in \Char_{G_0}(\fg_1)$ be a character sheaf. If $\cF$ is non-cuspidal, then~\autoref{prop:supercuspidal} and the adjunction~\eqref{eqn-adjunctions-parabolic} imply that $\Hom(\cF, \IC(\chi)[k]) = 0$ for $k\in \Z$ and thus the statement holds. Assume that $\cF$ is cuspidal. Let $\widecheck\cO\subseteq\fg_1$ be the supporting stratum of $\cF$ and $\cS\subseteq \fg_1$ the stratum containing $\rc$. \autoref{thm:cuspidal-distinguished}(i) and \autoref{cor:dim est-biorbitalcuspidal} imply that $\dim f(\widecheck\cO) = \dim \fg_1 - \dim \fg_0 = \dim f(\cS)$. If $f(\widecheck\cO) \neq f(\cS)$, then $\rc\cap \supp \cF = \emptyset$, so the cleanness of $\chi$ implies 
 \[
     \Hom(\cF, \IC(\chi)[k]) = \Hom(\cF|_{\rc}, \chi[\dim \rc + k]) =0
 \]
 for $k\in \Z$. Therefore, we may assume that $f(\cS) = f(\widecheck\cO)$.  This implies that for any given $x\in \rc$ we have $\widecheck\cO\cap (x_s + \fh_{x,1}^{\nil}) \neq \emptyset$ . It follows from \autoref{thm:cuspidal-distinguished}(ii) that the restriction $\cF|_{x_s + \fh_{x, 1}^{\nil}}[-\codim_{\fg_1}\fh_{x,1}]$ is a sum of bi-orbital supercuspidal sheaves. This implies by \eqref{restriction} that $\cF|_{f^{-1}f(x)}[-N]$ is also a sum of clean perverse sheaves and thus $\cF|_{\rc}$ is a local system concentrated in degree $-N -\dim \rc$. Now,
 \[
\Hom(\cF, \IC(\chi)[N + k])  = \Hom(\cF|_{\rc}[-\dim \rc -N], \chi[k])=0 \ \ 
\text{for}\ \   k\neq 0\]
 because,  as $x$ is $G_0$-distinguished and $Z(G)_0 = 1$, the category $\Db_{G_0}(\rc)$ is semisimple.\qed

\subsection{Double cover}\label{ssec:double-cover} 
In this subsection, we introduce a slight variant of the nearby-cycle construction by passing to the double cover. Using this description, the monodromy becomes unipotent, and this will allow us to use Beilisnson's unipotent nearby-cycle construction, which we review in the next subsection. 

Let $\rc$ be a $G_0$-orbit and $\chi$ a $G_0$-equivariant local system on $\rc$. Let $x\in \fg_{  1}$. Recall the objects $\cZ$, $\cL_{\chi}$ and $\cP_{\chi}$ introduced in~\eqref{eqn-cz}. Set 
     \[
      \cZ' = \left\{ (z, r)\in \fg_{  1}\times \C\;|\; f(z) = r^2f(x) \right\}.
     \]
     This is nothing but the base change of $\cZ$ along the two-fold ramified cover $\C\to \C: r\mapsto r^2$. The variety $\cZ'$ has a $(G_{  0}\times \Cc)$-action given by $(g, t)(z, r) = (t^{-2}\Ad_g z, t^{-1}r)$. There is a $(G_{  0}\times\Cc)$-equivariant finite morphism $\cZ'\to \cZ$ given by $(z, r)\mapsto (z, r^2)$. Let $\cL'_{\chi}\in \Perv_{G_0}(\cZ')$ denote the inverse image of $\cL_{\chi}$ under this map. We write $$\pi':\cZ'\to \C$$ for the projection to the second factor. As there is a canonical isomorphism $\psi_{\pi'}\cL'_{\chi}\cong \psi_\pi \cL_{\chi}=\cP_{\chi}$, we can replace $\psi_\pi \cL_{\chi}$  with  $\psi_{\pi'}\cL'_{\chi}$. \par 
 \begin{lemm}\label{lem:monodromy-unipotent}
 The monodromy action on $\psi_{\pi'}\cL'_{\chi}$ is unipotent.
 \end{lemm}
 \begin{proof}
The exponential map $\exp:\C\to \Cc$ induces an action of $G_0\times \C$ on $\cZ'$, given by $(g, t)(z, r) = (e^{-2t}\Ad_g z, e^{-t}r)$. Similarly, let $\C$ act on $\C$ by $t\cdot r = e^{-t}r$, so that the projection $\pi'$ is $\C$-equivariant. Then the local system $\cI$ on $\Cc$ associated with the regular representation of $\pi_1(\Cc, 1)$ is $\C$-equivariant. This induces a $(G_0\times \C)$-equivariant structure on the complex $\scrHom(\cI_{\pi'}, \cL'_{\chi})$. It follows that the perverse sheaf $\psi_{\pi'}\cL'_{\chi}$ is $(G_{  0}\times \C)$-equivariant. The subgroup $1\times 2\pi i\Z\subseteq G_{  0}\times \C$ acts trivially on the variety $\cZ'$ and yields a homomorphism $2\pi i \Z \to \Aut(\psi_{\pi'}\cL'_{\chi})$, which is by definition the monodromy action. \par 
The semisimplification $(\psi_{\pi'}\cL'_{\chi})^{\mathrm{ss}}$ is a direct sum of simple $(G_{  0}\times \C)$-equivariant perverse sheaves on $\pi'^{-1}(0) \cong \fg_{  1}^{\nil}$. However, every $G_0$-equivariant simple perverse sheaf on $\fg_1^{\nil}$ admits a unique $(G_{  0}\times \Cc)$-equivariant enhancement (see the proof of \autoref{lem:orbital-monodromic}). In other words, the monodromy action on $(\psi_{\pi'}\cL'_{\chi})^{\mathrm{ss}}$ is trivial, which means that the monodromy action on $\psi_{\pi'}\cL'_{\chi}$ is unipotent.
 \end{proof}

\subsection{Beilinson's construction}\label{ssec:beilinson}

To prove~\autoref{prop:hom-psi-2}, it is convenient for us to use Beilinson's construction of the unipotent nearby cycles \cite{B87}, which we review now. 
Let $\cZ$ be a complex algebraic variety equipped with a function $\pi\in \Gamma(\cZ, \cO_{\cZ})$, which defines a morphism $\pi:\cZ\to \C$.  Let $j:\cU = \pi^{-1}(\Cc) \hookrightarrow \cZ$ and $i:\pi^{-1}(0) \hookrightarrow \cZ$ be the inclusions. We write $\psi^{\mathrm{uni}}_\pi$ for the unipotent part of the nearby-cycle functor $\psi_\pi$.

Let $\gamma\in \pi_1(\Cc, 1)$ be the generator satisfying $\int_{\gamma} z^{-1}dz = 2\pi i$. For $n\in \N$, let $\cI^{(n)}$ be the local system on $\Cc$ corresponding to the representation $\C[\gamma^{\pm 1}] / (1 - \gamma)^n$ of $\pi_1(\Cc, 1)$. Then, there is a natural quotient map $\cI^{(a)}\to \cI^{(a)} / (1 - \gamma)^b \cong \cI^{(b)}$ for $b \le a$. Set $\cI_\pi^{(n)} = (\pi|_{\cU})^*\cI^{(n)}$. 

\begin{lemm}\label{lem:psi-unip}
There is a canonical bi-functorial isomorphism for $\cL\in \Dbc(\cU)$ and $\cK\in \Dbc(\pi^{-1}(0))$:
\[
\varinjlim_{n\in \N}\Hom\left(i^*j_*(\cL\otimes \cI^{(n)}_\pi), \cK\right) \cong \Hom(\psi^{\mathrm{uni}}_\pi \cL, \cK).
\]
\end{lemm}

\begin{proof}
We reproduce Beilinson's argument. Recall that there is a functorial exact triangle for $\cG\in \Dbc(\cU)$ (see~\cite[(8.6.7)]{KS90}):
\[
\psi^{\mathrm{uni}}_\pi\cG\xrightarrow{1 - \gamma} \psi^{\mathrm{uni}}_\pi\cG\to i^*j_* \cG \to .
\]
Let $\cL\in \Dbc(\cU)$. Applying this triangle to $\cL\otimes \cI^{(n)}_{\pi}$ for $n\in \N$, we obtain
\[
i^*j_*(\cL\otimes \cI^{(n)}_\pi) \cong \Cone(\psi^{\mathrm{uni}}_\pi(\cL\otimes \cI^{(n)}_\pi)\xrightarrow{\beta^{(n)}} \psi^{\mathrm{uni}}_\pi(\cL\otimes \cI^{(n)}_\pi)),
\]
where $\beta^{(n)} = 1 - \gamma$. We also have
\[
    \psi^{\mathrm{uni}}_\pi(\cL\otimes \cI^{(n)}_\pi) \cong (\psi^{\mathrm{uni}}_\pi\cL)\otimes \C[\gamma] / (1 - \gamma)^{n},
\]
where the monodromy operator $\gamma$ acts diagonally on the right-hand side. When $n \ge 0$ satisfies $(1 - \gamma)^n\psi^{\mathrm{uni}}_\pi\cL = 0$, we have a well-defined morphism
\[
\mu^{(n)}: (\psi^{\mathrm{uni}}_\pi\cL)\otimes \C[\gamma] / (1 - \gamma)^{n}\to \psi^{\mathrm{uni}}_\pi\cL,\quad m\otimes \gamma^k \mapsto \gamma^{-k}m.
\]
As, by construction, $\mu^{(n)}\beta^{(n)} = 0$, we obtain an induced morphism 
\[
    \bar\mu^{(n)}: \Cone(\beta^{(n)})\to \psi^{\mathrm{uni}}_\pi\cL.
\]
For $b \le a$, the quotient morphism $\cI^{(a)}\to \cI^{(b)}$ induces a morphism $\Cone(\beta^{(a)})\to \Cone(\beta^{(b)})$ which intertwines $\bar\mu^{(a)}$ and $\bar\mu^{(b)}$. 

Given $\cK\in \Dbc(\pi^{-1}(0))$,  $\bar\mu^{(n)}$ induces a map
\[
  \Hom(\psi^{\mathrm{uni}}_\pi\cL, \cK)\to \varinjlim_{n\in\N}\Hom(\Cone(\beta^{(n)}), \cK) \cong  \varinjlim_{n\in\N}\Hom(i^*j_*(\cL\otimes \cI^{(n)}_\pi), \cK). 
\]
We will show next that the first map is an isomorphism. This will complete the proof.
\par

By d\'evissage, it suffices to prove the claim for objects of $\Perv(\cU)$. Let $\cL\in \Perv(\cU)$. The source and target of $\beta^{(n)}$ are perverse by the t-exactness of the nearby-cycle functor. \par
Keeping the assumption $(1 - \gamma)^n\psi^{\mathrm{uni}}_\pi\cL = 0$, we see that $\bar \mu^{(n)}$ induces an isomorphism 
\[
\coker(\beta^{(n)}) = \pH^0(\Cone(\beta^{(n)})) \cong \psi^{\mathrm{uni}}_\pi\cL.
\]
Indeed, we have $\mu^{(n)}((\psi^{\mathrm{uni}}_\pi\cL)\otimes 1) = \psi^{\mathrm{uni}}_\pi\cL$ and one verifies easily that 
\[
(\psi^{\mathrm{uni}}_\pi\cL)\otimes\C[\gamma] / (1 - \gamma)^{n} = (\psi^{\mathrm{uni}}_\pi\cL)\otimes 1 + \on{im}(\beta^{(n)}).
\]

Moreover, for a given $b\in \N$, the map $\ker (\beta^{(a)})\to \ker (\beta^{(b)})$ is zero when $a \gg b$. For any $\cK\in \Dbc(\pi^{-1}(0))$, we have an exact sequence
 \begin{align*}
  &\varinjlim_{n\in\N}\Hom(\ker(\beta^{(n)}), \cK[-1])\to \varinjlim_{n\in\N}\Hom(\Cone(\beta^{(n)}), \cK) \to \varinjlim_{n\in\N}\Hom(\coker(\beta^{(n)}), \cK) \\
  &\to \varinjlim_{n\in\N}\Hom(\ker(\beta^{(n)}), \cK).
 \end{align*}
 By the above observations on transition maps, we have 
 \begin{align*}
  &\varinjlim_{n\in\N}\Hom(\ker(\beta^{(n)}), \cK[-1])= 0,\quad \varinjlim_{n\in\N}\Hom(\ker(\beta^{(n)}), \cK) = 0, \\
  &\varinjlim_{n\in\N}\Hom(\coker(\beta^{(n)}), \cK) \cong \Hom(\psi^{\mathrm{uni}}_\pi\cL, \cK).
 \end{align*}
 This proves the claim.
 \end{proof}

 \subsection{Cuspidal blocks of \texorpdfstring{$\Z$}{Z}-graded Lie algebras}
 In this subsection, let $\ul L$ be a connected reductive group whose Lie algebra $\ul\fl = \Lie \ul L$ is equipped with a $\Z$-grading given by the weight spaces of a cocharacter $\Cc\to \Aut(\ul L)$ and suppose that $\cC$ is a $\Lz$-cuspidal local system on the open dense $\ul L_0$-orbit in $\flo$. Let $\Db_{\ul L_0}(\ul\fl_1)_{\cC}$ denote the block of $\Db_{\ul L_0}(\ul\fl_1)$ containing $\IC(\cC)$. 
 We will make use of the following lemmas:
 \begin{lemm}\label{lem:block-cuspidal}
  The block $\Db_{\Lz}(\flo)_{\cC}$ is triangle-equivalent to the derived category of  $\rmH^{\bullet}_{Z(\ul L)^{\circ}}\dgmod$.
 \end{lemm}
 \begin{proof}
 Let $(e, h, f)$ be an $\fsl_2$-triple such that $e$ is an element of the open dense $\ul L_0$-orbit in $\ul\fl_1$, $h\in \ul\fl_0$ and $f\in \ul\fl_{-1}$. Let $\varphi\in \bX_*(\ul L_0)$ be the cocharacter satisfying $d\varphi(1) = h$. Then, \cite[\S 4.4]{L95} shows that $\ul\fl_n = \li^\varphi_{2n}{\ul\fl}$ for $n\in \Z$. Let $\ul\fq_{n} = \li^{\varphi}_{\ge n}{\ul\fl}$ and let $\ul Q_0\subseteq \ul L$ be the subgroup such that $\Lie \ul Q_0 = \ul\fq_0$. We may view $\ul\fq_*$ as a $1$-spiral of the ungraded Lie algebra $\ul\fl$ with splitting $\ul\fl_*$. \par
 By the definition of cuspidal pairs in \textit{loc. cit.}, there is a cuspidal $\ul L$-equivariant local system $\tilde\cC$ on an $\ul L$-conjugacy class $O\subseteq \ul\fl$ extending $\cC$. Let $\Db_{\ul L}(\ul\fl^{\nil})_{\tilde\cC}$ denote the thick subcategory of $\Db_{\ul L}(\ul\fl^{\nil})$ spanned by $\IC(\tilde\cC)$. The spiral induction 
 \[
     \Ind^{\ul\fl}_{\ul\fq_1}: \Db_{\ul L_0}(\ul\fl_1)\to \Db_{\ul L}(\ul\fl^{\nil})
 \]
yields $\Ind^{\ul\fl}_{\ul\fq_1}\IC(\cC) \cong \IC(\tilde\cC)$. Moreover, the canonical map
\[
    \Hom(\IC(\cC), \IC(\cC)[k]) \to \Hom(\Ind^{\ul\fl}_{\ul\fq_1}\IC(\cC), \Ind^{\ul\fl}_{\ul\fq_1}\IC(\cC)[k]) \cong\Hom(\IC(\tilde\cC), \IC(\tilde\cC)[k])
\]
is an isomorphism for $k\in \Z$ because both sides are canonically isomorphic to $\rmH^{k}_{Z(\ul L)^{\circ}}$  (see \cite[\S 4.7]{L88}). As $\Db_{\ul L_0}(\ul\fl_1)_{\cC}$ is the thick subcategory spanned by $\IC(\cC)$ and $\Db_{\ul L_0}(\ul\fl^{\nil})_{\tilde\cC}$ is spanned by $\IC(\tilde\cC)$, the spiral induction restricts to an equivalence on the blocks:
 \[
     \Ind^{\ul\fl}_{\ul\fq_1}: \Db_{\ul L_0}(\ul\fl_1)_{\cC}\xrightarrow{\sim} \Db_{\ul L}(\ul\fl^{\nil})_{\tilde\cC}.
 \]
In \cite[Proposition 2.4]{RR21}, it is shown that the latter category is equivalent to the derived category of $\rmH^{\bullet}_{Z(\ul L)^{\circ}}$-dg-modules via  $\RHom(\IC(\tilde\cC), \relbar)$. The required equivalence is then given  by the functor $\RHom(\IC(\cC), \relbar) \cong \RHom(\IC(\tilde\cC), \relbar)\circ\Ind^{\ul\fl}_{\ul\fq_1}$. 
 \end{proof}
 \begin{lemm}\label{lem:duality-cuspidal}
 Given $\cL, \cK\in \Db_{\Lz}(\flo)_{\cC}$, there is an isomorphism:
 \[
 \RHom(\cK, \cL)\cong \RHom(\cL, \cK)^{\vee}.
 \]
 where $M^{\vee} = \RHom_{\rmH^{\bullet}_{Z(\ul L)^{\circ}}}(M, \rmH^{\bullet}_{Z(\ul L)^{\circ}})$. 
 \end{lemm}
 \begin{proof}
The derived category of $\rmH^{\bullet}_{Z(\ul L)^{\circ}}\dgmod$ is equivalent to the homotopy category of finitely generated semi-free $\rmH^{\bullet}_{Z(\ul L)^{\circ}}$-dg-modules, i.e.,  $\rmH^{\bullet}_{Z(\ul L)^{\circ}}$-dg-modules which are finitely generated and free as graded $\rmH^{\bullet}_{Z(\ul L)^{\circ}}$-modules. If $M, N$ are semi-free $\rmH^{\bullet}_{Z(\ul L)^{\circ}}$-dg-modules, then there is a trace pairing
\[
\dgHom(M, N)\times \dgHom(N, M) \xrightarrow{\circ} \dgHom(N, N)\xrightarrow{\operatorname{tr}} \rmH^{\bullet}_{Z(\ul L)^{\circ}},
\]
which induces a natural quasi-isomorphism $\dgHom(M, N) \cong \dgHom(\dgHom(N, M), \rmH^{\bullet}_{Z(\ul L)^{\circ}})$. It becomes an isomorphism when we pass to the derived category. This duality is transferred to the category $\Db_{\Lz}(\flo)_{\cC}$ via~\autoref{lem:block-cuspidal}.
 \end{proof}
 \begin{lemm}\label{lem:duality-torsion}
 Suppose that $M\in \rmH^{\bullet}_{Z(\ul L)^\circ}\dgmod$ is cohomologically bounded. Then, there is a natural isomorphism for each $k\in \Z$:
 \[
 \rmH^{k}( M^{\vee}) \cong \rmH^{-\dim Z(\ul L) - k}(M)^*.
 \]
 \end{lemm}
 \begin{proof}
 This is a special case of Grothendieck's local duality theorem, which can be adapted for dg-modules. Set $A = \rmH^{\bullet}_{Z(\ul L)^\circ}$ and $d = \dim Z(\ul L)^\circ$ and pick homogeneous free generators $X_1, ..., X_d \in A$ of degree $2$, so that $A = \C[X_1, \cdots, X_d]$. For each $l \in \N$, consider the Koszul complex $K(X_1^l, \cdots, X_d^l)\in A\dgmod$. There is a natural quasi-isomorphism $K(X_1^l, \cdots, X_d^l)\xrightarrow{\sim} A / \fm_l$, where $\fm_l\subseteq A$ is the homogeneous ideal generated by $\{X_1^l, \cdots, X_d^l\}$. The local cohomology is defined to be the functor
 \[
 \Gamma_0:A\dgMod \to A\dgMod, \quad M\mapsto \varinjlim_{l\in \N}\dgHom_A(K(X_1^l, \cdots, X_d^l), M).
 \]
 There is a natural morphism:
 \[
 \Gamma_0(M)\to \dgHom_A(\dgHom_A(M, A), \Gamma_0(A))
 \]
 which is a quasi-isomorphism for $M = A$. Note that the functor $\dgHom_A(\relbar, A)$ preserves finitely generated semi-free dg-modules.
 Passing to the derived category, we obtain a natural morphism for $M\in A\dgmod$
 \[
 \Gamma_0(M)\to \RHom_A(M^{\vee}, \Gamma_0(A)).
 \]
 It is an isomorphism because $A$ generates the derived category of $A\dgmod$ as thick subcategory. Finally, one verifies easily that there is a quasi-isomorphism $\Gamma_0(A) \to A^*[d]$, where $A^*$ is the graded $\C$-linear dual of the regular right $A$-module, and hence
 \[
 \RHom_A(M^{\vee}, \Gamma_0(A))\cong
 \RHom_A(M^{\vee}, A^*[d]) \cong \RHom_\C(M^{\vee}, \C)[d].
 \]
 The lemma follows because the natural morphism $\Gamma_0(M) \to M$ is a quasi-isomorphism when $M\in A\dgmod$ is cohomologically bounded.
 \end{proof}
 
 \subsection{Proof of~\autoref{prop:hom-psi-2}}

 In this section, we regard $\cP_{\chi}$ as arising from the construction of \autoref{ssec:double-cover} via Beilinson's construction. 
 
 In what follows, we will omit the prime from the notations $\cZ'$, $\pi'$ etc. Set $D = \dim\fg_{  1} - \dim\fg_0 - \dim\flo + \dim\flz$ and $d = \dim Z(\ul L)$. We have the following identity (see the last paragraph of the proof of~\autoref{lem:gdim}):
\beqn
    \dim\ul\fl_0-\dim\ul\fl_1=d,
\eeqn
so that $N=\dim \fg_1 - \dim \fg_0 = D - d$. \par
 We treat first the case where $x$ is nilpotent. In this case, we have $\cP_{\chi} = \IC(\chi)$. By \eqref{eqn-adjunctions} and \autoref{prop:hyperbolic-spiral}, we have
 \[
 \begin{gathered}
\Hom\left(\IC(\chi), \Ind^{\fg_{  1}}_{\ul\fp_1}\IC(\cC)[k]\right)  
\cong \Hom\left(\Res^{\fg_{  1}}_{\ul\fp_1}\IC(\chi), \IC(\cC)[k]\right) \\
\cong \Hom\left(\res^{\fg_{  1}}_{\ul\fp'_1}\IC(\chi), \IC(\cC)[k-D]\right).
\end{gathered}
 \]
 The duality \autoref{lem:duality-cuspidal} yields
 \[
 \RHom\left(\res^{\fg_{  1}}_{\ul\fp'_1}\IC(\chi), \IC(\cC)[k]\right) \cong  \RHom\left(\IC(\cC), \res^{\fg_{  1}}_{\ul\fp'_1}\IC(\chi)[-k]\right)^{\vee}.
 \]
 The $G_0$-distinguishedness of $x$ and the hypothesis $Z(G)_0 = 1$ implies that $Z_{G_0}(x)^{\circ}$ is unipotent; therefore, the cleanness of $\chi$ implies that the term $\Hom\left(\IC(\chi), \Ind^{\fg_{  1}}_{\ul\fp_1}\IC(\cC)[k]\right)$ vanishes for all but finitely many $k$. Thus  \autoref{lem:duality-torsion} is applicable, so that
 \[
 \Hom\left(\res^{\fg_{  1}}_{\ul\fp'_1}\IC(\chi), \IC(\cC)[k-D]\right) \cong  \Hom\left(\Ind^{\fg_{  1}}_{\ul\fp'_1}\IC(\cC), \IC(\chi)[D-d-k]\right)^*
 \]
 and the statement follows. \par
We assume now that $x$ is non-nilpotent. Then the projection $h:\cZ\to \fg_1$ is proper and $h^{-1}(\fg^{\nil}_1) = \pi^{-1}(0)$. In particular, if $\iota: \fg^{\nil}_1\to \fg_1$ and $i: \fg^{\nil}_1\to \cZ$ denote the inclusions, then we have $\iota^*h_* \cong i^*$. By \autoref{lem:monodromy-unipotent}, we have $\cP_\chi = \psi^{\mathrm{uni}}_{\pi}\cL_{\chi}$. Applying the adjunction~\eqref{eqn-adjunctions}, \autoref{prop:hyperbolic-spiral} and~\autoref{lem:psi-unip}, we obtain
\begin{equation}\label{eq:hom1}
 \begin{gathered}
\Hom\left(\cP_\chi, \Ind^{\fg_{  1}}_{\ul\fp_1}\IC(\cC)[k]\right) \cong \varinjlim_n\Hom\left(h_*j_*(\cL_{\chi}\otimes\cI^{(n)}_{\pi}), \Ind^{\fg_{  1}}_{\ul\fp_1}\IC(\cC)[k]\right) \\
\cong \varinjlim_n\Hom\left(\Res^{\fg_{  1}}_{\ul\fp_1}h_*j_*(\cL_{\chi}\otimes\cI^{(n)}_{\pi}), \IC(\cC)[k]\right) \\
\cong \varinjlim_n\Hom\left(\res^{\fg_{  1}}_{\ul\fp'_1}h_*j_*(\cL_{\chi}\otimes\cI^{(n)}_{\pi}), \IC(\cC)[k-D]\right).
  \end{gathered}
 \end{equation}
 Making use of  \autoref{lem:duality-cuspidal} and adjunction, we obtain
    \begin{equation}\label{eq:hom2}
        \begin{gathered}
            \RHom\left(\res^{\fg_{  1}}_{\ul\fp'_1}h_*j_*(\cL_{\chi}\otimes\cI^{(n)}_{\pi}), \IC(\cC)\right) \cong
            \RHom\left(\IC(\cC), \res^{\fg_{  1}}_{\ul\fp'_1}h_*j_*(\cL_{\chi}\otimes\cI^{(n)}_{\pi})\right)^{\vee} \\
            \cong \RHom\left(\Ind^{\fg_{  1}}_{\ul\fp'_1}\IC(\cC), h_*j_*(\cL_{\chi}\otimes\cI^{(n)}_{\pi})\right)^{\vee} \,.
        \end{gathered}
    \end{equation}
    The spirally induced complex $\Ind^{\fg_{  1}}_{\ul\fp'_1}\IC(\cC)\in \Db_{G_{  0}}(\fg_{  1})$ has a $(G_0\times \Cc)$-equivariant enhancement $\cK\in \Db_{G_{  0}\times \Cc}(\fg_{  1})$.
 
 Let $\cZ_1 = \pi^{-1}(1)$ and let $q: \cU\to \cU / \Cc\cong \cZ_1$ denote the projection; then we have canonical isomorphisms $j^*h^*\cK \cong q^*(h^*\cK|_{\cZ_1})$ and $\cL_{\chi} \cong q^*\IC(\chi)[1]$, which yields for $k\in \Z$:
    \begin{equation}\label{eq:hom3}
        \bega
            \Hom(j^*h^*\cK,  \cL_\chi\otimes \cI^{(n)}_{\pi}[k]) \cong \Hom(q^*(h^*\cK|_{\cZ_1}),  \cL_\chi\otimes \cI^{(n)}_{\pi}[k])\\
            \cong \Hom(h^*\cK|_{\cZ_1},  q_*(q^*\IC(\chi)\otimes \cI^{(n)}_{\pi})[k+1]) \cong \Hom(\cK, \IC(\chi)\otimes \mathrm{R}\Gamma(\Cc, \cI^{(n)}_{\pi})[k+1]).
        \eega
    \end{equation}
    Since $\chi$ is clean, we have for $k\in \Z$:
    \[
        \Hom(\cK, \IC(\chi)[k]) \cong \Hom_{\Db_{Z_{G_0}(x)}(\{x\})}(\cK_x, \chi_x[k + \dim \rc]) \cong \rmH^{k + \dim \rc}_{Z_{G_0}(x)^{\circ}}((\cK_x)^*\otimes\chi_x)^{Z_{G_0}(x) / Z_{G_0}(x)^{\circ}}.
    \]
    The $G_0$-distinguishedness of $x$ and the hypothesis $Z(G)_0 = 1$ implies that $Z_{G_0}(x)^{\circ}$ is unipotent. It follows from the cleanness of $\chi$ that the last term of \eqref{eq:hom3} vanishes for all but finitely many $k\in \Z$.  Now 
    \autoref{lem:duality-torsion} implies that for every $k\in \Z$:
    \[
          \rmH^k(\RHom\left(\cK, h_*j_*(\cL_{\chi}\otimes\cI^{(n)}_{\pi})\right)^{\vee})
           \cong \Hom\left(\cK, h_*j_*(\cL_{\chi}\otimes\cI^{(n)}_{\pi})[-d - k]\right)^*,
    \]
    where $d = \dim Z(\ul L)$. Combining \eqref{eq:hom1}, \eqref{eq:hom2}, \eqref{eq:hom3} and the last line, we have
    \beq\label{eq:hom4}
         \Hom\left(\cP_\chi, \Ind^{\fg_{  1}}_{\ul\fp_1}\IC(\cC)[k]\right) \cong \varinjlim_n
         \Hom\left(\cK, \IC(\chi)\otimes \mathrm{R}\Gamma(\Cc, \cI^{(n)})[D-d-k+1]\right)^*\,.
    \eeq
    To calculate the last term, we have $\mathrm{R}\Gamma(\Cc, \cI^{(n)}) \cong \C \oplus \C[-1]$ and the transition map $\rmH^i(\Cc, \cI^{(n)}) \to \rmH^i(\Cc, \cI^{(n-1)})$ is an isomorphism for $i = 1$ and is $0$ for $i = 0$; hence 
    \[
    \varprojlim_{n\in \N}\rmH^1(\Cc, \cI^{(n)}) \cong \C,\quad  \varprojlim_{n\in \N}\rmH^i(\Cc, \cI^{(n)}) = 0\quad \forall i\neq 1.
    \]
    It follows that
    \beqn
        \varinjlim_n \Hom\left(\cK, \IC(\chi)\otimes \mathrm{R}\Gamma(\Cc, \cI^{(n)})[D-d-k+1]\right)^*  \cong \Hom\left(\cK, \IC(\chi)[D-d-k]\right)^*\,.
    \eeqn
    Combining~\eqref{eq:hom4} with the last line and recalling that $N= D - d$, we conclude the proof. \qed

\section{Applications to homogeneous affine Springer fibres and double affine Hecke algebras}\label{sec:applications}
In this section, we assume $G$ to be simply connected and almost simple and let $\theta: G\to G$ be an automorphism of order $m > 0$. 
\subsection{Graded Lie algebras and loop Lie algebras}
We may choose a pinning $E= (B_0, T_0, a)$ for $G$, where $a: \Ga \to U_0/[U_0,U_0]$ and $U_0 = [B_0,B_0]$, such that $\theta = \sigma\Ad_t$ with $\sigma\in \Aut_E(G)$ and $t\in T := T_0^{\sigma}$. Set $e = \ord \sigma$. Let $\tau$ be the $\zeta_e$-loop rotation on the loop group $L^{[e]}G = G(\!(\varpi^{1/e})\!)$: $\tau(f)(\varpi^{1/e}) =f(\zeta_e\varpi^{1/e})$. The twisted-loop group and Lie algebra are defined as
\[
L^{\sigma}G = (L^{[e]}G)^{\tau^{-1}\sigma},\quad L^{\sigma}\fg = (\fg\otimes \C(\!(\varpi^{1/e})\!))^{\tau^{-1}\sigma}.
\]
Let $\Gm^{\on{rot}} = \Cc$ be a rank-one torus acting on $L^{\sigma}G$ by loop rotation:
\[
    \Gm^{\on{rot}}\times L^{\sigma}G\to L^{\sigma}G,\quad (t, g(\varpi^{1/e})) = g(t\varpi^{1/e}). 
\]
Since $(\Ad_t)^m = \id_G$, we may find $\lambda_0: \Gm\to T$ such that $t = \lambda_0(\zeta_{mc})$, where $c = \#Z(G)$. We put
\[
\lambda = (\lambda_0, mc/e): \Gm\to T\times \Gm^{\on{rot}}. 
\]
Then, $\lambda$ yields $\Gm$-actions on $L^{\sigma}G, L^{\sigma}\fg$ and related spaces, and the evaluation $\varpi\mapsto 1$ induces isomorphisms
\begin{equation}\label{eq:loop-graded}
(L^{\sigma}G)^{\lambda}\xrightarrow{\sim} G_0 = G^{\theta},\quad \li^{\lambda}_i(L^{\sigma}\fg) \xrightarrow{\sim} \fg_{i + m\Z} \quad \text{for $i\in \Z$}. 
\end{equation}
Let $P\subset L^{\sigma}G$ be a $(T\times \Gm^{\on{rot}})$-stable parahoric subgroup (where $T$ acts by conjugacy on $L^{\sigma}G$). It admits a unique $(T\times \Gm^{\on{rot}})$-stable Levi factor, denoted by $L$. 
\begin{lemm}\label{lem:parahoric-spiral}
Under the isomorphism~\eqref{eq:loop-graded}, the $\lambda$-homogeneous components $\{\li^{\lambda}_i\fp\}_{i\in \Z}$ of the Lie algebra $\fp = \Lie P$ is identified with a $(-1)$-spiral $\fp_*$ of $\fg_*$ and $\{\li^{\lambda}_i\fl\}_{i\in \Z}$ is identified with a splitting of $\fp_*$ (see~\autoref{sec:LY}). 
\end{lemm}
\subsection{Homogeneous affine Springer resolutions}
The pinning $E$ determines a standard Iwahori subgroup $B\subset L^{\sigma}G$. Given a triple $\xi = (P, \cO, \cC)$ formed by a standard parahoric subgroup $B\subseteq P\subset L^{\sigma}G$ and a cuspidal pair $(\cO, \cC)$ (in the sense of~\cite{L84}) on the Lie algebra $\fl$ of the $(T\times\Gm^{\on{rot}})$-stable Levi factor $L$, we let $q:\fp\to \fl$ denote projection and form the affine analogue of the Springer resolution by
\[
\widetilde{L^{\sigma}\fg} := \{(x, gP)\in L^{\sigma}\fg\times (L^\sigma G/P)\mid g^{-1}x\in q^{-1}(\cO)\}\xrightarrow{\pi} L^{\sigma}\fg,\quad (x, gP)\mapsto x.
\] 
The cuspidal local system $\cC$ induces a local system on $\widetilde{L^{\sigma}\fg}$, denoted by $\tilde\cC$. Given any regular semisimple topologically nilpotent element $\gamma\in L^{\sigma}\fg$, the space $\on{Sp}_{\gamma} = \pi^{-1}(\gamma)$ is called the generalised affine Springer fibre at $\gamma$. \cite[\S 3.1]{KL88} shows that $\on{Sp}_{\gamma}$ is an ind-scheme locally of finite type. The cohomology $\rmH^*_c(\on{Sp}_{\gamma}, \tilde\cC)$ is closely related to orbital integrals on $p$-adic Lie algebras. \par

The homogeneous generalised affine Springer resolution of degree $i\in\Z$ is given by the homogeneous component of $\pi$ under the action of $\lambda$:
\[
\li^{\lambda}_i(\widetilde{L^{\sigma}\fg}) = \{(x, gP)\in \li^{\lambda}_i(L^{\sigma}\fg)\times (L^\sigma G/P)^{\lambda}\mid g^{-1}x\in q^{-1}(\cO)\}\xrightarrow{\pi^\lambda_i} \li^{\lambda}_i(L^{\sigma}\fg),\quad (x, gP)\mapsto x.
\]
By~\autoref{lem:parahoric-spiral}, setting $\fl_i = \li^{\lambda}_i \fl$, we may identify $\{\fl_i\}_{i\in \Z}$ with a splitting of a spiral of $\fg_*$. \par
Fix $i\in \Z$, we set $\cI_{\xi} := (\li^{\lambda}_i\pi)_! \tilde\cC|_{\li^{\lambda}_i(\widetilde{L^{\sigma}\fg})}$. 
Given $\gamma\in \li^{\lambda}_i(L^{\sigma}\fg)$, the homogeneous affine Springer fibre at $\gamma$ is defined to be $$\on{Sp}_{\gamma}^{\lambda} = (\li^{\lambda}_i\pi)^{-1}(\gamma).$$ It can be identified with the $\lambda$-fixed points of $\on{Sp}_{\gamma}$. The proper base change theorem implies that
\[
    (\cI_{\xi})_{\gamma} \cong \rmH^*_c(\on{Sp}_{\gamma}^{\lambda}, \tilde\cC).
\]
\begin{lemm}
The complex $\cI_{\xi}$ can be identified with the following (infinite) direct sum of cohomological shifts of spiral inductions (see~\autoref{sec:LY}):
\[
\cI_{\xi} \cong \bigoplus_{\ul{\fp}'_*} \Ind_{\ul{\fp}'_i}^{\fg_i} \IC(\cC)[d],
\]
where $\ul{\fp}'_*$ runs over all $G_0$-conjugacy classes of $(-1)$-spirals of $\fg_*$ which contain $\fl_*$ as splitting factor, and $d\in \Z$ depends on $\ul{\fp}'_*$. 
\end{lemm}
Therefore, $\cI_\xi$ is an infinite direct sum of semisimple anti-orbital complexes on $\fg_i$ when $i > 0$. The simple perverse constituents of $\cI_{\xi}$ generate an LY-block $\Db_{G_0}(\fg_i)_\xi\subseteq \Db_{G_0}(\fg_i)^{\nil}$ of the category of anti-orbital complexes on $\fg_i$ (see~\cite{LY17a}).  
\subsection{Character sheaves and double affine Hecke algebras}
In what follows, we will consider for simplicity the case $i=1$ (the discussion remains valid as long as $i > 0$). 
Via the geometric extension algebra $\End(\cI_{\xi})$, we obtain in~\cite{Liu23} an action of a degenerate double affine Hecke algebra $\cH_{\xi}$ on $\cI_{\xi}$, which induces a bijection
\begin{equation}\label{eq:springer}
    \Irr\on{O}_{\lambda}(\mathcal{H}_{\xi}) \xrightarrow{\sim} \Char_{G_{0}}(\fg_1)_{\xi}.
\end{equation}
We recall briefly the definitions of the DAHA $\cH_{\xi}$ and the category $\on{O}_{\lambda}(\mathcal{H}_{\xi})$.

Set $\ft_{\R} = \bX_*(T)\otimes \R$ and $\tilde\ft_{\R} = \bX_*(T\times \Gm^{\on{rot}})\otimes \R$. Let $\delta: \Gm^{\on{rot}}\to \Cc$ be $e$ times the defining character of $\Gm^{\on{rot}}$ and we regard $\delta\in \tilde\ft_{\R}^*$ so that $\delta(\ft_{\R}) = 0$. Let $\Phi = \Phi(L^{\sigma}\fg, \tilde\ft)\subset \tilde\ft_{\R}^*$ be the set of (real) affine roots of $L^{\sigma}\fg$. The Iwahori subgroup $B\subset L^{\sigma}G$ corresponds to a fundamental alcove $\kappa\subset \tilde\ft_{\R}|_{\delta=1}$. The standard parahoric subgroup $P\subset L^{\sigma}G$ corresponds to a facet $\nu\subset \bar{\kappa}$. Let $\tilde\ft_{\xi,\R}\subseteq \tilde\ft_{\R}$ be the linear span of $\nu$ and $\tilde\ft_{\R}|_{\delta=1}$ the affine span of $\nu$. Then, $\Phi$ restricts to a possibly reduced affine root system $\Phi_{\xi}\subset\tilde\ft_{\xi,\R}^*$ on $\tilde\ft_{\xi,\R}|_{\delta=1}$ and $P$ corresponds to a base $\Pi_{\xi}\subset \Phi_{\xi}$. The orthogonal reflections of $\tilde\ft_{\xi,\R}$ with respect to the root hyperplanes of $\Phi_{\xi}$ generate the relative affine Weyl group $W_{\xi}$.  \par
Set $\tilde\ft_{\xi} = \ft_{\xi,\R}\otimes_{\R}\C$. The DAHA $\cH_{\xi}$ is an associative $\C$-algebra with underlying vector space $\C W_{\xi}\otimes \C[\tilde\ft_{\xi}]/(\delta - 1)$. It depends on a set of parameters indexed by the $W_{\xi}$-conjugacy classes in $\Phi_{\xi}$ (see~\cite[\S 2.5]{LY18} and~\cite[\S 5.1]{Liu23}).

The cocharacter $\lambda$ can be identified with a point in $\tilde\ft_{\R}|_{\delta=1}$, which under orthogonal projection yields a point $\lambda_{\xi}\in \tilde\ft_{\xi, \R}|_{\delta=1}$. The category $\on{O}_{\lambda}(\mathcal{H}_{\xi})$ is the category of finitely generated $\mathcal{H}_{\xi}$-modules on which the subalgebra $\C[\tilde\ft_{\xi}]/(\delta - 1)$ act locally finitely with eigenvalues lying in the orbit $W_{\xi}\lambda_{\xi}\subseteq \tilde\ft_{\xi}|_{\delta=1}$.

The affine Weyl group $W_{\xi}$ acts simply transitively on the set of alcoves of $(\tilde\ft_{\xi,\R}|_{\delta=1}, \Phi_{\xi})$. Each $w\in W_{\xi}$ corresponds to a $(T\times\Gm)$-stable parahoric subgroup $wPw^{-1}\subseteq L^{\sigma}G$ which has $L$ as Levi factor. We will denote by $\{w\fp_i\}_{i\in \Z}$ the $(-1)$-spiral of $\fg_*$ corresponding to the parahoric subgroup $wPw^{-1}$ under the correspondence~\autoref{lem:parahoric-spiral}. Given $E\in \on{O}_{\lambda}(\mathcal{H}_{\xi})$, we let $E_{\lambda'}$ denote the generalised $\lambda'$-eigenspace in $E$ for $\lambda'\in W_{\xi}\lambda_{\xi}$. Then, given any $w\in W_{\xi}$, and any character sheaf $\cF\in \Char_{G_0}(\fg_1)_{\xi}$, we have
\begin{equation}\label{eq:dimE}
    \dim E_{w^{-1}\lambda_{\xi}} = [\Ind^{\fg_1}_{\li^{\lambda}_1(w\fp)}\IC(\cC): \cF]
\end{equation}
if $E$ corresponds to $\cF$ under~\eqref{eq:springer}, see~\cite[\S 9.9]{Liu23}.  

\subsection{Cuspidality and finite-dimensional modules}
 
\begin{prop}\label{prop:fd-cuspidal}
	The bijection~\eqref{eq:springer} restricts to a one-to-one correspondence between finite-dimensional simple modules in $\on{O}_{\lambda}(\mathcal{H}_{\xi})$ and cuspidal character sheaves in $\Char_{G_{0}}(\fg_1)_{\xi}$. 
\end{prop}
\begin{proof}
The proof is similar to~\cite[Proposition 6.3.1]{Liu24}.
Let $E\in \on{O}_{\lambda}(\cH_{\xi})$ be a simple module which corresponds to a character sheaf $\cF\in \Char_{G_0}(\fg_1)_{\xi}$ under~\eqref{eq:springer}.
Recall from~\cite[\S 3.4]{VV09} and~\cite[\S5.5.4]{OY16} that the intertwiners of $\cH_{\xi}$ yield a decomposition $W_\xi\lambda_\xi = \bigsqcup_{j\in I}\Lambda_j$ for some finite index set $I$ (called \emph{clans}), such that there is a functorial linear isomorphism $E_{\lambda'} \cong E_{\lambda''}$ whenever $E\in \on{O}_{\lambda}(\cH_{\xi})$ and $\{\lambda',\lambda''\}\subseteq \Lambda_j$ for some $j\in I$. 

It follows that $E$ is infinite dimensional if and only if $E_{\lambda'} \neq 0$ for some $j\in I$ and $\lambda'\in \Lambda_j$ such that $\#\Lambda_j = \infty$. Given such $j\in I$, \eqref{eq:dimE} implies that $[\Ind^{\fg_1}_{w\fp_1}\IC(\tilde\cC):\cF]\neq 0$ for $w\in W_{\xi}$ such that $w^{-1}\lambda_{\xi}\in \Lambda_j$. 
We can find $0\neq \mu\in \bX_*(Z(L)^{\circ})$ such that $\Ind^{\fg_1}_{w\fp_1}$ factorises through the parabolic induction $\Ind^{\fg_1}_{\fq_1}$ for the parabolic subalgebra $\fq = \li^{\mu}_{\ge 0}\fg$ whenever $w^{-1}\lambda_{\xi}\in \Lambda_j$. This implies that $\cF$ is non-cuspidal. 

Conversely, if $\cF$ is non-cuspidal, we can find a proper $(T, \theta)$-stable parabolic subgroup $Q\subseteq G$  such that $\cF$ is a constituent of the parabolic induction $\Ind^{\fg_1}_{\fq_1}\cG$ for some $\cG\in \Char_{M_0}(\fm_1)$, where $M$ is the unique $(T, \theta)$-stable Levi component of $Q$. Since $\cF$ lies in the $\xi$-block, up to replacing $(Q, M, \cG)$ with a $G_0$-conjugate, we may assume that $\fl_1\subseteq \fm_{i+m\Z}$ for $i\in \Z$ and $\cG\in \Perv_{M_0}(\fm_1)_{\xi}$ ($\xi$ regarded as a LY cuspidal support on $\fm_1$). Since $Q$ is $\theta$-stable, we may choose $\mu\in \bX_*(Z(M)^{\sigma})\subset\ft_{\xi,\R}$ such that $\fq = \li^{\mu}_{\ge 0}\fg$. If $\bar P\subseteq L^{\sigma}M$ is a $(T\times \Gm^{\on{rot}})$-stable parahoric subgroup such that $\cG\subseteq \Ind^{\fm_1}_{\bar\fp_1}\IC(\tilde\cC)$, then $\cF\subseteq\Ind^{\fg_1}_{\fq_1}\Ind^{\fm_1}_{\bar\fp_1}\IC(\tilde\cC)$. If $w\in W_{\xi}$ satisfies $\Ind^{\fg_1}_{w\fp_1}\cong\Ind^{\fg_1}_{ \fq_1}\Ind^{\fm_1}_{\bar\fp_1}$, then $\Ind^{\fg_1}_{(wt_{n\mu}\fp)_1}\cong\Ind^{\fg_1}_{\fq_1}\Ind^{\fm_1}_{\bar\fp_1}$ for any $n \ge 0$. Therefore, $\cF\subseteq\Ind^{\fg_1}_{(wt_{n\mu})\fp_1}\IC(\tilde\cC)$ for $n\ge 0$. It follows from~\eqref{eq:dimE} that $E$ is infinite dimensional.  
\end{proof}

\subsection{Elliptic Springer fibres and finite-dimensional modules}
A regular semisimple element $\gamma \in L^{\sigma}\fg$ is called \emph{elliptic} if the centraliser $T_{\gamma} = Z_{L^{\sigma}G}(\gamma)$ is finite. Under the identification~\eqref{eq:loop-graded}, the existence of a homogeneous elliptic regular semisimple element $\gamma\in L^{\sigma}\fg$ is equivalent to the GIT-stability of the action $G_0\curvearrowright\fg_1$, which can occur only when $m$ is a (twisted) elliptic regular number of the pair $(W(G, T_0), \sigma)$ in the sense of Springer~\cite{S74}. The following result is an extension of~\cite[Theorem 3.3.1]{VV09} to non-principal series and non-spherical modules.
\begin{theo}\label{theo:DAHA-fd}
Suppose there is a homogeneous elliptic regular semisimple topologically nilpotent element $\gamma\in \li^{\lambda}_1(L^{\sigma}\fg)$. Then, the cohomology of homogeneous affine Springer fibre $\rmH^*_c(\on{Sp}_{\gamma}^{\lambda}, \tilde\cC)$ is a semisimple finite-dimensional $\cH_{\xi}$-module lying in $\cO_{\lambda}(\cH_{\xi})$. Moreover, every finite-dimensional simple module in $\cO_{\lambda}(\cH_{\xi})$ appears as a constituent of $\rmH^*_c(\on{Sp}_{\gamma}^{\lambda}, \tilde\cC)$.
\end{theo}
 
\begin{proof}
 
The standard geometric Ext-algebra formalism in~\cite[Ch. 8]{CG10} endows the stalk $(\cI_{\xi})_\gamma = \rmH^*_c(\on{Sp}_{\gamma}^{\lambda}, \tilde\cC)$ with an $\cH_{\xi}$-module structure which lies in $\cO_{\lambda}(\cH_{\xi})$, as shown in~\cite[\S 7.6]{Va05} and~\cite[\S 9.14]{Liu23}. \par
The decomposition theorem implies the semisimplicity of $\cI_{\xi}$. Therefore, we can write 
\[
\cI_{\xi} = \bigoplus_{\cG\in \Char_{G_0}(\fg_1)_{\xi}}E_{\cG}\otimes\cG,
\]
where $E_{\cG}$ is a graded vector space, and the construction of~\cite{Liu23} gives $E_{\cG}$ the structure of a simple $\cH_{\xi}$-module in $\on{O}_{\lambda}(\cH_{\xi})$. 

The elliptic regular semisimple elements form a dense open subset $U\subset (L^{\sigma}\fg)^{\lambda}_1$ and $\gamma\in U$. Under the identification~\eqref{eq:loop-graded}, $U$ is a stratum of $\fg_1$ in the sense of~\autoref{stratification}. We can write
\[
    \rmH^*_c(\on{Sp}_{\gamma}^{\lambda}, \tilde\cC) = (\cI_{\xi})_{\gamma} = \bigoplus_{\substack{\cG\in \Char_{G_0}(\fg_1)_{\xi} \\ \on{supp}\cG = \fg_1}}E_{\cG}\otimes\cG_{\gamma},
\]
and $\cG_{\gamma}$ is concentrated in degree $-\dim\fg_1$. \autoref{coro:GIT-stable} implies that the character sheaves $\cG$ in the last direct sum are exactly the cuspidal character sheaves in $\Char_{G_0}(\fg_1)_{\xi}$. By~\autoref{prop:fd-cuspidal}, the simple $\cH_{\xi}$-module $E_{\cG}\in\cO_{\lambda}(\cH_{\xi})$ parametrised by a character sheaf $\cG\in \Char_{G_0}(\fg_1)_{\xi}$ is finite dimensional if and only if $\cG$ is of full support. Moreover, every such $E_{\cG}$ appears in $\rmH^*_c(\on{Sp}_{\gamma}^{\lambda}, \tilde\cC)$ with multiplicity space $\cG_{\gamma}$ (see also~\cite[\S 3.3.1]{VV09}). Since the number of isomorphism classes character sheaves on $\fg_1$ is finite, the cohomology $\rmH^*_c(\on{Sp}_{\gamma}^{\lambda}, \tilde\cC)$ is finite dimensional.
\end{proof}

\printbibliography

@incollection {B87,
    AUTHOR = {Beilinson, Alexander},
     TITLE = {How to glue perverse sheaves},
 BOOKTITLE = {{$K$}-theory, arithmetic and geometry ({M}oscow, 1984--1986)},
    SERIES = {Lecture Notes in Math.},
    VOLUME = {1289},
     PAGES = {42--51},
 PUBLISHER = {Springer, Berlin},
      YEAR = {1987},
      ISBN = {3-540-18571-2},
   MRCLASS = {14F99 (18E25 32C38)},
  MRNUMBER = {923134},
MRREVIEWER = {Jean-Luc\ Brylinski},
       DOI = {10.1007/BFb0078366},
       URL = {https://doi.org/10.1007/BFb0078366},
}

@Book{BBDG18,
 Author = {Beilinson, Alexander and Bernstein, Joseph and Deligne, Pierre and Gabber, Ofer},
 Title = {Faisceaux pervers. {Actes} du colloque ``{Analyse} et {Topologie} sur les {Espaces} {Singuliers}''. {Partie} {I}},
 Edition = {2nd edition},
 FSeries = {Ast{\'e}risque},
 Series = {Ast{\'e}risque},
 ISSN = {0303-1179},
 Volume = {100},
 ISBN = {978-2-85629-878-7},
 Year = {2018},
 Publisher = {Paris: Soci{\'e}t{\'e} Math{\'e}matique de France},
 Language = {French},
 zbMATH = {6868966},
 Zbl = {1390.14055}
}

@Article{B03,
 Author = {Braden, Tom},
 Title = {Hyperbolic localization of intersection cohomology},
 FJournal = {Transformation Groups},
 Journal = {Transform. Groups},
 ISSN = {1083-4362},
 Volume = {8},
 Number = {3},
 Pages = {209--216},
 Year = {2003},
 Language = {English},
 DOI = {10.1007/s00031-003-0606-4},
 zbMATH = {2005214},
 Zbl = {1026.14005}
}

@incollection {Bry86,
    AUTHOR = {Brylinski, Jean-Luc},
     TITLE = {Transformations canoniques, dualit\'{e} projective,
              th\'{e}orie de {L}efschetz, transformations de {F}ourier et
              sommes trigonom\'{e}triques},
      Booktitle = {G\'{e}om\'{e}trie et analyse microlocales},
     FSeries = {Ast{\'e}risque},
     Series = {Ast{\'e}risque},
    VOLUME = {140-141},
      YEAR = {1986},
     PAGES = {3--134},
      ISSN = {0303-1179,2492-5926},
   MRCLASS = {32C38 (14F05 32C42 58G07)},
  MRNUMBER = {864073},
}

@book{C93,
    AUTHOR = {Carter, Roger W.},
     TITLE = {Finite groups of {L}ie type},
    SERIES = {Wiley Classics Library},
      NOTE = {Conjugacy classes and complex characters,
              Reprint of the 1985 original,
              A Wiley-Interscience Publication},
 PUBLISHER = {John Wiley \& Sons, Ltd., Chichester},
      YEAR = {1993},
     PAGES = {xii+544},
      ISBN = {0-471-94109-3},
   MRCLASS = {20C33 (20-02 20G40)},
  MRNUMBER = {1266626},
}

@book{CG10,
 author = {Chriss, Neil and Ginzburg, Victor},
 title = {Representation theory and complex geometry},
 edition = {Reprint of the 1997 original},
 isbn = {978-0-8176-4937-1},
 year = {2010},
 publisher = {Boston, MA: Birkh{\"a}user},
 language = {English},
 doi = {10.1007/978-0-8176-4938-8},
 zbMATH = {5622606},
 Zbl = {1185.22001}
}

@Article{Grin98,
 Author = {Grinberg, Mikhail},
 Title = {A generalization of {Springer} theory using nearby cycles},
 FJournal = {Representation Theory},
 Journal = {Represent. Theory},
 ISSN = {1088-4165},
 Volume = {2},
 Pages = {410--431},
 Year = {1998},
 Language = {English},
 zbMATH = {1267661},
 Zbl = {0938.22011}
}

@misc{GVX20,
      title={Nearby Cycle Sheaves for Stable Polar Representations}, 
      author={Grinberg, Mikhail and Vilonen, Kari and Xue, Ting},
      year={2024},
      eprint={2012.14522v2},
      archivePrefix={arXiv},
      primaryClass={math.AG}
}

@Article{GVX23,
 Author = {Grinberg, Mikhail and Vilonen, Kari and Xue, Ting},
 Title = {Nearby cycle sheaves for symmetric pairs},
 FJournal = {American Journal of Mathematics},
 Journal = {Am. J. Math.},
 ISSN = {0002-9327},
 Volume = {145},
 Number = {1},
 Pages = {1--63},
 Year = {2023},
 Language = {English},
 DOI = {10.1353/ajm.2023.0000},
 zbMATH = {7653690},
 Zbl = {1528.22011}
}

@article{KL88,
 author = {Kazhdan, David and Lusztig, George},
 title = {Fixed point varieties on affine flag manifolds},
 fjournal = {Israel Journal of Mathematics},
 journal = {Isr. J. Math.},
 issn = {0021-2172},
 volume = {62},
 number = {2-3},
 pages = {129--168},
 year = {1988},
 language = {English},
 doi = {10.1007/BF02787119},
 zbMATH = {4075419},
 Zbl = {0658.22005}
}

@Book{KS90,
 Author = {Kashiwara, Masaki and Schapira, Pierre},
 Title = {Sheaves on manifolds. {With} a short history ``{Les} d{\'e}buts de la th{\'e}orie des faisceaux'' by {Christian} {Houzel}},
 FSeries = {Grundlehren der Mathematischen Wissenschaften},
 Series = {Grundlehren Math. Wiss.},
 ISSN = {0072-7830},
 Volume = {292},
 ISBN = {3-540-51861-4},
 Year = {1990},
 Publisher = {Springer-Verlag},
 Language = {English},
 zbMATH = {47944},
 Zbl = {0709.18001}
}

@article{Liu23,
 author = {Liu, Wille},
 title = {Generalized {Springer} correspondence for {{\(\mathbf{Z}/m\)}}-graded {Lie} algebras},
 fjournal = {Annales Scientifiques de l'{\'E}cole Normale Sup{\'e}rieure. Quatri{\`e}me S{\'e}rie},
 journal = {Ann. Sci. {\'E}c. Norm. Sup{\'e}r. (4)},
 issn = {0012-9593},
 volume = {56},
 number = {5},
 pages = {1449--1515},
 year = {2023},
 language = {English},
 doi = {10.24033/asens.2559},
 zbMATH = {7827707},
 Zbl = {1542.20033}
}

@misc{Liu24,
      title={Bi-orbital sheaves and affine Hecke algebras at roots of unity}, 
      author={Liu, Wille},
      year={2024},
      eprint={1911.11587v2},
      archivePrefix={arXiv},
      primaryClass={math.RT}
}

@misc{LVX,
      title={Cuspidal character sheaves on graded Lie algebras II}, 
      author={Wille Liu and Kari Vilonen and Ting Xue},
      year={2025},
      eprint={2512.20472},
      archivePrefix={arXiv},
      primaryClass={math.RT},
      url={https://arxiv.org/abs/2512.20472}, 
}

@Article{L84,
 Author = {Lusztig, George},
 Title = {Intersection cohomology complexes on a reductive group},
 FJournal = {Inventiones Mathematicae},
 Journal = {Invent. Math.},
 ISSN = {0020-9910},
 Volume = {75},
 Pages = {205--272},
 Year = {1984},
 Language = {English},
 DOI = {10.1007/BF01388564},
 zbMATH = {3871643},
 Zbl = {0547.20032}
}

@Article{L85a,
 Author = {Lusztig, George},
 Title = {Character sheaves. {I}},
 FJournal = {Advances in Mathematics},
 Journal = {Adv. Math.},
 ISSN = {0001-8708},
 Volume = {56},
 Pages = {193--237},
 Year = {1985},
 Language = {English},
 DOI = {10.1016/0001-8708(85)90034-9},
 zbMATH = {3939580},
 Zbl = {0586.20018}
}

@incollection {L87,
    AUTHOR = {Lusztig, George},
     TITLE = {Fourier transforms on a semisimple {L}ie algebra over {${\bf
              F}_q$}},
 BOOKTITLE = {Algebraic groups {U}trecht 1986},
    SERIES = {Lecture Notes in Math.},
    VOLUME = {1271},
     PAGES = {177--188},
 PUBLISHER = {Springer, Berlin},
      YEAR = {1987},
   MRCLASS = {17B45 (20G40)},
  MRNUMBER = {911139},
MRREVIEWER = {James E. Humphreys},
       DOI = {10.1007/BFb0079237},
       URL = {https://doi.org/10.1007/BFb0079237},
}

@Article{L88,
 Author = {Lusztig, George},
 Title = {Cuspidal local systems and graded {Hecke} algebras. {I}},
 FJournal = {Publications Math{\'e}matiques},
 Journal = {Publ. Math., Inst. Hautes {\'E}tud. Sci.},
 ISSN = {0073-8301},
 Volume = {67},
 Pages = {145--202},
 Year = {1988},
 Language = {English},
 DOI = {10.1007/BF02699129},
 zbMATH = {4146302},
 Zbl = {0699.22026}
}

@article{L95,
	title={Study of perverse sheaves arising from graded Lie algebras},
	author={Lusztig, George},
	journal={Adv. Math.},
	volume={112},
	number={2},
	pages={147--217},
	year={1995},
	publisher={Elsevier}
}

@article{LY17a,
	AUTHOR = {Lusztig, George and Yun, Zhiwei},
	TITLE = {{$\mathbf{Z}/m$}-graded {L}ie algebras and perverse sheaves,
		{I}},
	JOURNAL = {Represent. Theory},
	VOLUME = {21},
	YEAR = {2017},
	PAGES = {277--321},
	ISSN = {1088-4165},
	MRCLASS = {17B70 (14F05 20G99)},
	MRNUMBER = {3697026},
	DOI = {10.1090/ert/500},
	URL = {https://doi.org/10.1090/ert/500},
}

@article{LY18,
	AUTHOR = {Lusztig, George and Yun, Zhiwei},
	TITLE = {{$\mathbf{Z}/m\mathbf{Z}$}-graded {L}ie algebras and perverse
		sheaves, {III}: {G}raded double affine {H}ecke algebra},
	JOURNAL = {Represent. Theory},
	VOLUME = {22},
	YEAR = {2018},
	PAGES = {87--118},
	ISSN = {1088-4165},
	MRCLASS = {20G25 (17B70 20C08)},
	MRNUMBER = {3829497},
	DOI = {10.1090/ert/515},
	URL = {https://doi.org/10.1090/ert/515},
}

@Article{MV87,
 Author = {Mirollo, Renato and Vilonen, Kari},
 Title = {Bernstein-{Gelfand}-{Gelfand} reciprocity on perverse sheaves},
 FJournal = {Annales Scientifiques de l'{\'E}cole Normale Sup{\'e}rieure. Quatri{\`e}me S{\'e}rie},
 Journal = {Ann. Sci. {\'E}c. Norm. Sup{\'e}r. (4)},
 ISSN = {0012-9593},
 Volume = {20},
 Number = {3},
 Pages = {311--323},
 Year = {1987},
 Language = {English},
 DOI = {10.24033/asens.1536},
 zbMATH = {4061384},
 Zbl = {0651.14010}
}

@Article{MV88,
 Author = {Mirkovi{\'c}, Ivan and Vilonen, Kari},
 Title = {Characteristic varieties of character sheaves},
 FJournal = {Inventiones Mathematicae},
 Journal = {Invent. Math.},
 ISSN = {0020-9910},
 Volume = {93},
 Number = {2},
 Pages = {405--418},
 Year = {1988},
 Language = {English},
 DOI = {10.1007/BF01394339},
 zbMATH = {4118629},
 Zbl = {0683.22012}
}

@InCollection{N17,
 Author = {Nakajima, Hiraku},
 Title = {Lectures on perverse sheaves on instanton moduli spaces},
 BookTitle = {Geometry of moduli spaces and representation theory. Lecture notes from the 2015 IAS/Park City Mathematics Institute (PCMI) Graduate Summer School, June 28 -- July 18, 2015},
 Pages = {381--436},
 Year = {2017},
 Publisher = {American Mathematical Society},
 Language = {English},
 zbMATH = {6854853},
 Zbl = {1403.14036}
}

@article{OY16,
 author = {Oblomkov, Alexei and Yun, Zhiwei},
 title = {Geometric representations of graded and rational {Cherednik} algebras},
 fjournal = {Advances in Mathematics},
 journal = {Adv. Math.},
 issn = {0001-8708},
 volume = {292},
 pages = {601--706},
 year = {2016},
 language = {English},
 doi = {10.1016/j.aim.2016.01.015},
 zbMATH = {6548175},
 Zbl = {1403.20007}
}

@Article{R19,
 Author = {Richarz, Timo},
 Title = {Spaces with {{\(\mathbb{G}_m\)}}-action, hyperbolic localization and nearby cycles},
 FJournal = {Journal of Algebraic Geometry},
 Journal = {J. Algebr. Geom.},
 ISSN = {1056-3911},
 Volume = {28},
 Number = {2},
 Pages = {251--289},
 Year = {2019},
 Language = {English},
 DOI = {10.1090/jag/710},
 zbMATH = {7018318},
 Zbl = {1444.14085}
}

@Article{RR21,
 Author = {Rider, Laura and Russell, Amber},
 Title = {Formality and {Lusztig}'s generalized {Springer} correspondence},
 FJournal = {Algebras and Representation Theory},
 Journal = {Algebr. Represent. Theory},
 ISSN = {1386-923X},
 Volume = {24},
 Number = {3},
 Pages = {699--714},
 Year = {2021},
 Language = {English},
 DOI = {10.1007/s10468-020-09966-w},
 zbMATH = {7374841},
 Zbl = {1484.14036}
}

@article{S74,
 author = {Springer, T. A.},
 title = {Regular elements of finite reflection groups},
 fjournal = {Inventiones Mathematicae},
 journal = {Invent. Math.},
 issn = {0020-9910},
 volume = {25},
 pages = {159--198},
 year = {1974},
 language = {English},
 doi = {10.1007/BF01390173},
 url = {https://eudml.org/doc/142286},
 zbMATH = {3450476},
 Zbl = {0287.20043}
}

@Article{S78,
 Author = {Springer, Tonny A.},
 Title = {A construction of representations of {Weyl} groups},
 FJournal = {Inventiones Mathematicae},
 Journal = {Invent. Math.},
 ISSN = {0020-9910},
 Volume = {44},
 Pages = {279--293},
 Year = {1978},
 Language = {English},
 DOI = {10.1007/BF01403165},
 zbMATH = {3585641},
 Zbl = {0376.17002}
}

@Article{V94,
 Author = {Vilonen, Kari},
 Title = {Perverse sheaves and finite dimensional algebras},
 FJournal = {Transactions of the American Mathematical Society},
 Journal = {Trans. Am. Math. Soc.},
 ISSN = {0002-9947},
 Volume = {341},
 Number = {2},
 Pages = {665--676},
 Year = {1994},
 Language = {English},
 DOI = {10.2307/2154577},
 zbMATH = {548948},
 Zbl = {0811.14016}
}

@Article{VX22,
 Author = {Vilonen, Kari and Xue, Ting},
 Title = {Character sheaves for classical symmetric pairs},
 FJournal = {Representation Theory},
 Journal = {Represent. Theory},
 ISSN = {1088-4165},
 Volume = {26},
 Pages = {1097--1144},
 Year = {2022},
 Language = {English},
 DOI = {10.1090/ert/622},
 URL = {hdl.handle.net/10138/355555},
 zbMATH = {7612834},
 Zbl = {1527.20077}
}

@article {VX23sl,
    AUTHOR = {Vilonen, Kari and Xue, Ting},
     TITLE = {Character sheaves for symmetric pairs: special linear groups},
   JOURNAL = {Trans. Amer. Math. Soc.},
  FJOURNAL = {Transactions of the American Mathematical Society},
    VOLUME = {376},
      YEAR = {2023},
    NUMBER = {2},
     PAGES = {837--853},
      ISSN = {0002-9947},
   MRCLASS = {20G20 (14L35 17B08)},
  MRNUMBER = {4531663},
MRREVIEWER = {Yu Chen},
       DOI = {10.1090/tran/8825},
       URL = {https://doi.org/10.1090/tran/8825},
}

@article {VX23st,
    AUTHOR = {Vilonen, Kari and Xue, Ting},
     TITLE = {Character sheaves for graded {L}ie algebras: stable gradings},
   JOURNAL = {Adv. Math.},
  FJOURNAL = {Advances in Mathematics},
    VOLUME = {417},
      YEAR = {2023},
     PAGES = {Paper No. 108935, 59},
      ISSN = {0001-8708},
   MRCLASS = {17B70 (14F08 20C08 22E57)},
  MRNUMBER = {4554668},
MRREVIEWER = {Volodymyr Mazorchuk},
       DOI = {10.1016/j.aim.2023.108935},
       URL = {https://doi.org/10.1016/j.aim.2023.108935},
}

@misc{VX21,
      title={Invariant systems and character sheaves for graded Lie algebras}, 
      author={Kari Vilonen and Ting Xue},
      year={2024},
      eprint={2111.08403},
      archivePrefix={arXiv},
      primaryClass={math.RT}
}

@article{Va05,
 author = {Vasserot, Eric},
 title = {Induced and simple modules of double affine {Hecke} algebras.},
 fjournal = {Duke Mathematical Journal},
 journal = {Duke Math. J.},
 issn = {0012-7094},
 volume = {126},
 number = {2},
 pages = {251--323},
 year = {2005},
 language = {English},
 doi = {10.1215/S0012-7094-04-12623-5},
 zbMATH = {2151092},
 Zbl = {1114.20002}
}

@Article{Vi77,
 Author = {Vinberg, {\`E}rnest Borisovich},
 Title = {The {Weyl} group of a graded {Lie} algebra},
 FJournal = {Mathematics of the USSR. Izvestiya},
 Journal = {Math. USSR, Izv.},
 ISSN = {0025-5726},
 Volume = {10},
 Pages = {463--495},
 Year = {1977},
 Language = {English},
 DOI = {10.1070/IM1976v010n03ABEH001711},
 zbMATH = {3577478},
 Zbl = {0371.20041}
}

@article{VV09,
 author = {Varagnolo, M. and Vasserot, E.},
 title = {Finite-dimensional representations of {DAHA} and affine {Springer} fibers: the spherical case.},
 fjournal = {Duke Mathematical Journal},
 journal = {Duke Math. J.},
 issn = {0012-7094},
 volume = {147},
 number = {3},
 pages = {439--540},
 year = {2009},
 language = {English},
 doi = {10.1215/00127094-2009-016},
 zbMATH = {5550766},
 Zbl = {1237.20008}
}

@Article{Wh65,
 Author = {Whitney, Hassler},
 Title = {Tangents to an analytic variety},
 FJournal = {Annals of Mathematics. Second Series},
 Journal = {Ann. Math. (2)},
 ISSN = {0003-486X},
 Volume = {81},
 Pages = {496--549},
 Year = {1965},
 Language = {English},
 DOI = {10.2307/1970400},
 zbMATH = {3244703},
 Zbl = {0152.27701}
}

@Article{X24,
 Author = {Xue, Ting},
 Title = {Character sheaves for classical graded {Lie} algebras},
 FJournal = {Acta Mathematica Sinica. English Series},
 Journal = {Acta Math. Sin., Engl. Ser.},
 ISSN = {1439-8516},
 Volume = {40},
 Number = {3},
 Pages = {870--884},
 Year = {2024},
 Language = {English},
 DOI = {10.1007/s10114-023-2079-9},
 zbMATH = {7815114}}

@article {X21,
    AUTHOR = {Xue, Ting},
     TITLE = {Character sheaves for symmetric pairs: spin groups},
   JOURNAL = {Pure Appl. Math. Q.},
  FJOURNAL = {Pure and Applied Mathematics Quarterly},
    VOLUME = {21},
      YEAR = {2025},
    NUMBER = {1},
     PAGES = {591--629},
      ISSN = {1558-8599},
   MRCLASS = {20G05 (14L35 17B08 20G15)},
  MRNUMBER = {4847246},
       DOI = {10.4310/pamq.241203044944},
       URL = {https://doi.org/10.4310/pamq.241203044944},
}

@article {Yu01,
    AUTHOR = {Yu, Jiu-Kang},
     TITLE = {Construction of tame supercuspidal representations},
   JOURNAL = {J. Amer. Math. Soc.},
  FJOURNAL = {Journal of the American Mathematical Society},
    VOLUME = {14},
      YEAR = {2001},
    NUMBER = {3},
     PAGES = {579--622 (electronic)},
      ISSN = {0894-0347},
   MRCLASS = {22E50},
  MRNUMBER = {1824988 (2002f:22033)},
MRREVIEWER = {Bertrand Lemaire},
       DOI = {10.1090/S0894-0347-01-00363-0},
       URL = {http://dx.doi.org/10.1090/S0894-0347-01-00363-0},
}

@misc{NY25,
      title={Character sheaves on loop Lie algebras: polar partition}, 
      author={Bao Chau Ngo and Zhiwei Yun},
      year={2025},
      eprint={2506.14584},
      archivePrefix={arXiv},
      primaryClass={math.RT},
      url={https://arxiv.org/abs/2506.14584}, 
note = {\url{https://arxiv.org/abs/2506.14584v1}},
}

\end{document}